\documentclass[a4paper,12pt]{article}
\usepackage{amsthm}
\usepackage{amsmath}
\allowdisplaybreaks[4]%
\usepackage{bm}%
\usepackage{amssymb}
\usepackage{tabularx}
\usepackage{indentfirst}
\usepackage{enumerate}
\usepackage{dsfont}
\usepackage{graphicx}
\usepackage{subfigure}
\usepackage{epsfig}
\usepackage{graphics}
\usepackage{cases}
\usepackage{caption}
\usepackage[compress]{cite}
\usepackage{txfonts}
\usepackage{geometry}
 \usepackage{epstopdf}
 \usepackage{lineno}
 \usepackage[justification=centering]{caption}
 \usepackage{color}
 \usepackage{hyperref}
\hypersetup{
colorlinks=true,
linkcolor=blue,
anchorcolor=blue,
citecolor=blue}
\newtheorem{theorem}{Theorem}[section]
\newtheorem{lemma}{Lemma}[section]

\newtheorem{remark}{Remark}[section]
\newtheorem{definition}{Definition}[section]

\numberwithin{equation}{section}

\begin{document}

\thispagestyle{empty}

\title{On a two-season faecal-oral model with impulsive intervention\thanks{The first author is supported by Postgraduate Research \& Practice Innovation Program of Jiangsu Province (KYCX24\_3711), the second author is supported by the National Natural Science Foundation of China (No. 12271470) , and the third author acknowledges the support of CNPq/Brazil Proc. $N^{o}$ $311562/2020-5$ and FAPDF grant 00193.00001133/2021-80.}}

\date{\empty}

\author{Qi Zhou$^1$, Zhigui Lin$^1 \thanks{Corresponding author. Email: zglin@yzu.edu.cn (Z. Lin).}$ and Carlos Alberto Santos$^{2}$\\
{\small $^1$ School of Mathematical Science, Yangzhou University, Yangzhou 225002, China}\\
{\small $^2$ Department of Mathematics, University of Brasilia, BR-70910900 Brasilia, DF, Brazil}
}
 \maketitle
\begin{quote}
\noindent
{\bf Abstract.}{\footnotesize\small~Rainfall is associated with the outbreak of certain waterborne faecal-oral diseases, driving the implementation of various human interventions for their control and prevention. Taking into account human intervention and temporal variation in rainfall, this paper develops a two-season switching faecal-oral model with impulsive intervention and free boundaries. In this model, the infection fronts are represented by fixed boundaries during the dry season and by moving boundaries during the wet season, with impulsive intervention occurring at the end of each wet season. The simultaneous introduction of impulsive intervention and seasonal switching creates new difficulties for mathematical analysis. We overcome these challenges through novel analytical techniques, resulting in a spreading-vanishing dichotomy and a sharp criteria governing this dichotomy. Finally, numerical simulations are presented to validate the theoretical results and to visually illustrate the influence of seasonal switching and impulsive intervention. Our results mathematically explain that two factors, the duration of the dry season and the intensity of impulsive intervention are both positively correlated with effective disease control.
}


\end{quote}

\section{Introduction}
The faecal-oral route can be used to describe the spread of infectious diseases like cholera, poliomyelitis, typhoid fever, infectious hepatitis, and hand-foot-mouth diseases \cite{Zhou-Lin-Santos-CFS}. Additionally, studies suggest that coronavirus disease 2019 (COVID-19) may also be transmitted via this route \cite{4,5}. The diseases transmitted via the faecal-oral route have posed a significant threat to human survival. Mogasale et al. estimated that typhoid disease results between 9.9 and 14.7 million infection cases and between 75,000 and 208,000 deaths annually \cite{Mogasale-Maskery}. Hay et al. estimated that 15.9 million deaths globally were attributable to the COVID-19 pandemic in 2020 and 2021 combined \cite{Hay-collaborators}.

The emergence of an expanding front is a common feature of natural spreading processes. In 2010, Du and Lin proposed a diffusive logistic model with a free boundary in one space dimension, in which the well-known Stefan condition was employed to characterize the spreading front of a new or invasive species \cite{Du-Lin-SIAM-JMA}.
A biological interpretation of the Stefan condition is provided in \cite{Bunting-Du-Krakowski-NHM}. Subsequently, this condition has been used to
describe the expanding boundary of infection in epidemic models \cite{Hadeler-DCDSB}; we refer to \cite{Kim-Lin-Zhang-NARWA,Li-Ni-Wang-DCDSB}
for susceptible-infective-recovered models, \cite{Ge-Kim-Lin-Zhu-JDE} for susceptible-infective-susceptible models,
\cite{Lin-Zhu-JMB,Wang-Nie-Du-JMB,Tarboush-Lin-Zhang-SCM} for West Nile virus models, \cite{Ahn-Baek-Lin-AMM,Zhao-Li-DCDSB,Wang-Du-DCDSB} for faecal-oral
models, and the references therein.

Notably, all the aforementioned faecal-oral transmission models are formulated within a continuous modeling framework. However, human intervention can render disease progression a non-continuous process. For example, routine environmental disinfection, one of the most effective measures for disease prevention and control, primarily functions by disrupting the most vulnerable link in the faecal-oral transmission chain: the persistence of pathogens in the environment. This intervention measure leads to a rapid decline in pathogenic bacterial load in the environment within a short time frame, thereby significantly influencing disease dynamics. Conventional continuous models fail to capture such abrupt, transient changes. The most prominent feature of impulse differential equation models is their capability to fully account for the impact of transient sudden phenomena on disease progression. To more accurately characterize disease dynamics, an impulsive faecal-oral model within an expanding infected environment has been proposed in \cite{Zhou-Lin-Pedersen-ARXIV}. Moreover, impulsive reaction-diffusion models have also been applied in other contexts; we refer to \cite{Li-Zhao-Cheng-CNSNS,Li-Zhao-JNS} for modeling mosquito control, \cite{Liang-Yan-Xiang-CNSNS} for modeling pest control and management, \cite{Lewis-Li-BMB,Fazly-Lews-Wang-SJAM,Lu-Wang-JMB,Lu-Wang-Zhao-SCM} for describing the population dynamics of species with distinct reproductive and dispersal stages, and \cite{Huang-Wang-Lewis-SJAM} for describing the dynamics of invasive freshwater mussels in rivers.

The primary environmental media involved in faecal-oral transmission include water, food, hands and skin surfaces, fomites, soil, and biological vectors, among which contaminated water is the most significant and typical medium. Seven common waterborne diseases are typhoid fever, cholera, giardia, dysentery, escherichia coli, hepatitis A, and salmonella \cite{Lifewater}. Annually, waterborne diseases afflict hundreds of millions of people, primarily those living without safe, accessible water in developing countries \cite{Lifewater}. A substantial body of research has demonstrated that the transmission of cholera, a common waterborne disease, exhibits distinct seasonal patterns in certain regions; see, for example, \cite{Ruiz-Moreno,Perez-Saez,Shackleton} and the references therein. Among these researches, Perez-Saez et al. found most consistent and positive correlations between cholera seasonality and precipitation \cite{Perez-Saez}, and Shackleton et al. identified rainfall as the main driver of cholera outbreaks during monsoon seasons \cite{Shackleton}. Therefore, it is highly necessary to consider rainfall when developing mathematical models for diseases transmitted via the faecal-oral route.

In this paper, we partition the year into wet and dry seasons. Then, we develop a two-season switching faecal-oral model to examine how impulsive intervention and seasonal switching together affect the spread of the diseases. In this model, the dry season is represented by a fixed boundary framework, while the wet season is described by a free boundary framework, with impulse interventions occurring at the end of each wet season. For some related two-season models, we refer to \cite{Peng-Zhao-DCDS} for a diffusive logistic model with a free boundary and seasonal succession, \cite{Wang-Zhang-Zhao-JDE} for a diffusive competition model with seasonal succession and different free boundaries, \cite{Li-Sun-JDE} for a free boundary problem with advection and seasonal succession, \cite{Li-Han-Sun-ANA} for a free boundary model with seasonal succession and impulsive harvesting, and the references therein.

The simultaneous introduction of impulsive interventions and seasonal switching renders the model discontinuous at impulse points and invalidates certain analytical approaches. Therefore, new analytical techniques are required to overcome these challenges. We are curious about whether a  spreading-vanishing dichotomy still holds, what a new criteria for spreading and vanishing are, into what form the smoothness of solution may degenerate, and how impulsive interventions and seasonal switching influence the dynamical behavior of disease transmitted by the faecal-oral route. These are the research motivations of this paper.

The rest of our paper is arranged as follows. In \autoref{Section-2}, we develop a two-season switching faecal-oral model with impulsive interventions and free boundaries, and establish some preparatory results. \autoref{Section-3} treats the associated fixed boundary problem and the corresponding eigenvalue problem. A spreading-vanishing dichotomy and a sharp criteria governing this dichotomy are presented in \autoref{Section-4}. \autoref{Section-5} employs numerical examples to validate the theoretical findings and to visually illustrate the influence of impulsive interventions and seasonal switching on the transmission dynamics of faecal-oral diseases. Finally, a brief conclusion of this paper and future work are provided in \autoref{section-6}.
\section{Model formulation and preliminaries}\label{Section-2}
In this section, we develop a two-season switching faecal-oral model with impulsive intervention as well as establish some preliminary results that include a prior estimate, the well-posedness of the solution, two comparison principles, and the definitions of the upper and lower solutions and spreading and vanishing.
\subsection{Model formulation}
In \cite{Wang-Du-DCDSB}, Wang and Du considered a diffusive faecal-oral model with free boundaries, which has the form
\begin{eqnarray}\label{Wang-Du-DCDSB}
\left\{
\begin{array}{ll}
u_{t}=d_{1}\Delta u-a_{11}u+a_{12}v,\; &\,  x\in(r(t),s(t)), t>0, \\[2mm]
v_{t}=d_{2}\Delta v-a_{22}v+f(u),\; &\,  x\in(r(t),s(t)), t>0,\\[2mm]
u(x,t)=0, v(x,t)=0,\; &\,  x\in\{r(t), s(t)\},t>0,\\[2mm]
r'(t)=-\mu_{1}u_{x}(r(t),t)-\mu_{2}v_{x}(r(t),t),\; &\, t>0,\\[2mm]
s'(t)=-\mu_{1}u_{x}(s(t),t)-\mu_{2}v_{x}(s(t),t),\; &\, t>0,\\[2mm]
r(0)=-s_{0}, u(x,0)=u_{0}(x),\; &\,x\in[-s_{0},s_{0}],\\[2mm]
s(0)=s_{0},~~ v(x,0)=v_{0}(x),\; &\,x\in[-s_{0},s_{0}]
\end{array} \right.
\end{eqnarray}
with suitable initial functions $u_{0}(x)$ and $u_{0}(x)$, where $\mu_{2}=\rho\mu_{1}$. In model \eqref{Wang-Du-DCDSB}, $u(x,t)$ represents the spatial density of the infectious agents (bacteria, virus, etc.) at position $x$ and time $t$, while $v(x,t)$ denotes the spatial density of the infective human population at position $x$ and time $t$. All constants are positive. The biological meanings of symbols in model \eqref{Wang-Du-DCDSB} are listed in \autoref{Tab:01}.
\begin{table}[!htb]
\centering
\caption{List of notations and their meanings in model \eqref{Wang-Du-DCDSB}}
\vspace{0.1cm}
\label{Tab:01}
\small
\begin{tabular}{ll}
\hline
 Notation& Biological meaning\\
\hline
$a_{11}$  & Natural death rate of the agents\\
$a_{12}$  & Growth rate of the agents contributed by the infective humans\\
$a_{22}$  & Fatality rate of the infective human population \\
$r(t)$ & Left boundary of the infected region at time $t$\\
$s(t)$ & Right boundary of the infected region at time $t$\\
$d_{1}$  & Diffusion coefficient of the infectious agents\\
$d_{2}$  & Diffusion coefficient of the infective human \\
$\mu_{1}$   &Expansion capacity of the infectious agents \\
$\mu_{2}$   &Expansion capacity of the infective human\\
$s_{0}$   &Length of the right boundary of the initial infection region\\
\hline
\end{tabular}
\end{table}
The function $f(u)$ denotes the infective rate of humans, and satisfies
\begin{eqnarray*}
\bf{(H):}
\left\{
\begin{array}{l}
f(0)=0, ~f\in \mathcal{C}^{1}([0,\infty))~\text{and}~f'(u)>0~\text{for}~u\geq 0, \\[2mm]
\frac{f(u)}{u}~\text{is~strictly~decreasing}~\text{for}~u>0~\text{and}~\lim\limits_{u\rightarrow+\infty}\frac{f(u)}{u}<\frac{a_{11}a_{22}}{a_{12}}.
\end{array} \right.
\end{eqnarray*}

By considering impulsive interventions in model \eqref{Wang-Du-DCDSB}, Zhou, Lin and Pedersen developed an impulsive faecal-oral model in a moving infected environment \cite{Zhou-Lin-Pedersen-ARXIV}, which reads as follows
\begin{eqnarray}\label{Zhou-Lin-Pedersen}
\left\{
\begin{array}{ll}
u_{t}=d_{1}\Delta u-a_{11}u+a_{12}v,\; &\,  x\in(r(t),s(t)), t\in((kT)^{+}, (k+1)T],\\[2mm]
v_{t}=d_{2}\Delta v-a_{22}v+f(u),\; &\,  x\in(r(t),s(t)),t\in((kT)^{+}, (k+1)T], \\[2mm]
u(x,t)=0, v(x,t)=0,\; &\, x\in\{r(t), s(t)\},t\in(kT, (k+1)T],\\[2mm]
u(x,(kT)^{+})=H(u(x,kT)), \; &\, x\in(r(kT), s(kT)), \\[2mm]
v(x,(kT)^{+})=v(x,kT), \; &\,  x\in(r(kT), s(kT)), \\[2mm]
r'(t)=-\mu_{1}u_{x}(r(t),t)-\mu_{2}v_{x}(r(t),t),\; &\, t\in(kT, (k+1)T],\\[2mm]
s'(t)=-\mu_{1}u_{x}(s(t),t)-\mu_{2}v_{x}(s(t),t),\; &\, t\in(kT, (k+1)T],\\[2mm]
r(0)=-s_{0}, u(x,0)=u_{0}(x),\; &\,x\in[-s_{0},s_{0}],\\[2mm]
s(0)=s_{0},~~~ v(x,0)=v_{0}(x),\; &\,x\in[-s_{0},s_{0}], ~k=0,1,2,3,4,\cdots
\end{array} \right.
\end{eqnarray}
with suitable initial functions $u_{0}(x)$ and $v_{0}(x)$, where $T$ is the time span between two neighboring impulsive interventions and $(kT)^{+}$ denotes the right limit.
The impulsive function $H(u)$ is assumed to satisfy
\begin{eqnarray*}
\bf{(H):}
\left\{
\begin{array}{l}
H(0)=0, ~H\in \mathcal{C}^{2}([0,\infty))~\text{and}~H'(u)>0~\text{for}~u\geq 0, \\[2mm]
\frac{H(u)}{u}~\text{is~strictly~decreasing}~\text{and}~0<\frac{H(u)}{u}\leq1~\text{for}~u>0.
\end{array} \right.
\end{eqnarray*}
In model \eqref{Zhou-Lin-Pedersen}, the other symbols have the same biological meaning as in model \eqref{Wang-Du-DCDSB}.

Based on model \eqref{Zhou-Lin-Pedersen}, this paper develops a two-season switching faecal-oral model with free boundaries, which has the form
\begin{eqnarray}\label{Fixed-1}
\left\{
\begin{array}{ll}
u_{t}=-\delta_{1} u ,\; &\,x\in(r(t),s(t)), t\in((mT)^{+}, mT+\tau],  \\[2mm]
v_{t}=d_{2}\Delta v-\delta_{2} v ,\; &\,x\in(r(t),s(t)), t\in((mT)^{+}, mT+\tau],  \\[2mm]
u(x,(mT)^{+})=H(u(x,mT)),\; &\, x\in(r(mT),s(mT)),  \\[2mm]
v(x,(mT)^{+})=v(x,mT),\; &\, x\in(r(mT),s(mT)),  \\[2mm]
u(x,t)=v(x,t)=0,\; &\, x\in\{r(t),s(t)\},  t>0,\\[2mm]
r(t)=r(mT),~s(t)=s(mT),\; &\,  t\in[mT, mT+\tau],  \\[2mm]
r(0)=-s_{0}, u(0,x)=u_{0}(x),\; &\,x\in[-s_{0},s_{0}],\\[2mm]
s(0)=s_{0},~~ v(0,x)=v_{0}(x),\; &\,x\in[-s_{0},s_{0}]
\end{array} \right.
\end{eqnarray}
for dry season, and
 \begin{eqnarray}\label{Fixed-2}
\left\{
\begin{array}{ll}
u_{t}=d_{1}\Delta u-a_{11}u+a_{12}v ,\; &\,x\in(r(t),s(t)), t\in(mT+\tau, (m+1)T],  \\[2mm]
v_{t}=d_{2}\Delta v-a_{22}v+f(u) ,\; &\,x\in(r(t),s(t)), t\in(mT+\tau, (m+1)T],  \\[2mm]
u(x,t)=v(x,t)=0,\; &\, x\in\{r(t),s(t)\},  t>0,\\[2mm]
r'(t)=-\mu_{1}u_{x}(r(t),t)-\mu_{2}v_{x}(r(t),t),\; &\,  t\in(mT+\tau, (m+1)T),  \\[2mm]
s'(t)=-\mu_{1}u_{x}(s(t),t)-\mu_{2}v_{x}(s(t),t),\; &\,  t\in(mT+\tau, (m+1)T) \\[2mm]
\end{array} \right.
\end{eqnarray}
for wet season. The initial functions $u_{0}(x)$ and $v_{0}(x)$ in \eqref{Fixed-1} are assumed to satisfy
\begin{eqnarray}\label{Fixed-3}
\left\{
\begin{array}{l}
u_{0}\in\mathcal{C}^{2}([-s_{0},s_{0}]), ~u_{0}(\pm s_{0})=0~\text{and}~u_{0}(x)>0~\text{in}~(-s_{0},s_{0}), \\[2mm]
v_{0}\in\mathcal{C}^{2}([-s_{0},s_{0}]), ~v_{0}(\pm s_{0})=0~\text{and}~v_{0}(x)>0~~\text{in}~(-s_{0},s_{0}).
\end{array} \right.
\end{eqnarray}
In \eqref{Fixed-1}, $\tau\in(0, T)$ denotes the duration of dry season, and $\delta_{1}$ and $\delta_{2}$ denote the natural death rate of the agents and the fatality rate of the infective human population in dry season, respectively. Here and in what follows, unless stated otherwise, we always take $m=0,1,2, \cdots$, and the initial, growth, and impulsive functions are always satisfied \eqref{Fixed-3}, assumption \textbf{(F)}, and assumption \textbf{(H)}. The other symbols of model \eqref{Fixed-1}-\eqref{Fixed-3}
have the same biological meaning as in model \eqref{Zhou-Lin-Pedersen}.
It is worth mentioning that the models in \cite{Ahn-Baek-Lin-AMM,Zhao-Li-DCDSB,Wang-Du-DCDSB,Zhou-Lin-Pedersen-ARXIV,7,Capasso-Maddalena-3} are special cases of model \eqref{Fixed-1}-\eqref{Fixed-3}.

\begin{remark}
  We note that in our model, the impulsive intervention acts only on the
$u$-component, while the $v$-component remains continuous in time. 
\end{remark}

\subsection{Preliminary results}
For convenience, we first introduce some notations. For given $m, \tau, T>0$, we define
\begin{equation*}
\mathcal{T}(m):=((mT)^{+}, mT+\tau]\cup (mT+\tau, (m+1)T],
\end{equation*}
\begin{equation*}
\mathcal{PC}^{1}(m):=\mathcal{C}^{1}([(mT)^{+}, mT+\tau]\cup [mT+\tau, (m+1)T]).
\end{equation*}
For $r(t),s(t)\in\mathcal{PC}^{1}(m)$, we define
\begin{equation}\label{pc}
\mathcal{PC}^{2,1}(r(t), s(t), m):=\mathcal{C}^{2,1}( [r(t), s(t)]\times \mathcal{T}(m)).
\end{equation}

By using parabolic regularity theory and the contraction mapping principle, we can easily obtain the existence and uniqueness of nonnegative classical solution $(u,v,r,s)$ of model \eqref{Fixed-1}-\eqref{Fixed-3}, which satisfies $r,s\in\mathcal{PC}^{1}(m)$ and $u,v\in\mathcal{PC}^{2,1}(r(t), s(t), m)$. In addition, we define
\begin{equation}\label{2-2}
C_{2}:=\max\Big\{u^{*}, \|u_{0}\|_{\mathcal{C}([-s_{0},s_{0}])}, \frac{a_{12}}{a_{11}}\|v_{0}\|_{\mathcal{C}([-s_{0},s_{0}])}\Big\}, ~~C_{3}:=\max\Big\{\|v_{0}\|_{\mathcal{C}([-s_{0},s_{0}])}, \frac{f(C_{2})}{a_{22}}\Big\},
\end{equation}
where $u^{*}=0$ if $\frac{a_{12}G'(0)}{a_{11}a_{22}}\leq1,$ and $u^{*}$ is uniquely determined by $\frac{G(u)}{u}=\frac{a_{11}a_{22}}{a_{12}}$ if $\frac{a_{12}G'(0)}{a_{11}a_{22}}>1$. It is easily seen that
\begin{equation}\label{2-2-0}
-a_{11}C_{2}+a_{12}C_{3}\leq 0 \text{~and~}-a_{22}C_{3}+f(C_{2})\leq 0.
\end{equation}
Then, it is also easy to obtain from the strong maximum principle and the Hopf boundary lemma for cooperative systems that there exists positive constants $C_{1}$ such that
\begin{equation}\label{2-1}
\begin{aligned}
0&<-r'(t), s'(t)< C_{1} \text{~for~}t\in(mT+\tau, (m+1)T],\\
0<u(x,t)&<C_{2}\text{~and}~0<v(x,t)< C_{3} \text{~for~}t\in[0, \infty)\text{~and~}x\in(r(t), s(t)),
\end{aligned}
\end{equation}
where $C_{2}$ and $C_{3}$ are defined in \eqref{2-2}.
For more details, see \cite[Lemma 2.2]{Zhou-Lin-Pedersen-ARXIV}.

Next, we present two parabolic-type comparison principles, which are the basis of this research.
\begin{lemma}\label{lemma 2-1}
Suppose that $M$ is a positive integer, $\overline{r}, \overline{s}\in \mathcal{C}((0, MT])\cap \mathcal{PC}^{1}(i)$, $\overline{u}, \overline{v}\in \mathcal{PC}^{2,1}(\overline{r}(t), \overline{s}(t), i)$, and
\begin{eqnarray}\label{2-3}
\left\{
\begin{array}{ll}
\overline{u}_{t}\geq-\delta_{1} \overline{u} ,\; &\,x\in(\overline{r}(t),\overline{s}(t)), t\in((iT)^{+}, iT+\tau],  \\[2mm]
\overline{v}_{t}\geq d_{2}\Delta \overline{v}-\delta_{2} \overline{v} ,\; &\,x\in(\overline{r}(t),\overline{s}(t)), t\in((iT)^{+}, iT+\tau],  \\[2mm]
\overline{u}_{t}\geq d_{1}\Delta \overline{u}-a_{11}\overline{u}+a_{12}\overline{v} ,\; &\,x\in(\overline{r}(t),\overline{s}(t)), t\in(iT+\tau, (i+1)T],  \\[2mm]
\overline{v}_{t}\geq d_{2}\Delta \overline{v}-a_{22}\overline{v}+f(\overline{u}) ,\; &\,x\in(\overline{r}(t),\overline{s}(t)), t\in(iT+\tau, (i+1)T],  \\[2mm]
\overline{r}(t)\leq \overline{r}(iT),~\overline{s}(t)\geq\overline{s}(iT),\; &\,  t\in[iT, iT+\tau],  \\[2mm]
\overline{r}'(t)\leq -\mu_{1}\overline{u}_{x}(\overline{r}(t),t)-\mu_{2}\overline{v}_{x}(\overline{r}(t),t),\; &\,  t\in(iT+\tau, (i+1)T),  \\[2mm]
\overline{s}'(t)\geq-\mu_{1}\overline{u}_{x}(\overline{s}(t),t)-\mu_{2}\overline{v}_{x}(\overline{s}(t),t),\; &\,  t\in(iT+\tau, (i+1)T),  \\[2mm]
\overline{u}(x,t)= \overline{v}(x,t)=0,\; &\, x\in\{\overline{r}(t),\overline{s}(t)\},  t\in(0,MT],\\[2mm]
\overline{u}(x,(iT)^{+})\geq H(\overline{u}(x,iT)),\; &\, x\in(\overline{r}(t),\overline{s}(t)),  \\[2mm]
\overline{v}(x,(iT)^{+})\geq \overline{v}(x,iT),\; &\, x\in(\overline{r}(t),\overline{s}(t)),
\end{array} \right.
\end{eqnarray}
where $i=0,1,2, \cdots, M-1$. If
\begin{equation*}
\overline{r}(0)\leq -s_{0},~ \overline{s}(0)\geq s_{0},
\end{equation*}
and
\begin{equation*}
\overline{u}(x,0)\geq u_{0}(x), ~\overline{v}(x,0)\geq v_{0}(x)\text{~for~}x\in[-s_{0}, s_{0}],
\end{equation*}
then the solution $(u,v,r,s)$ of model \eqref{Fixed-1}-\eqref{Fixed-3} satisfies
\begin{equation*}
r(t)\geq \overline{r}(t),~s(t)\leq \overline{s}(t)\text{~in~}[0, MT],
\end{equation*}
and
\begin{equation*}
u(x,t)\leq \overline{u}(x,t), v(x,t)\leq \overline{v}(x,t)\text{~for~}x\in[r(t), s(t)]\text{~and~}t\in[0, MT].
\end{equation*}
\begin{proof}
We first assume that $\overline{r}(0)< -s_{0}<s_{0}<\overline{s}(0)$, and prove that $\overline{r}(t)< r(t)<s(t)<\overline{s}(t)$ for $t\in[0, MT]$. If our assertion does not hold, then we can find a first $t^{*}\leq MT$ such that
$\overline{r}(t^{*})= r(t^{*})$ or $\overline{s}(t^{*})= s(t^{*})$. Specifically, we have from the fifth equation of \eqref{2-3} that the above defined $t^{*}\in (i^{*}T+\tau, (i^{*}+1)T]$ for some  $i^{*}\in \{0, 1, \cdots, M-1 \}$.  Without loss of generality, we assume that $\overline{s}(t^{*})= s(t^{*})$ and $\overline{r}(t^{*})\leq r(t^{*})$, which implies that
\begin{equation}\label{2-6}
s'(t^{*})\geq\overline{s} '(t^{*}),~s(t)<\overline{s}(t),~ r(t)>\overline{r}(t)\text{~for~}t\in[0, t^{*}).
\end{equation}

Now, we go to achieve the contradiction. Firstly, consider the following problem
\begin{eqnarray*}
\left\{
\begin{array}{ll}
\overline{u}_{t}\geq-\delta_{1} \overline{u} ,\; &\,x\in(\overline{r}(t),\overline{s}(t)), t\in((jT)^{+}, jT+\tau],  \\[2mm]
\overline{v}_{t}\geq d_{2}\Delta \overline{v}-\delta_{2} \overline{v} ,\; &\,x\in(\overline{r}(t),\overline{s}(t)), t\in((jT)^{+}, jT+\tau],  \\[2mm]
\overline{u}_{t}\geq d_{1}\Delta \overline{u}-a_{11}\overline{u}+a_{12}\overline{v} ,\; &\,x\in(\overline{r}(t),\overline{s}(t)), t\in\Sigma(t^{*}),  \\[2mm]
\overline{v}_{t}\geq d_{2}\Delta \overline{v}-a_{22}\overline{v}+f(\overline{u}) ,\; &\,x\in(\overline{r}(t),\overline{s}(t)), t\in\Sigma(t^{*}),  \\[2mm]
\overline{u}(x,t)= \overline{v}(x,t)=0,\; &\, x\in\{\overline{r}(t),\overline{s}(t)\},  t\in(0,t^{*}],\\[2mm]
\overline{u}(x,(jT)^{+})\geq H(\overline{u}(x,jT)),\; &\, x\in(\overline{r}(t),\overline{s}(t)),  \\[2mm]
v(x,(iT)^{+})\geq v(x,iT),\; &\, x\in(\overline{r}(t),\overline{s}(t)),\\[2mm]
\overline{r}(0)\leq -s_{0},~ \overline{s}(0)\geq s_{0},\; &\ \\[2mm]
\overline{u}(x,0)\geq u_{0}(x), ~\overline{v}(x,0)\geq v_{0}(x),\; &\, x\in [-s_{0}, s_{0}],
\end{array} \right.
\end{eqnarray*}
where $j=0,1,2, \cdots, i^{*}$ and $\Sigma(t^{*})=((j-1)T+\tau, jT]\cup (jT+\tau, t^{*}]$. Then, it follows from the maximum principle that $\overline{u}(x,t)\geq 0$ and $\overline{v}(x,t)\geq 0$ for $x\in[r(t), s(t)]$ and $t\in(0, t^{*}]$. We define $\phi:=\overline{u}-u$ and $\psi:=\overline{v}-v$, and then it follows that
\begin{eqnarray}\label{2-4}
\left\{
\begin{array}{ll}
\phi_{t}\geq-\delta_{1} \phi ,\; &\,x\in(r(t),s(t)), t\in(0^{+}, \tau],  \\[2mm]
\psi_{t}\geq d_{2}\Delta \psi-\delta_{2} \psi ,\; &\,x\in(r(t),s(t)), t\in(0^{+}, \tau],  \\[2mm]
\phi_{t}\geq d_{1}\Delta \phi-a_{11}\phi+a_{12}\psi ,\; &\,x\in(r(t),s(t)), t\in(\tau, (i+1)T],  \\[2mm]
\psi_{t}\geq d_{2}\Delta \psi-a_{22}\psi+f'(\theta)\phi ,\; &\,x\in(r(t),s(t)), t\in(\tau, (i+1)T],  \\[2mm]
\phi(x,t)\geq 0,~~\psi(x,t)\geq 0,\; &\, x\in\{r(t),s(t)\},  t\in(0,T],\\[2mm]
\phi(x,0^{+})\geq 0,~~\psi(x,0^{+})\geq 0\; &\, x\in(r(t),s(t)),  \\[2mm]
\phi(x,0)\geq 0, ~\psi(x,0)\geq 0,\; &\, x\in [-s_{0}, s_{0}],
\end{array} \right.
\end{eqnarray}
where $\min\{\overline{u},u\}\leq \theta\leq\max\{\overline{u},u\} $ for $x\in(r(t),s(t))$ and $t\in(\tau, (i+1)T]$.
The maximum principle yields that $\phi(x,t)\geq 0$ and $\psi(x,t)\geq 0$ for $x\in[r(t), s(t)]$ and $t\in(0, T]$.
The same procedure gives that $\phi(x,t)\geq 0$ and $\psi(x,t)\geq 0$ for $x\in[r(t), s(t)]$ and $t\in(0,t^{*}]$, which implies that
\begin{equation}\label{2-5}
u(x,t)\leq \overline{u}(x,t)\text{~and~} v(x,t)\leq \overline{v}(x,t)
\end{equation}
in $[r(t), s(t)]\times (0,t^{*}]$. Since $\overline{r}(t)<r(t)< s(t)<\overline{s}(t)$ for $t\in(0, t^{*})$, we have that $\phi(x,t)\not\equiv 0$ and $\psi(x,t)\not\equiv  0$ in $[r(t), s(t)]\times (0,t^{*}]$. Additionally, it follows from \eqref{2-4} and \eqref{2-5} that
\begin{eqnarray*}
\left\{
\begin{array}{ll}
\phi_{t}\geq-\delta_{1} \phi ,\; &\,x\in(r(t),s(t)), t\in((jT)^{+}, jT+\tau],  \\[2mm]
\psi_{t}\geq d_{2}\Delta \psi-\delta_{2} \psi ,\; &\,x\in(r(t),s(t)), t\in((jT)^{+}, jT+\tau],  \\[2mm]
\phi_{t}\geq d_{1}\Delta \phi-a_{11}\phi ,\; &\,x\in(r(t),s(t)), t\in\Sigma(t^{*}),  \\[2mm]
\psi_{t}\geq d_{2}\Delta \psi-a_{22}\psi ,\; &\,x\in(r(t),s(t)), t\in\Sigma(t^{*}),  \\[2mm]
\phi(x,t)\geq 0,~~\psi(x,t)\geq 0,\; &\, x\in\{r(t),s(t)\},  t\in(0,t^{*}],\\[2mm]
\phi(s(t^{*}),t^{*})=\psi(s(t^{*}),t^{*})=0,\; &\ \\[2mm]
\phi(x,(jT)^{+})\geq 0,~~\psi(x,(jT)^{+})\geq 0\; &\, x\in(r(t),s(t)),  \\[2mm]
\phi(x,0)\geq 0, ~\psi(x,0)\geq 0,\; &\, x\in [-s_{0}, s_{0}].
\end{array} \right.
\end{eqnarray*}
Then, we have from the strong maximum principle and Hopf boundary lemma that $\phi(x,t), \psi(x,t)>0$ in $[r(t), s(t)]\times (0,t^{*}]\setminus\{(r(t^{*}), t^{*}), (s(t^{*}), t^{*})\}$ and $\phi_{x}(s(t^{*}), t^{*}), \psi_{x}(s(t^{*}), t^{*})<0$, which implies that $s'(t^{*})<\overline{s} '(t^{*})$. But this contradicts \eqref{2-6},
which proves our assertion that $\overline{r}(t)< r(t)<s(t)<\overline{s}(t)$ for $t\in[0, MT]$.

Now, we may apply the usual comparison principle over $[r(t), s(t)]\times[0, MT]$ to conclude that $u(x,t)\leq \overline{u}(x,t)$ and $v(x,t)\leq \overline{v}(x,t)$ in $[r(t), s(t)]\times[0, MT]$. By employing an approximation argument, the proof of the general case $\overline{r}(0)\leq -s_{0}<s_{0}\leq \overline{s}(0)$ can be completed. Since the procedure is standard, we omit the details here and refer the interested reader to \cite[lemma 2.3]{Wang-Du-DCDSB} for further details. This ends the proof.
\end{proof}
\end{lemma}
\begin{lemma}\label{lemma 2-1-B}
Suppose that $M$ is a positive integer, $\overline{s}(t)\equiv s(t)$ for $t\in[0, MT]$, $\overline{r}\in \mathcal{C}((0, MT])\cap \mathcal{PC}^{1}(i)$, $\overline{u}, \overline{v}\in \mathcal{PC}^{2,1}(\overline{r}(t), s(t), i)$, and
\begin{eqnarray}\label{2-3-B}
\left\{
\begin{array}{ll}
\overline{u}_{t}\geq-\delta_{1} \overline{u} ,\; &\,x\in(\overline{r}(t),s(t)), t\in((iT)^{+}, iT+\tau],  \\[2mm]
\overline{v}_{t}\geq d_{2}\Delta \overline{v}-\delta_{2} \overline{v} ,\; &\,x\in(\overline{r}(t),s(t)), t\in((iT)^{+}, iT+\tau],  \\[2mm]
\overline{u}_{t}\geq d_{1}\Delta \overline{u}-a_{11}\overline{u}+a_{12}\overline{v} ,\; &\,x\in(\overline{r}(t),s(t)), t\in(iT+\tau, (i+1)T],  \\[2mm]
\overline{v}_{t}\geq d_{2}\Delta \overline{v}-a_{22}\overline{v}+f(\overline{u}) ,\; &\,x\in(\overline{r}(t),s(t)), t\in(iT+\tau, (i+1)T],  \\[2mm]
\overline{u}(\overline{r}(t),t)\geq u(\overline{r}(t),t),~~ \overline{v}(s(t),t)=0,\; &\, t\in(0,MT],\\[2mm]
\overline{u}(x,(iT)^{+})\geq H(\overline{u}(x,iT)),\; &\, x\in(\overline{r}(t),s(t)),  \\[2mm]
\overline{v}(x,(iT)^{+})\geq \overline{v}(x,iT),\; &\, x\in(\overline{r}(t),s(t)),
\end{array} \right.
\end{eqnarray}
where $i=0,1,2, \cdots, M-1$. If
\begin{equation*}
\overline{r}(t)\geq r(t)\text{~for~}t\in[0, MT],
\end{equation*}
and
\begin{equation*}
\overline{u}(x,0)\geq u_{0}(x), ~\overline{v}(x,0)\geq v_{0}(x)\text{~for~}x\in[\overline{r}(0), s_{0}],
\end{equation*}
then the solution $(u,v,r,s)$ of model \eqref{Fixed-1}-\eqref{Fixed-3} satisfies
\begin{equation*}
u(x,t)\leq \overline{u}(x,t), v(x,t)\leq \overline{v}(x,t)\text{~for~}x\in[\overline{r}(t), s(t)]\text{~and~}t\in[0, MT].
\end{equation*}
\begin{proof}
The proof of this lemma can be obtained by a minor modification of the proof of \autoref{lemma 2-1}. Since the arguments are quite similar, we omit the details here.
\end{proof}
\end{lemma}
\begin{remark}\label{remark 2.1}
We call $(\overline{u},\overline{v},\overline{r},\overline{s})$ in \autoref{lemma 2-1} and \autoref{lemma 2-1-B} an upper solution of model \eqref{Fixed-1}-\eqref{Fixed-3}.
We can define a lower solution by reversing all the inequalities, and there is an analogue of \autoref{lemma 2-1} and \autoref{lemma 2-1-B} for lower solutions.
\end{remark}

It follows from the first inequality of \eqref{2-1} that
\begin{equation*}
\lim\limits_{t\rightarrow\infty} r(t)=:r_{\infty}\in [-\infty, -s_{0})\text{~and~}\lim\limits_{t\rightarrow\infty} s(t)=:s_{\infty}\in (s_{0}, +\infty]
\end{equation*}
are always well-defined. At the end of this section, we provide the definitions of spreading and vanishing.
\begin{definition}\label{definition 2.1}
We call the unique solution $(u,v, r, s)$ of model \eqref{Fixed-1}-\eqref{Fixed-3}  vanishing if
\begin{equation*}
s_{\infty}-r_{\infty}<\infty~\text{and}~\lim\limits_{t\rightarrow\infty}\|u(x,t)\|_{\mathcal{C}([r_{\infty}, s_{\infty}])}+\lim\limits_{t\rightarrow\infty}\|v(x,t)\|_{\mathcal{C}([r_{\infty}, s_{\infty}])}=0,
\end{equation*}
and spreading if
\begin{equation*}
s_{\infty}-r_{\infty}=\infty~\text{and}~\lim\limits_{t\rightarrow\infty}u(x,t)=Z(t)~\text{as~well~as}~\lim\limits_{t\rightarrow\infty}v(x,t)=W(t)
\end{equation*}
locally uniformly in $\mathcal{R}$ and uniformly in $[0,T]$, where $(Z(t),W(t))$ denotes the solution of \eqref{3-42}.
\end{definition}
\section{The corresponding fixed boundary problem}\label{Section-3}
This section investigates the fixed boundary problem associated with model \eqref{Fixed-1}-\eqref{Fixed-3}, which is given by
\begin{eqnarray}\label{3-1}
\left\{
\begin{array}{ll}
u_{t}=-\delta_{1} u ,\; &\,x\in(l_{1},l_{2}), t\in((mT)^{+}, mT+\tau],  \\[2mm]
v_{t}=d_{2}\Delta v-\delta_{2} v ,\; &\,x\in(l_{1},l_{2}), t\in((mT)^{+}, mT+\tau],  \\[2mm]
u_{t}=d_{1}\Delta u-a_{11}u+a_{12}v ,\; &\,x\in(l_{1},l_{2}), t\in(mT+\tau, (m+1)T],  \\[2mm]
v_{t}=d_{2}\Delta v-a_{22}v+f(u) ,\; &\,x\in(l_{1},l_{2}), t\in(mT+\tau, (m+1)T],  \\[2mm]
u(x,t)=v(x,t)=0,\; &\, x\in\{l_{1},l_{2}\},  t\in(0, \infty),\\[2mm]
u(x,(mT)^{+})=H(u(x,mT)),\; &\, x\in(l_{1},l_{2}),  \\[2mm]
v(x,(mT)^{+})=v(x,mT),\; &\, x\in(l_{1},l_{2}),  \\[2mm]
u(0,x)=u_{0}(x),v(0,x)=v_{0}(x),\; &\,x\in[l_{1},l_{2}],\\[2mm]
\end{array} \right.
\end{eqnarray}
where $-l_{1}, l_{2}\in\mathds{R}^{+}$, and the initial functions $u_{0}(x)$ both $v_{0}(x)$ satisfy \eqref{Fixed-3} with $-s_{0}$ and $s_{0}$ replaced by $l_{1}$ and $l_{2}$,
respectively. A good understanding of problem \eqref{3-1} plays an important role in analysing the long time dynamical behaviours of model \eqref{Fixed-1}-\eqref{Fixed-3}.
However, problem \eqref{3-1} has not been studied in previous work. Therefore, we examine problem \eqref{3-1} here. For the sake of convenience in studying
problem \eqref{3-1} later, we present some properties of its solution. Since their proofs are simple, we omit the proof process.
\begin{lemma}\label{lemma 3.0.1}
Problem \eqref{3-1} admits a unique nonnegative global classical solution
\begin{equation*}
(u,v)\in \mathcal{PC}^{2,1}(l_{1},l_{2}, m),
\end{equation*}
where $\mathcal{PC}^{2,1}(l_{1},l_{2}, m)$ is defined in \eqref{pc}. Moreover, the solution of problem \eqref{3-1} satisfies that
\begin{equation*}
0<u(x,t)<C_{2}\text{~and}~0<v(x,t)< C_{3} \text{~for~}(x,t)\in(l_{1},l_{2})\times[0, \infty),
\end{equation*}
where $C_{2}$ and $C_{3}$ are defined in \eqref{2-2}.
\end{lemma}
\begin{lemma}\label{lemma 3.0.2}
Suppose that $M$ is a positive integer, $(\overline{u}, \overline{v})\in \mathcal{PC}^{2,1}(l_{1},l_{2}, i)$, and
\begin{eqnarray*}
\left\{
\begin{array}{ll}
\overline{u}_{t}\geq -\delta_{1} \overline{u} ,\; &\,x\in(l_{1},l_{2}), t\in((iT)^{+}, iT+\tau],  \\[2mm]
\overline{v}_{t}\geq d_{2}\Delta \overline{v}-\delta_{2} \overline{v} ,\; &\,x\in(l_{1},l_{2}), t\in((iT)^{+}, iT+\tau],  \\[2mm]
\overline{u}_{t}\geq d_{1}\Delta \overline{u}-a_{11}\overline{u}+a_{12}\overline{v} ,\; &\,x\in(l_{1},l_{2}), t\in(iT+\tau, (i+1)T],  \\[2mm]
\overline{v}_{t}\geq d_{2}\Delta\overline{v}-a_{22}\overline{v}+f(\overline{u}) ,\; &\,x\in(l_{1},l_{2}), t\in(iT+\tau, (i+1)T],  \\[2mm]
\overline{u}(x,t)=\overline{v}(x,t)=0,\; &\, x\in\{l_{1},l_{2}\},  t\in(0, MT],\\[2mm]
\overline{u}(x,(iT)^{+})\geq H(\overline{u}(x,iT)),\; &\, x\in(l_{1},l_{2}),  \\[2mm]
\overline{v}(x,(iT)^{+})\geq \overline{v}(x,iT),\; &\, x\in(l_{1},l_{2}),  \\[2mm]
\overline{u}(0,x)\geq u_{0}(x),\overline{v}(0,x)\geq v_{0}(x),\; &\,x\in[l_{1},l_{2}],
\end{array} \right.
\end{eqnarray*}
where $i=1,2, \cdots, M-1$, and $\mathcal{PC}^{2,1}(l_{1},l_{2}, i)$ is defined in \eqref{pc}. Then, the solution of problem \eqref{3-1} satisfies
\begin{equation*}
u(x,t)\leq \overline{u}(x,t) \text{~and~}v(x,t)\leq \overline{v}(x,t) \text{~for~}(x,t)\in[l_{1},l_{2}]\times [0, MT].
\end{equation*}
\end{lemma}

\begin{remark}\label{remark 3.0.1}
The pair $(\overline{u}, \overline{v})$ in \autoref{lemma 3.0.2} is usually called an upper solution of problem \eqref{3-1}. We can also define a
lower solution $(\underline{u}, \underline{v})$ by reversing all the inequalities, and there is an analogue of \autoref{lemma 3.0.2} for lower solutions.
\end{remark}
\subsection{An analysis of the eigenvalue problem}\label{Section-3-1}
The dynamical behaviours of problem \eqref{3-1} are classified by its eigenvalue problem linearised at $(0,0)$, which has the form
\begin{eqnarray}\label{3-2}
\left\{
\begin{array}{ll}
\phi_{t}=-\delta_{1} \phi+\lambda \phi ,\; &\,x\in(l_{1},l_{2}), t\in(0^{+}, \tau],  \\[2mm]
\psi_{t}=d_{2}\Delta \psi-\delta_{2} \psi+\lambda \psi ,\; &\,x\in(l_{1},l_{2}), t\in(0^{+}, \tau],  \\[2mm]
\phi_{t}=d_{1}\Delta \phi-a_{11}\phi+a_{12}\psi+\lambda \phi ,\; &\,x\in(l_{1},l_{2}), t\in(\tau, T],  \\[2mm]
\psi_{t}=d_{2}\Delta \psi-a_{22}\psi+f'(0) \phi+\lambda \psi  ,\; &\,x\in(l_{1},l_{2}), t\in(\tau, T],  \\[2mm]
\phi(x,t)=\psi(x,t)=0,\; &\, x\in\{l_{1},l_{2}\},  t\in(0, T],\\[2mm]
\phi(x,0^{+})=H'(0)\phi(x,0),\; &\, x\in(l_{1},l_{2}),  \\[2mm]
\psi(x,0^{+})=\psi(x,0),\; &\, x\in(l_{1},l_{2}),  \\[2mm]
\phi(0,x)=\phi(T,x),\psi(0,x)=\psi(T,x),\; &\,x\in[l_{1},l_{2}].
\end{array} \right.
\end{eqnarray}
For ease of analyzing the principal eigenvalue in problem \eqref{3-2}, we equivalently rewrite problem \eqref{3-2} as the following problem
\begin{eqnarray}\label{3-3}
\left\{
\begin{array}{ll}
\phi_{t}=-\delta_{1} \phi+\lambda \phi ,\; &\,x\in(l_{1},l_{2}), t\in(0, \tau],  \\[2mm]
\psi_{t}=d_{2}\Delta \psi-\delta_{2} \psi+\lambda \psi ,\; &\,x\in(l_{1},l_{2}), t\in(0, \tau],  \\[2mm]
\phi_{t}=d_{1}\Delta \phi-a_{11}\phi+a_{12}\psi+\lambda \phi ,\; &\,x\in(l_{1},l_{2}), t\in(\tau, T],  \\[2mm]
\psi_{t}=d_{2}\Delta \psi-a_{22}\psi+f'(0) \phi+\lambda \psi  ,\; &\,x\in(l_{1},l_{2}), t\in(\tau, T],  \\[2mm]
\phi(x,t)=\psi(x,t)=0,\; &\, x\in\{l_{1},l_{2}\},  t\in(0, T],\\[2mm]
\phi(x,0)=H'(0)\phi(x,T),\; &\, x\in[l_{1},l_{2}],  \\[2mm]
\psi(x,0)=\psi(x,T),\; &\, x\in[l_{1},l_{2}].
\end{array} \right.
\end{eqnarray}

For any $\textbf{x}=(x_{1},x_{2})^{T}\in \mathds{R}^{2}$, we put
\begin{equation*}
\|\textbf{x} \|_{\mathds{R}^{2}}=\sqrt{x_{1}^{2}+x^{2}_{2}}.
\end{equation*}
Let
\begin{equation*}
(\mathds{R}^{2})^{+}=\{\textbf{x}=(x_{1},x_{2})^{T}\in \mathds{R}^{2}|  x_{1},x_{2}\geq 0  \},
\end{equation*}
and
\begin{equation*}
(\mathds{R}^{2})^{++}=\{\textbf{x}=(x_{1},x_{2})^{T}\in \mathds{R}^{2}|  x_{1},x_{2}> 0  \}.
\end{equation*}
For any $\alpha,\beta\in \mathds{R}^{+}$, we set
\begin{equation*}
\mathds{H}^{\alpha,\beta}=\mathcal{C}^{\alpha,\beta}( [l_{1},l_{2}] \times [0,\tau],\mathds{R}^{2}) \cap \mathcal{C}^{\alpha,\beta}(  [l_{1},l_{2}] \times [\tau,T],\mathds{R}^{2})
\end{equation*}
with the norm
\begin{equation*}
\|\textbf{u}\|_{\mathds{H}^{\alpha,\beta}}=\sqrt{\|u_{1}\|^{2}_{\mathcal{C}^{\alpha,\beta}( [l_{1},l_{2}] \times [0,\tau])} +\|u_{2}\|^{2}_{\mathcal{C}^{\alpha,\beta}( [l_{1},l_{2}] \times [0,\tau])} }+\sqrt{\|u_{1}\|^{2}_{\mathcal{C}^{\alpha,\beta}( [l_{1},l_{2}] \times [\tau,T])} +\|u_{2}\|^{2}_{\mathcal{C}^{\alpha,\beta}( [l_{1},l_{2}] \times [\tau,T])} }
\end{equation*}
for all $\textbf{u}=(u_{1}, u_{2})\in \mathds{H}^{\alpha,\beta}$,
and then define
\begin{equation*}
\begin{split}
\mathds{E}=&\big\{\textbf{u}=(u_{1}, u_{2})\in \mathds{H}^{1,0}|~ \textbf{u}(x,t)=\textbf{0}\text{~for~all~} (x,t)\in\{l_{1},l_{2}\}\times [0, T],\\
&~~~~~~~~~~\text{and~}(u_{1}(x,0), u_{2}(x,0))=(H'(0)u_{1}(x,T),u_{2}(x,T)) \text{~for~all~}  x\in [l_{1},l_{2}]\big\}.
\end{split}
\end{equation*}
In fact, $\mathds{E}$ is a real Banach space with a positive cone
\begin{equation*}
\mathds{E}^{+}=\big\{\textbf{u}\in \mathds{E} |~ \textbf{u}\in (\mathds{R}^{2})^{+}\text{~for~all~} (x,t)\in Q^{T}_{0}\text{~and~}-\nabla \textbf{u} \cdot \mathbf{n}\in (\mathds{R}^{2})^{+}\text{~for~all~} (x,t)\in \partial Q^{T}_{0} \big\},
\end{equation*}
where $Q^{T}_{0}=(l_{1},l_{2})\times [0,T]$, $\partial Q^{T}_{0}=\{l_{1},l_{2}\}\times [0, T]$, and $\mathbf{n}$ denotes the unit outward normal vector on $\{l_{1},l_{2}\}$. The positive cone $\mathds{E}^{+}$ has nonempty interior
\begin{equation*}
\begin{split}
\mathds{E}^{++}=\big\{&\textbf{u}\in \mathds{E} |~ \textbf{u}\in (\mathds{R}^{2})^{++}~\forall (x,t)\in Q^{T}_{0} \text{~and~}-\nabla \textbf{u} \cdot \mathbf{n}\in (\mathds{R}^{2})^{++}~\forall (x,t)\in \partial Q^{T}_{0}\big\}.
\end{split}
\end{equation*}
For any given $\textbf{u}, \textbf{v}\in \mathcal{X}$, we write
\begin{equation*}
\begin{split}
&\textbf{u}\succeq\textbf{v}\text{~if~}\textbf{u}-\textbf{v}\in \mathcal{X}^{+},\\
&\textbf{u}\succ\textbf{v}\text{~if~}\textbf{u}-\textbf{v}\in \mathcal{X}^{+} -\textbf{0},\\
&\textbf{u}\succ \succ \textbf{v}\text{~if~}\textbf{u}-\textbf{v}\in \mathcal{X}^{++},
\end{split}
\end{equation*}
where $ \mathcal{X}=\mathds{R}^{2}$ or $\mathds{E}$. Let $\mathcal{A}$ be a continuous self-mapping of $\mathds{E}$. We say that $\mathcal{A}$ is
\begin{equation*}
\text{strongly~order-preserving~if~}\textbf{u}\succ\textbf{v}\Rightarrow \mathcal{A}(\textbf{u})\succ \succ\mathcal{A}(\textbf{v}).
\end{equation*}
With the above preparations, we now present the main results of this subsection.
\begin{theorem}\label{theorem 3.1}
The boundary value problem \eqref{3-3} has exactly one real eigenvalue, which will be denoted by $\lambda_{1}$, to a positive vector-valued eigenfunction.
This eigenvalue is algebraically simple and there is no other eigenvalue to a positive vector-valued eigenfunction. Moreover, the positive vector-valued eigenfunction associated with $\lambda_{1}$, say $(\phi_{1}, \psi_{1})$, satisfies $(\phi_{1}, \psi_{1})\in \mathds{E}^{++}$.
\begin{proof}
For any given $\textbf{u}(x,t)=(u_{1}(x,t),u_{2}(x,t) )\in \mathds{E}$, we claim that the linear problem
\begin{eqnarray}\label{3-4}
\left\{
\begin{array}{ll}
\zeta_{t}=-\delta_{1} \zeta-\kappa\zeta+u_{1} ,\; &\,x\in(l_{1},l_{2}), t\in(0, \tau],  \\[2mm]
\eta_{t}=d_{2}\Delta \eta-\delta_{2} \eta-\kappa\eta+u_{2},\; &\,x\in(l_{1},l_{2}), t\in(0, \tau],  \\[2mm]
\zeta_{t}=d_{1}\Delta \zeta-a_{11}\zeta+a_{12}\eta-\kappa\zeta+u_{1} ,\; &\,x\in(l_{1},l_{2}), t\in(\tau, T],  \\[2mm]
\eta_{t}=d_{2}\Delta \eta-a_{22}\eta+f'(0) \zeta-\kappa\eta+u_{2}  ,\; &\,x\in(l_{1},l_{2}), t\in(\tau, T],  \\[2mm]
\zeta(x,t)=\eta(x,t)=0,\; &\, x\in\{l_{1},l_{2}\},  t\in(0, T],\\[2mm]
\zeta(x,0)=H'(0)\zeta(x,T),\; &\, x\in[l_{1},l_{2}],  \\[2mm]
\eta(x,0)=\eta(x,T),\; &\, x\in[l_{1},l_{2}]
\end{array} \right.
\end{eqnarray}
admits a unique solution $\textbf{v}(x,t)=(v_{1}(x,t),v_{2}(x,t) )\in \mathds{H}^{1+\nu,\frac{1+\nu}{2}}\cap\mathds{E}$ with any $0<\nu<1$, where
\begin{equation*}
\kappa=a_{12}\max\{H'(0), 1\}+f'(0)+1+\ln[\max(1,1/H'(0))].
\end{equation*}
Indeed, we define the space
\begin{equation*}
\mathcal{BC}=\big\{u_{i}\in \mathcal{C}[l_{1},l_{2}];~ |u_{i}(x)|\leq C~\forall x\in [l_{1},l_{2}]\text{~and~}u_{i}(l_{j})=0~\forall i,j=1,2\big\},
\end{equation*}
where $C=\max\Big\{\max\limits_{\overline{Q}^{T}_{0}}u_{1}(x,t),\max\limits_{\overline{Q}^{T}_{0}}u_{2}(x,t)\Big\}$. Clearly, the space $\mathcal{BC}$ is a
bounded, closed, and convex set. For any given $(\zeta_{0}(x), \eta_{0}(x))\in \mathcal{BC}$, we define the
operator $\mathcal{B}(\zeta_{0}(x), \eta_{0}(x))=(H'(0)\zeta(T, x),\eta(T, x))$, where $(\zeta(x,t), \eta(x,t))$ is the unique solution to the following initial-boundary
value problem
\begin{eqnarray*}
\left\{
\begin{array}{ll}
\zeta_{t}=-\delta_{1} \zeta-\kappa\zeta+u_{1} ,\; &\,x\in(l_{1},l_{2}), t\in(0, \tau],  \\[2mm]
\eta_{t}=d_{2}\Delta \eta-\delta_{2} \eta-\kappa\eta+u_{2},\; &\,x\in(l_{1},l_{2}), t\in(0, \tau],  \\[2mm]
\zeta_{t}=d_{1}\Delta \zeta-a_{11}\zeta+a_{12}\eta-\kappa\zeta+u_{1} ,\; &\,x\in(l_{1},l_{2}), t\in(\tau, T],  \\[2mm]
\eta_{t}=d_{2}\Delta \eta-a_{22}\eta+f'(0) \zeta-\kappa\eta+u_{2}  ,\; &\,x\in(l_{1},l_{2}), t\in(\tau, T],  \\[2mm]
\zeta(x,t)=\eta(x,t)=0,\; &\, x\in\{l_{1},l_{2}\},  t\in(0, T],\\[2mm]
\zeta(x,0)=\zeta_{0}(x),~\eta(x,0)=\eta_{0}(x),\; &\, x\in[l_{1},l_{2}].
\end{array} \right.
\end{eqnarray*}
By using parabolic regularity theory, it is easily seen that
\begin{equation*}
(\zeta(x,t), \eta(x,t))\in \mathcal{C}^{0,1}(\overline{Q}^{\tau}_{0}, \mathds{R}^{2})\cap\mathcal{C}^{1+\nu,\frac{1+\nu}{2}}(\overline{Q}^{T}_{\tau}, \mathds{R}^{2}),
\end{equation*}
where $\nu\in(0,1)$. This implies that $\mathcal{B}(\zeta_{0}(x), \eta_{0}(x))\in \mathcal{C}^{1+\nu}([l_{1},l_{2}], \mathds{R}^{2})$. By using Arzela-Ascoli theorem and Sobolev embedding theorem, we have that operator $\mathcal{B}$ is compact. Additionally, the comparison principle for cooperative systems yields that $\mathcal{B}$ is a continuous self-mapping in $\mathcal{BC}$. Therefore, Schauder fixed point theorem show that $\mathcal{B}$ admits at least one fixed point. Specifically, problem \eqref{3-4} has at least one solution $(\zeta, \eta)\in \mathds{H}^{1+\nu,\frac{1+\nu}{2}}$, where $\nu\in(0,1)$.

Next, we will prove that the solution of problem \eqref{3-4} is unique. In order to prove this, we only need to prove that there is only zero solution to the problem
\begin{eqnarray}\label{3-5-ZQ}
\left\{
\begin{array}{ll}
\zeta_{t}=-\delta_{1} \zeta-\kappa\zeta,\; &\,x\in(l_{1},l_{2}), t\in(0, \tau],  \\[2mm]
\eta_{t}=d_{2}\Delta \eta-\delta_{2} \eta-\kappa\eta,\; &\,x\in(l_{1},l_{2}), t\in(0, \tau],  \\[2mm]
\zeta_{t}=d_{1}\Delta \zeta-a_{11}\zeta+a_{12}\eta-\kappa\zeta,\; &\,x\in(l_{1},l_{2}), t\in(\tau, T],  \\[2mm]
\eta_{t}=d_{2}\Delta \eta-a_{22}\eta+f'(0) \zeta-\kappa\eta,\; &\,x\in(l_{1},l_{2}), t\in(\tau, T],  \\[2mm]
\zeta(x,t)=\eta(x,t)=0,\; &\, x\in\{l_{1},l_{2}\},  t\in(0, T],\\[2mm]
\zeta(x,0)=H'(0)\zeta(x,T),\; &\, x\in[l_{1},l_{2}],  \\[2mm]
\eta(x,0)=\eta(x,T),\; &\, x\in[l_{1},l_{2}].
\end{array} \right.
\end{eqnarray}
Otherwise, we assume that
\begin{equation*}
\min\bigg\{\min\limits_{[l_{1},l_{2}]\times [0,T] }\zeta(x,t), \min\limits_{[l_{1},l_{2}]\times [0,T]  }\eta(x,t)\bigg\}<0
\end{equation*}
or
\begin{equation*}
\max\bigg\{\max\limits_{[l_{1},l_{2}]\times [0,T] }\zeta(x,t), \max\limits_{[l_{1},l_{2}]\times [0,T]  }\eta(x,t)\bigg\}>0,
\end{equation*}
and show a contradiction. We assume, without loss of generality, that there exists $(x_{0}, t_{0})\in (l_{1},l_{2})\times [0,T]$ such that
\begin{equation*}
\max\bigg\{\max\limits_{[l_{1},l_{2}]\times [0,T] }\zeta(x,t), \max\limits_{[l_{1},l_{2}]\times [0,T]  }\eta(x,t)\bigg\}=\zeta(x_{0}, t_{0})>0.
\end{equation*}
When $(x_{0}, t_{0})\in (l_{1},l_{2})\times (0,\tau]$, we have that $\zeta_{t}|_{(x,t)=(x_{0}, t_{0})}=0$ and $-\delta_{1} \zeta(x_{0}, t_{0})-\kappa\zeta(x_{0}, t_{0})<0$.
This contradicts the first equation in problem \eqref{3-5-ZQ}. When $(x_{0}, t_{0})\in (l_{1},l_{2})\times (\tau,T]$, we have that
\begin{equation*}
\zeta_{t}|_{(x,t)=(x_{0}, t_{0})}-d_{1}\Delta \zeta|_{(x,t)=(x_{0}, t_{0})}\geq 0.
\end{equation*}
On the other hand, one can obtain that
\begin{equation*}
-a_{11}\zeta(x_{0}, t_{0})+a_{12}\eta(x_{0}, t_{0})-\kappa\zeta(x_{0}, t_{0})<0.
\end{equation*}
This contradicts the third equation in problem \eqref{3-5-ZQ}. When $(x_{0}, t_{0})\in (l_{1},l_{2})\times \{0\}$, it follows from the sixth equation in problem \eqref{3-5-ZQ} that $\zeta(x_{0},0)\leq \zeta(x_{0},T)$. If $\zeta(x_{0},0)< \zeta(x_{0},T)$, this contradicts with $(x_{0}, t_{0})\in (l_{1},l_{2})\times \{0\}$. If $\zeta(x_{0},0)= \zeta(x_{0},T)$, the case of $(x_{0}, t_{0})\in (l_{1},l_{2})\times (\tau,T]$ shows that this is not possible. Hence, the solution of problem \eqref{3-5-ZQ} is unique.

Now, we define the operator $\mathcal{A}\textbf{u}=\textbf{v}$. Since $\mathds{H}^{1+\nu,\frac{1+\nu}{2}}$ is compactly embedded in $\mathds{H}^{1,0}$, $\mathcal{A}$ is
a linear compact operator. Moreover, it follows from the strong maximum principle and Hopf boundary lemma for cooperative systems that $\mathcal{A}$ is a strongly
order-preserving with respect to $\mathds{E}$. Finally, a strong version of the Krein-Rutman theorem (see, e.g., \cite[Theorem 19.3]{DeimIing}) yields that the desired
conclusion stands. More precisely, the principal eigenvalue $\lambda_{1}$ of problem \eqref{3-3} may be expressed as follows
\begin{equation*}
\lambda_{1}=\frac{1}{r(\mathcal{A})}-\kappa,
\end{equation*}
where $r(\mathcal{A})$ denotes the spectral radius of $\mathcal{A}$. This ends the proof.
\end{proof}
\end{theorem}
\autoref{theorem 3.1} above establishes the existence and uniqueness of the principal eigenvalue to problem \eqref{3-3}. Additionally, its corresponding eigenfunction is algebraically simple and strongly positive. It is well known that the eigenfunction plays an important role in constructing the upper and lower solutions of the model. Next,
we investigate the property of the eigenfunction further using the method of separation of variables. Let
\begin{equation}\label{3-46}
\phi_{1}(x,t)=\chi_{1}(x)\Phi(t) \text{~and~}  \psi_{1}(x,t)=\chi_{1}(x)\Psi(t)
\end{equation}
for $x\in[l_{1}, l_{2}]$ and $t\in[0, T]$, where $\chi_{1}(x)$ is the eigenfunction corresponding to the principal eigenvalue $\kappa_{1}$ of problem
 \begin{eqnarray}\label{3-30}
\left\{
\begin{array}{ll}
-\Delta \chi(x)=\kappa  \chi(x)\text{~~in~~} (l_{1}, l_{2}),  \\[2mm]
~~~~~~ \chi(x)=0  \text{~~~~~~~on~~} \{l_{1}, l_{2}\},
\end{array} \right.
\end{eqnarray}
and $(\Phi(t), \Psi(t))$ satisfies the following problem
\begin{eqnarray}\label{3-15}
\left\{
\begin{array}{ll}
\Phi'(t)=[\lambda-\delta_{1} ]\Phi(t),\; &\,  t\in(0, \tau], \\[2mm]
\Psi'(t)=[\lambda-\delta_{2}-d_{2}\kappa_{1} ]\Psi(t),\; &\,  t\in(0, \tau], \\[2mm]
\Phi'(t)=[\lambda-d_{1}\kappa_{1} -a_{11}]\Phi(t)+a_{12}\Psi(t),\; &\,  t\in(\tau, T], \\[2mm]
\Psi'(t)=[\lambda-d_{2}\kappa_{1} -a_{22}]\Psi(t)+f'(0)\Phi(t),\; &\,  t\in(\tau, T], \\[2mm]
\Phi(0)=H'(0)\Phi(T), ~~~~~~\Psi(0)=\Psi(T).
\end{array} \right.
\end{eqnarray}
Specifically, when $l_{2}=-l_{1}=l$ we have that $\kappa_{1}=\frac{\pi^{2}}{4l^{2}}$ and $\chi_{1}(x)=\cos(\frac{\pi}{2l}x)$ for $x\in[-l, l]$. By observing problems \eqref{3-30} and \eqref{3-15}, it can be seen that $\kappa_{1}$, $\chi_{1}(x)$, $\Phi(t)$, and $\Psi(t)$ depend on $l$. Therefore, we also write $\kappa_{1}=\kappa^{l}_{1}$, $\chi_{1}(x)=\chi^{l}_{1}(x)$, $\Phi(t)=\Phi_{l}(t)$, and $\Psi(t)=\Psi_{l}(t)$ in order to reflect this dependency.
\begin{lemma}\label{lemma 3.1.4}
Let $l_{2}=-l_{1}=l$ and $l\geq s_{0}$. The eigenvalue problem \eqref{3-15} has a unique principal eigenvalue $\lambda_{1}$ with an eigenfunction
$(\Phi(t),\Psi(t))\succ\succ \textbf{0}$ for $t\in[0,T]$. Moreover, there exist positive constants $n_{1}$, $n_{2}$, $n_{3}$, and $n_{4}$, independent of $l$, such that
\begin{eqnarray*}
\begin{array}{l}
n_{1}\leq \Phi_{l}(0), ~~\Phi_{l}(t)\leq n_{2}, ~~~~t\in[0,T], \\[2mm]
n_{3}\leq \Psi_{l}(0), ~~\Psi_{l}(t)\leq n_{4},~~~~t\in[0,T].
\end{array}
\end{eqnarray*}
\begin{proof}
We first prove the first assertion. The third and fourth equations of \eqref{3-15} are abbreviated as
\begin{equation*}
\Bigg(
\begin{matrix}
   \Phi'(t)\\
   \Psi'(t)\\
  \end{matrix}
\Bigg)=\Bigg(
\begin{matrix}
   \lambda-d_{1}\kappa_{1}-a_{11}&a_{12}\\
   f'(0)&   \lambda-d_{2}\kappa_{1}-a_{22}\\
  \end{matrix}
\Bigg)
\Bigg(
\begin{matrix}
   \Phi(t) \\
   \Psi(t) \\
  \end{matrix}
\Bigg):=\mathbf{A}
\Bigg(
\begin{matrix}
   \Phi(t) \\
   \Psi(t) \\
  \end{matrix}
\Bigg)
\end{equation*}
for any $t \in (\tau,T]$. Then, we see from the characteristic equation $|\mathbf{A}-\mu \mathbf{E}|=0$ that
\begin{equation*}
\mu_{1,2}=\frac{2\lambda-(d_{1}+d_{2})\kappa_{1}-a_{11}-a_{22}\pm\sqrt{[a_{22}+(d_{2}-d_{1})\kappa_{1}-a_{11}]^{2}+4a_{12}f'(0)}}{2}:=\lambda+c_{1,2}.
\end{equation*}
Direct calculation yields that
\begin{eqnarray}\label{3-16}
\begin{array}{l}
\min\Big\{2\sqrt{a_{12}f'(0)}, |a_{22}+(d_{2}-d_{1})\kappa_{1}-a_{11}|\Big\}\leq c_{1}-c_{2}\\
~~~~~~~~~~~~~~~~~~~~~~~~~~~~~~~~~~~~~~~~~~~~~~~~~~~~~~~~~~~~~~~~~~~\leq 2\sqrt{a_{12}f'(0)}+|a_{22}+(d_{2}-d_{1})\kappa_{1}-a_{11}|, \\[2mm]
\min\bigg\{\frac{\sqrt{2a_{12}f'(0)}}{2}, \frac{\sqrt{a_{12}f'(0)}}{2\sqrt{a_{22}+(d_{2}+d_{1})\kappa^{s_{0}}_{1}+a_{11}}}\bigg\}\leq a_{11}+d_{1}\kappa_{1}+c_{1}=-\big(a_{22}+d_{2}\kappa_{1}+c_{2}\big)\\
~~~~~~~~~~~~~~~~~~~~~~~~~~~~~~~~~~~~~~~~~~~~~~~~~~~~~~~~~~\leq a_{11}+a_{22}+(d_{2}+d_{1})\kappa_{1}+\sqrt{a_{12}f'(0)}.
\end{array}
\end{eqnarray}
The linearly independent eigenvectors $(k_{11}, k_{12})$ and $(k_{21}, k_{22})$ associated with eigenvalues $\mu_{1}$ and $\mu_{2}$ satisfy
\begin{equation*}
\big(
\begin{matrix}
   k_{11}&k_{12}\\
  \end{matrix}
\big)
\Bigg(
\begin{matrix}
   \lambda-d_{1}\kappa_{1}-a_{11}-\mu_{1}\\
   f'(0) \\
  \end{matrix}
\Bigg)=0\text{~and~}
\big(
\begin{matrix}
   k_{21}&k_{22}\\
  \end{matrix}
\big)
\Bigg(
\begin{matrix}
  a_{12}\\
  \lambda-d_{2}\kappa_{1}-a_{22}-\mu_{2}\\
  \end{matrix}
\Bigg)=0,
\end{equation*}
respectively. Then, it is easily seen that
\begin{equation*}
(k_{11}, k_{12})=\big(f'(0), a_{11}+d_{1}\kappa_{1}-\lambda+\mu_{1}\big)=\big(f'(0), a_{11}+d_{1}\kappa_{1}+c_{1}\big)
\end{equation*}
and
\begin{equation*}
(k_{21}, k_{22})=\big(a_{22}+d_{2}\kappa_{1}-\lambda+\mu_{2}, a_{12}\big)=\big(a_{22}+d_{2}\kappa_{1}+c_{2}, a_{12}\big).
\end{equation*}
In the following, we consider the corresponding algebraic equations
\begin{equation*}
\Bigg(
\begin{matrix}
   e^{\mu_{1}t} \\
   ke^{\mu_{2}t} \\
  \end{matrix}
\Bigg)=
\Bigg(
\begin{matrix}
   k_{11}&k_{12}\\
   k_{21}&k_{22}\\
  \end{matrix}
\Bigg)
\Bigg(
\begin{matrix}
   \Phi(t)\\
   \Psi(t) \\
  \end{matrix}
\Bigg):=\mathbf{B}
\Bigg(
\begin{matrix}
   \Phi(t)\\
   \Psi(t) \\
  \end{matrix}
\Bigg),
\end{equation*}
and a simple calculation yields that
\begin{equation}\label{3-17}
(\Phi(t), \Psi(t))=\Bigg(\frac{a_{12}e^{\mu_{1}t}-\big(a_{11}+d_{1}\kappa_{1}+c_{1}\big)ke^{\mu_{2}t}}{|\mathbf{B}|},
\frac{f'(0)ke^{\mu_{2}t}-\big(a_{22}+d_{2}\kappa_{1}+c_{2}\big)e^{\mu_{1}t}}{|\mathbf{B}|}\Bigg)
\end{equation}
for $t\in[\tau, T]$, where
\begin{equation*}
\begin{aligned}
|\mathbf{B}|=&a_{12}f'(0)-\big(a_{22}+d_{2}\kappa_{1}+c_{2}\big)\big(a_{11}+d_{1}\kappa_{1}+c_{1}\big)\\
=&a_{12}f'(0)+\big(a_{22}+d_{2}\kappa_{1}+c_{2}\big)^{2}>0.
\end{aligned}
\end{equation*}
Using \eqref{3-17} and the first two equations of problem \eqref{3-15}, we have that
\begin{equation}\label{3-18}
\begin{cases}
\begin{aligned}
&a_{12}e^{\mu_{1}\tau}-\big(a_{11}+d_{1}\kappa_{1}+c_{1}\big)ke^{\mu_{2}\tau}=H'(0)e^{(\lambda-\delta_{1})\tau}\big[a_{12}e^{\mu_{1}T}-\big(a_{11}
+d_{1}\kappa_{1}+c_{1}\big)ke^{\mu_{2}T}\big],\\
&f'(0)ke^{\mu_{2}\tau}-\big(a_{22}+d_{2}\kappa_{1}+c_{2}\big)e^{\mu_{1}\tau}=e^{(\lambda-\delta_{2}-d_{2}\kappa_{1})\tau}\big[f'(0)ke^{\mu_{2}T}
-\big(a_{22}+d_{2}\kappa_{1}+c_{2}\big)e^{\mu_{1}T}\big].
\end{aligned}
\end{cases}
\end{equation}
For convenience, we denote
\begin{equation*}
\begin{aligned}
&y=e^{\lambda T}, ~~~~~~~~b_{11}=a_{12}e^{c_{1}\tau},  ~~~~~~~~~b_{12}=\big(a_{11}+d_{1}\kappa_{1}+c_{1}\big)e^{c_{2}\tau},~~b_{13}=a_{12}e^{c_{1}T},\\
& b_{14}=\big(a_{11}+d_{1}\kappa_{1}+c_{1}\big)e^{c_{2}T}, ~~~~~~~~b_{22}=\big(a_{22}+d_{2}\kappa_{1}+c_{2}\big)e^{c_{1}\tau}, ~~b_{21}=f'(0)e^{c_{2}\tau},\\
&b_{23}=f'(0)e^{c_{2}T},\vartheta_{1}=H'(0)e^{-\delta_{1}\tau},b_{24}=\big(a_{22}+d_{2}\kappa_{1}+c_{2}\big)e^{c_{1}T},~  ~\vartheta_{2}=e^{(-\delta_{2}-d_{2}\kappa_{1})\tau},
\end{aligned}
\end{equation*}
and then \eqref{3-18} can be rewritten as
\begin{equation}\label{3-19}
\begin{cases}
b_{11}-b_{12}k=\vartheta_{1}(b_{13}-b_{14}k)y,\\
b_{21}k-b_{22}=\vartheta_{2}(b_{23}k-b_{24})y.
\end{cases}
\end{equation}
When $k=\frac{b_{13}}{b_{14}}$ or $\frac{b_{24}}{b_{23}}$, we have that $\Phi(t),\Psi(t)>0$ does not hold for all $t\in[0,T]$. Therefore, \eqref{3-19} can be further written as
\begin{equation}\label{3-20}
\begin{cases}
y=(b_{11}-b_{12}k)/[\vartheta_{1}(b_{13}-b_{14}k)],\\
y=(b_{21}k-b_{22})/[\vartheta_{2}(b_{23}k-b_{24})].
\end{cases}
\end{equation}
Next, we show that \eqref{3-20} admits a unique solution $(k_{0}, y_{0})$ such that $\Phi(t),\Psi(t)>0$ for all $t\in[0,T]$. For this purpose, we need to consider the following three cases.
\begin{figure}[!htb]
\centering
\subfigure[$\frac{b_{12}}{\vartheta_{1}b_{14}}=\frac{b_{21}}{\vartheta_{2}b_{23}}$]{ {
\includegraphics[width=0.47\textwidth]{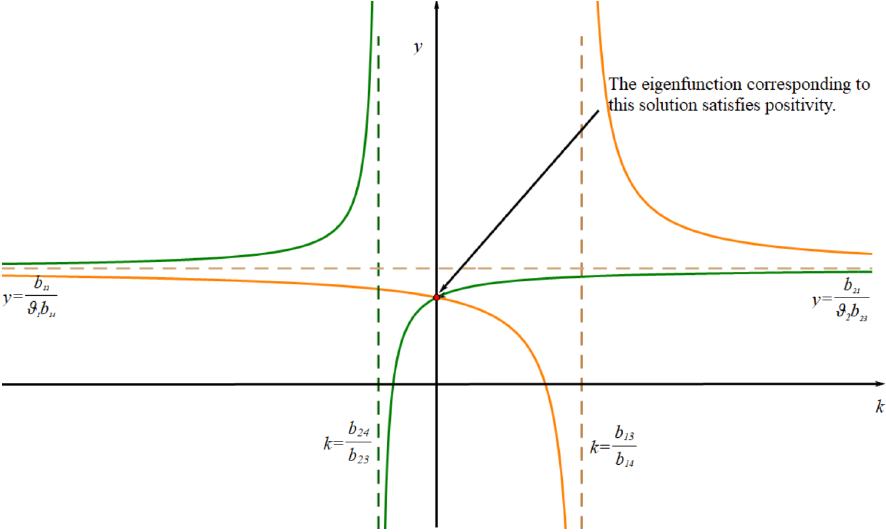}
} }\\
\subfigure[$\frac{b_{12}}{\vartheta_{1}b_{14}}>\frac{b_{21}}{\vartheta_{2}b_{23}}$]{ {
\includegraphics[width=0.47\textwidth]{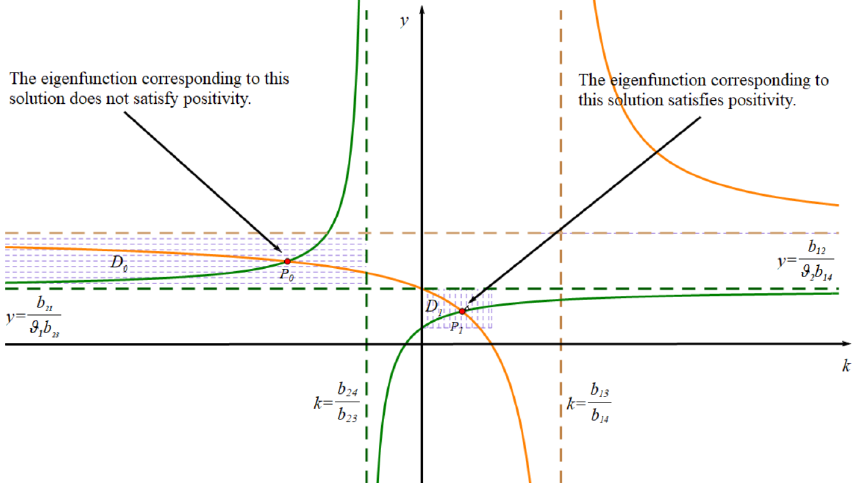}
} }
\subfigure[$\frac{b_{12}}{\vartheta_{1}b_{14}}<\frac{b_{21}}{\vartheta_{2}b_{23}}$]{ {
\includegraphics[width=0.47\textwidth]{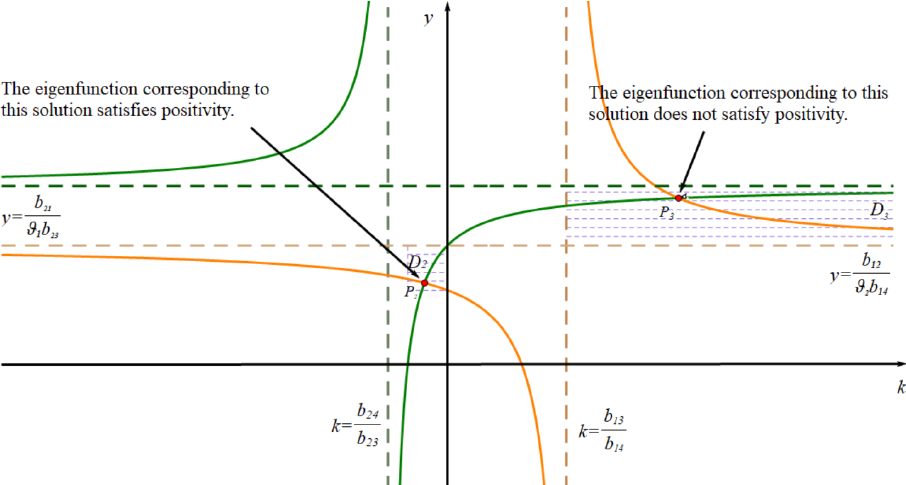}
} }
\caption{Distribution graph of the unique solution of \eqref{3-20}}
\label{A}
\end{figure}

\textbf{Case 1}. If $\frac{b_{12}}{\vartheta_{1}b_{14}}=\frac{b_{21}}{\vartheta_{2}b_{23}}$, it follows from \autoref{A}\textcolor[rgb]{0.00,0.00,1.00}{(a)} that \eqref{3-20} admits a unique solution
\begin{equation*}
p=:\big(0, \frac{b_{11}}{\vartheta_{1}b_{13}}\big)=\big(0, \frac{b_{22}}{\vartheta_{2}b_{24}}\big).
\end{equation*}
Substituting the coordinate of point $p$ into \eqref{3-17}, we obtain that $\Phi(t),\Psi(t)>0$ for all $t\in[0,T]$. Additionally, the principal eigenvalue $\lambda_{1}$ can be expressed explicitly as follows
\begin{equation}\label{3-21}
\lambda_{1}=\frac{c_{1}(\tau-T)+(\delta_{2}+d_{2}\kappa_{1})\tau}{T}.
\end{equation}

\textbf{Case 2}. If $\frac{b_{12}}{\vartheta_{1}b_{14}}>\frac{b_{21}}{\vartheta_{2}b_{23}}$, one can obtain from \autoref{A}\textcolor[rgb]{0.00,0.00,1.00}{(b)} that \eqref{3-20} exists two solutions $p_{0}$ and $p_{1}$, and
\begin{equation*}
\begin{split}
&p_{0}\in D_{0}= \big(-\infty, \frac{b_{24}}{b_{23}}\big)\times \big(\frac{b_{21}}{\vartheta_{2}b_{23}}, \frac{b_{12}}{\vartheta_{1}b_{14}}\big),\\
&~~~p_{1}\in D_{1}=\big(0, \frac{b_{11}}{b_{12}}\big)\times \big(\frac{b_{22}}{\vartheta_{2}b_{24}}, \frac{b_{11}}{\vartheta_{1}b_{13}}\big).
\end{split}
\end{equation*}
In region $D_{1}$, we have from problem \eqref{3-15} and \eqref{3-17} that $\Phi(t),\Psi(t)>0$ for all $t\in[0,T]$. However, it is not true that $\Psi(T)>0$ in region $D_{0}$. Therefore, \eqref{3-20} admits a unique solution such that $\Phi(t),\Psi(t)>0$ for all $t\in[0,T]$.

\textbf{Case 3}. If $\frac{b_{12}}{\vartheta_{1}b_{14}}<\frac{b_{21}}{\vartheta_{2}b_{23}}$, \autoref{A}\textcolor[rgb]{0.00,0.00,1.00}{(c)} indicates that \eqref{3-20} has two solutions $p_{2}$ and $p_{3}$, and
\begin{equation*}
\begin{split}
&~p_{2}\in D_{2}= \big(\frac{b_{22}}{b_{21}}, 0\big)\times \big(\frac{b_{11}}{\vartheta_{1}b_{13}},\frac{b_{22}}{\vartheta_{2}b_{24}}\big),\\
&p_{3}\in D_{3}=\big(\frac{b_{13}}{b_{14}},+\infty\big)\times \big(\frac{b_{12}}{\vartheta_{1}b_{14}}, \frac{b_{21}}{\vartheta_{2}b_{23}}\big).
\end{split}
\end{equation*}
A discussion similar to that in \textbf{Case 2} yields that \eqref{3-20} admits a unique solution such that $\Phi(t),\Psi(t)>0$ for all $t\in[0,T]$.
This completes the proof of the first assertion.

Next, we prove the second assertion. We first  show that there exists a positive constant $n_{1}$, independent of $l$, such that  $n_{1}\leq \Phi_{l}(0)$.
It follows from \eqref{3-15} and \eqref{3-17} that
\begin{equation}\label{3-45}
\Phi_{l}(0)=e^{(\delta_{1}- \lambda)\tau}\Phi_{l}(\tau)=\frac{e^{\delta_{1}\tau}[a_{12}e^{c_{1}\tau}-\big(a_{11}+d_{1}\kappa_{1}+c_{1}\big)ke^{c_{2}\tau}]}{|\mathbf{B}|}.
\end{equation}
For \textbf{Cases 1} and \textbf{3}, the result is straightforward; we therefore focus solely on proving \textbf{Case 2}. Let
\begin{equation*}
x_{0}=:e^{\delta_{2}\tau+[\tau-T]\sqrt{[a_{22}+(d_{2}+d_{1})\kappa^{s_{0}}_{1}+a_{11}]^{2}+4a_{12}f'(0)}}.
\end{equation*}
We claim that there exists a sufficiently small positive constant
\begin{equation*}
\epsilon_{0}=\min\Big\{a_{12}e^{-(a_{11}+a_{22}+d_{1}\kappa^{s_{0}}_{1}+d_{2}\kappa^{s_{0}}_{1})\tau},x_{0}a_{12}H'(0)e^{-\delta_{1}\tau}e^{-(a_{11}+a_{22}\\
+d_{1}\kappa^{s_{0}}_{1}+d_{2}\kappa^{s_{0}}_{1})T}[1-e^{-2\sqrt{a_{12}f'(0)}(T-\tau)}]\Big\}
\end{equation*}
such that
\begin{equation*}
0<k_{0}=:\frac{b_{11}-\epsilon_{0}}{b_{12}},
\end{equation*}
where $k_{0}$ denotes the $k$-coordinate of point $p_{1}$. By regular calculations, we have that
\begin{equation*}
f(\epsilon)=:\frac{b_{21}(b_{11}-\epsilon)-b_{12}b_{22}}{\vartheta_{2}[b_{23}(b_{11}-\epsilon)-b_{12}b_{24}]}
\end{equation*}
is strictly monotonically decreasing with respect to $\epsilon$. Then, it follows that
\begin{equation*}
f(\epsilon_{0})\geq f(b_{11})=\frac{b_{22}}{\vartheta_{2}b_{24}}\geq e^{\delta_{2}\tau+[\tau-T]\sqrt{[a_{22}+(d_{2}+d_{1})\kappa^{s_{0}}_{1}+a_{11}]^{2}+4a_{12}f'(0)}}=x_{0}.
\end{equation*}
Therefore, one can obtain that
\begin{equation*}
\frac{b_{21}k_{0}-b_{22}}{\vartheta_{2}(b_{23}k_{0}-b_{24})}-\frac{b_{11}-b_{12}k_{0}}{\vartheta_{1}(b_{13}-b_{14}k_{0})}\geq x_{0}-\frac{b_{12}\epsilon_{0}}{\vartheta_{1}(b_{12}b_{13}-b_{14}b_{11}}\geq 0,
\end{equation*}
which implies that the claim holds. With the help of this claim, it can be immediately concluded from \eqref{3-16} and \eqref{3-45} that
\begin{equation*}
\Phi_{l}(0)\geq \frac{\epsilon_{0}e^{\delta_{1}\tau}}{a_{12}f'(0)+\Big(a_{11}+a_{22}+(d_{2}+d_{1})\kappa^{s_{0}}_{1}+\sqrt{a_{12}f'(0)}\Big)^{2}}:=n_{1}.
\end{equation*}
A similar method can be used to prove that there exists a positive constant $n_{3}$, independent of $l$, such that $\Psi_{l}(0)\geq n_{3}$. To conserve space, we omit its proof here.

Finally, we show the remaining part. For \textbf{Case 1}, it can be immediately obtained from \eqref{3-17} and \autoref{lemma 3.1.3}\textcolor{blue}{(2)} that
\begin{equation*}
\Phi_{l}(t)\leq \frac{a_{12}e^{\mu_{1}t}}{|\mathbf{B}|}\leq \frac{e^{\mu_{1}t}}{f'(0)}\leq \frac{e^{\big[\max\{\lambda_{1}(s_{0}),1\}+\sqrt{[a_{22}+(d_{2}+d_{1})\kappa^{s_{0}}_{1}+a_{11}]^{2}+4a_{12}f'(0)}\big]T}}     {f'(0)}:=n_{21}
\end{equation*}
for $t\in(\tau, T]$. Then, it follows from the first equation of problem \eqref{3-15} that
\begin{equation*}
\Phi_{l}(t)=\Phi_{l}(0)e^{[\lambda_{1}(l)-\delta_{1}]t}\leq H'(0)n_{21} e^{\max\{\lambda_{1}(s_{0}),1\}\tau}:=n_{22}
\end{equation*}
for $t\in(0, \tau]$. Similarly, when \textbf{Case 2} happens, one can find positive constants $n_{23}$ and $n_{24}$, independent of $l$, such that $\Phi_{l}(t)\leq n_{23}$ for $t\in(\tau, T]$ and $\Phi_{l}(t)\leq n_{24}$ for $t\in(0, \tau]$. For \textbf{Case 3}, we have from \eqref{3-16} and \eqref{3-17} that
\begin{equation*}
\begin{split}
\Phi_{l}(t)\leq& \frac{f'(0)a_{12}e^{\mu_{1}t}+ \big(a_{11}+d_{1}\kappa_{1}+c_{1}\big)^{2}e^{\mu_{2}t+c_{1}\tau-c_{2}\tau}}{f'(0)|\mathbf{B}|}   \leq \frac{e^{\mu_{1}t}[a_{12}f'(0)+(a_{11}+d_{1}\kappa_{1}+c_{1})^{2}]}{a_{12}[f'(0)]^{2}}  \\
\leq& \frac{e^{\mu_{1}t}\big[a_{12}f'(0)+\big(a_{11}+a_{22}+(d_{2}+d_{1})\kappa^{s_{0}}_{1}+\sqrt{a_{12}f'(0)}\big)^{2}\big]}{a_{12}[f'(0)]^{2}} :=n_{25}
\end{split}
\end{equation*}
for $t\in(\tau, T]$. Then, it follows from the first equation of problem \eqref{3-15} that
\begin{equation*}
\Phi_{l}(t)=\Phi_{l}(0)e^{[\lambda_{1}(l)-\delta_{1}]t}\leq H'(0) n_{25}e^{\max\{\lambda_{1}(s_{0}),1\}\tau}:=n_{26}
\end{equation*}
for $t\in(0, \tau]$. By taking
\begin{equation*}
n_{2}=\max\{n_{21},n_{22},n_{23},n_{24},n_{25},n_{26} \},
\end{equation*}
we can prove the conclusion that there exists a positive constant $n_{2}$, independent of $l$, such that $\Phi_{l}(t)\leq n_{2}$ for $t\in[0, T]$. By a case analysis similar to that used in the preceding result, one can obtain that $\Psi_{l}(t)\leq n_{4}$ for $t\in[0, T]$, where $n_{4}$ is a positive constant independent of $l$. This proof is completed.
\end{proof}
\end{lemma}
Since this eigenvalue cannot be expressed explicitly, it becomes challenging to determine its sign. To address this issue, we shall present two lemmas that provide estimates for this principal eigenvalue. Before proceeding, we first introduce the adjoint problem corresponding to problem \eqref{3-3}, which has the form
\begin{eqnarray}\label{3-6}
\left\{
\begin{array}{ll}
-\zeta_{t}=-\delta_{1} \zeta+\mu \zeta ,\; &\,x\in(l_{1},l_{2}), t\in(0, \tau],  \\[2mm]
-\eta_{t}=d_{2}\Delta \eta-\delta_{2} \eta+\mu \eta ,\; &\,x\in(l_{1},l_{2}), t\in(0, \tau],  \\[2mm]
-\zeta_{t}=d_{1}\Delta \zeta-a_{11}\zeta+f'(0)\eta+\mu \zeta ,\; &\,x\in(l_{1},l_{2}), t\in(\tau, T],  \\[2mm]
-\eta_{t}=d_{2}\Delta \eta-a_{22}\eta+ a_{12}\zeta+\mu \eta  ,\; &\,x\in(l_{1},l_{2}), t\in(\tau, T],  \\[2mm]
\zeta(x,t)=\eta(x,t)=0,\; &\, x\in\{l_{1},l_{2}\},  t\in[0, T],\\[2mm]
\zeta(x,0)=\frac{1}{H'(0)}\zeta(x,T),\; &\, x\in[l_{1},l_{2}],  \\[2mm]
\eta(x,0)=\eta(x,T),\; &\, x\in[l_{1},l_{2}].
\end{array} \right.
\end{eqnarray}
Using a method similar to the proof of \autoref{theorem 3.1}, we can conclude that the adjoint problem\eqref{3-6} has a unique algebraically simple principal eigenvalue, denoted by $\mu_{1}$, to a positive vector-valued eigenfunction. Additionally, the positive vector-valued eigenfunction associated with $\mu_{1}$, denoted by $(\zeta_{1}, \eta_{1})$, satisfies $(\zeta_{1}, \eta_{1})\in \mathds{E}^{++}$. In fact, the principal eigenvalue $\lambda_{1}$ of the problem \eqref{3-3} and the principal
eigenvalue $\mu_{1}$ of the problem \eqref{3-6} are equal. This conclusion can be proven by a method similar to \cite[Theorem 5.1]{Zhou-Lin-Santos-CFS}.
\begin{lemma}\label{lemma 3.1.1}
If there exists a vector-valued eigenfunction $(\phi, \psi)\in \mathds{E}^{+}-\{\textbf{0}\}$ and a number $\overline{\lambda}$ satisfying
\begin{eqnarray}\label{3-7}
\left\{
\begin{array}{ll}
\phi_{t}\leq-\delta_{1} \phi+\overline{\lambda} \phi ,\; &\,x\in(l_{1},l_{2}), t\in(0, \tau],  \\[2mm]
\psi_{t}\leq d_{2}\Delta \psi-\delta_{2} \psi+\overline{\lambda} \psi ,\; &\,x\in(l_{1},l_{2}), t\in(0, \tau],  \\[2mm]
\phi_{t}\leq d_{1}\Delta \phi-a_{11}\phi+a_{12}\psi+\overline{\lambda} \phi ,\; &\,x\in(l_{1},l_{2}), t\in(\tau, T],  \\[2mm]
\psi_{t}\leq d_{2}\Delta \psi-a_{22}\psi+f'(0) \phi+\overline{\lambda} \psi  ,\; &\,x\in(l_{1},l_{2}), t\in(\tau, T],  \\[2mm]
\phi(x,t)=\psi(x,t)=0,\; &\, x\in\{l_{1},l_{2}\},  t\in(0, T],\\[2mm]
\phi(x,0)\leq H'(0)\phi(x,T),\; &\, x\in[l_{1},l_{2}],  \\[2mm]
\psi(x,0)\leq \psi(x,T),\; &\, x\in[l_{1},l_{2}],
\end{array} \right.
\end{eqnarray}
then $\lambda_{1}\leq \overline{\lambda}$ holds. Additionally, $\lambda_{1}= \overline{\lambda}$ if and only if the equalities in \eqref{3-7} all hold.
\begin{proof}
We first prove the first assertion. Multiplying the first equation in \eqref{3-7} by $\zeta_{1}$ and the first equation in \eqref{3-6} by $\phi$ yield that
\begin{eqnarray}\label{3-8}
\left\{
\begin{array}{ll}
\phi_{t}\zeta_{1}\leq-\delta_{1} \phi\zeta_{1}+\overline{\lambda} \phi\zeta_{1} , \\[2mm]
\zeta_{1t}\phi=\delta_{1} \zeta_{1}\phi-\mu_{1} \zeta_{1}\phi.
\end{array}
\right.
\end{eqnarray}
In the following, integrating both sides of the equations in \eqref{3-8} over $(l_{1},l_{2})\times (0, \tau]$, and summing the two equations give that
\begin{equation}\label{3-9}
\int_{l_{1}}^{l_{2}}\phi\zeta_{1}|_{0}^{\tau}dx\leq (\overline{\lambda}-\mu_{1})\int_{l_{1}}^{l_{2}}\int_{0}^{\tau}\zeta_{1}\phi dtdx.
\end{equation}
A similar method to that used to obtain \eqref{3-9} yields that
\begin{equation}\label{3-10}
\int_{l_{1}}^{l_{2}}\phi\zeta_{1}|_{\tau}^{T}dx\leq \int_{l_{1}}^{l_{2}}\int_{\tau}^{T}(a_{12}\zeta_{1}\psi-f'(0)\eta_{1}\phi) dtdx+(\overline{\lambda}-\mu_{1})\int_{l_{1}}^{l_{2}}\int_{\tau}^{T}\zeta_{1}\phi dtdx.
\end{equation}
Adding \eqref{3-9} and \eqref{3-10} gives that
\begin{equation}\label{3-11}
\int_{l_{1}}^{l_{2}}\phi\zeta_{1}|_{0}^{T}dx\leq \int_{l_{1}}^{l_{2}}\int_{\tau}^{T}(a_{12}\zeta_{1}\psi-f'(0)\eta_{1}\phi) dtdx+(\overline{\lambda}-\mu_{1})\int_{l_{1}}^{l_{2}}\int_{0}^{T}\zeta_{1}\phi dtdx.
\end{equation}
By employing methods analogous to those used in deriving \eqref{3-11}, we obtain that
\begin{equation}\label{3-12}
\int_{l_{1}}^{l_{2}}\psi\eta_{1}|_{0}^{T}dx\leq \int_{l_{1}}^{l_{2}}\int_{\tau}^{T}(f'(0)\eta_{1}\phi-a_{12}\zeta_{1}\psi) dtdx+(\overline{\lambda}-\mu_{1})\int_{l_{1}}^{l_{2}}\int_{0}^{T}\eta_{1}\psi dtdx.
\end{equation}
Adding \eqref{3-11} and \eqref{3-12} yields that
\begin{equation}\label{3-13}
 (\overline{\lambda}-\lambda_{1})\int_{l_{1}}^{l_{2}}\int_{0}^{T}(\eta_{1}\psi +\zeta_{1}\phi )dtdx\geq 0.
\end{equation}
Since $(\zeta_{1}, \eta_{1})\in \mathds{E}^{++}$ and $(\phi, \psi)\in \mathds{E}^{+}-\{\textbf{0}\}$, it follows from \eqref{3-13} that the first assertion holds.

Next, we prove the second assertion. If the equalities in \eqref{3-7} all hold, the assertion $\overline{\lambda}=\lambda_{1}$ can be obtained by repeating the proof of the first assertion. Finally, it remains to show that if  $\overline{\lambda}=\lambda_{1}$, then the equalities in \eqref{3-7} all hold. Arguing indirectly, we assume that there exists an inequality in \eqref{3-7}. Then, repeating the proof of the first assertion yields that $\lambda_{1}< \overline{\lambda}$, which is a contradiction to $\overline{\lambda}=\lambda_{1}$. The proof is now completed.
\end{proof}
\end{lemma}
\begin{lemma}\label{lemma 3.1.2}
If there exists a vector-valued eigenfunction $(\phi, \psi)\in \mathds{E}^{+}-\{\textbf{0}\}$ and a number $\underline{\lambda}$ satisfying
\begin{eqnarray}\label{3-14}
\left\{
\begin{array}{ll}
\phi_{t}\geq-\delta_{1} \phi+\underline{\lambda} \phi ,\; &\,x\in(l_{1},l_{2}), t\in(0, \tau],  \\[2mm]
\psi_{t}\geq d_{2}\Delta \psi-\delta_{2} \psi+\underline{\lambda} \psi ,\; &\,x\in(l_{1},l_{2}), t\in(0, \tau],  \\[2mm]
\phi_{t}\geq d_{1}\Delta \phi-a_{11}\phi+a_{12}\psi+\underline{\lambda} \phi ,\; &\,x\in(l_{1},l_{2}), t\in(\tau, T],  \\[2mm]
\psi_{t}\geq d_{2}\Delta \psi-a_{22}\psi+f'(0) \phi+\underline{\lambda} \psi  ,\; &\,x\in(l_{1},l_{2}), t\in(\tau, T],  \\[2mm]
\phi(x,t)=\psi(x,t)=0,\; &\, x\in\{l_{1},l_{2}\},  t\in[0, T],\\[2mm]
\phi(x,0)\geq H'(0)\phi(x,T),\; &\, x\in[l_{1},l_{2}],  \\[2mm]
\psi(x,0)\geq \psi(x,T),\; &\, x\in[l_{1},l_{2}],
\end{array} \right.
\end{eqnarray}
then $\lambda_{1}\geq \underline{\lambda}$ holds. Additionally, $\lambda_{1}= \underline{\lambda}$ if and only if the equalities hold in \eqref{3-14}.
\begin{proof}
The proof for this lemma is similar to \autoref{lemma 3.1.1}, and we omit it here.
\end{proof}
\end{lemma}
In order to understand the effects of pulse intensity, region size, and dry season duration on the dynamical behavior of problem \eqref{3-1}, it suffices to analyze their influence on the principal eigenvalue $\lambda_{1}$.
\begin{lemma}\label{lemma 3.1.3}
Suppose that $[l_{1},l_{2}]=[-l,l]$. Let $H'(0)=\theta$. In order to emphasize the dependence of the principal eigenvalue $\lambda_{1}$ on $l$, $\theta$, and $\tau$, we denote $\lambda_{1}$ by $\lambda_{1}(l)$, $\lambda_{1}(\theta)$, and $\lambda_{1}(\tau)$, respectively. Then, we have the following assertions:
\begin{enumerate}
\item[$(1)$]
$\lambda_{1}(\theta)$ is strongly decreasing and continuous with respect to $\theta$;
\item[$(2)$]
$\lambda_{1}(l)$ is strongly decreasing and continuous with respect to $l$;
\item[$(3)$]
$\lambda_{1}(\tau)$ is strongly increasing and continuous with respect to $\tau$ if $\delta_{1}\geq a_{11}+\frac{d_{1}\pi^{2}}{4l^{2}}$ and $\delta_{2}\geq a_{22}$.
\end{enumerate}
\begin{proof}
(1) By employing \textcolor{blue}{Lemmas} \ref{lemma 3.1.1} and \ref{lemma 3.1.2}, the assertion can be established following an approach similar to that in \cite[Lemma 3.4(1)]{Zhou-Pedersen-Lin-ANA}. Since the proof is analogous, we omit the details here.

(2) We first prove the monotonicity of $\lambda_{1}(l)$ with respect to $l$. For $0<l_{1}<l_{2}$, let $(\phi_{l_{2}}, \psi_{l_{2}})$ be a strongly positive eigenfunction pair associated with $\lambda_{1}(l_{2})$. By using the method of separation of variables, we have that
\begin{equation*}
\phi_{l_{2}}(x,t)=\cos\big(\frac{\pi}{2l_{2}}x\big)\Phi_{l_{2}}(t) \text{~and~}  \psi_{l_{2}}(x,t)=\cos\big(\frac{\pi}{2l_{2}}x\big)\Psi_{l_{2}}(t).
\end{equation*}
Let
\begin{equation*}
\phi(x,t)=\cos\big(\frac{\pi}{2l_{1}}x\big)\Phi_{l_{2}}(t) \text{~and~}  \psi(x,t)=\cos\big(\frac{\pi}{2l_{1}}x\big)\Psi_{l_{2}}(t).
\end{equation*}
Next, we verify that the triple $(\phi, \psi, \lambda_{1}(l_{2}))$ satisfies \eqref{3-14}. For $x\in(-l_{1}, l_{1})$ and $t\in(\tau, T]$, we have that
\begin{equation*}
\begin{split}
&\phi_{t}- d_{1}\Delta \phi+a_{11}\phi-a_{12}\psi-\lambda_{1}(l_{2}) \phi \\
= &\cos\big(\frac{\pi}{2l_{1}}x\big)\big[\Phi'_{l_{2}}(t)-[\lambda_{1}(l_{2})-d_{1}\kappa_{1}(l_{1}) -a_{11}]\Phi_{l_{2}}(t)-a_{12}\Psi_{l_{2}}(t)\big]\\
> &\cos\big(\frac{\pi}{2l_{1}}x\big)\big[\Phi'_{l_{2}}(t)-[\lambda_{1}(l_{2})-d_{1}\kappa_{1}(l_{2}) -a_{11}]\Phi_{l_{2}}(t)-a_{12}\Psi_{l_{2}}(t)\big]\\
=&0.
\end{split}
\end{equation*}
Similarity, the first, second, and fourth inequalities in \eqref{3-14} can also be obtained. The remainder of \eqref{3-14} can be easily verified. Therefore, it follows from \autoref{lemma 3.1.2} that  $\lambda_{1}(l_{1})> \lambda_{1}(l_{2})$.

Next we prove the continuity of $\lambda_{1}(l)$ with respect to $l>0$. To obtain this, we will prove that for any $l_{0}>0$ and small $\epsilon>0$, there exists some constant $\delta>0$ such that
\begin{equation}\label{3-22}
|\lambda_{1}(l)-\lambda_{1}(l_{0})|<\epsilon\text{~if~}|l-l_{0}|<\delta.
\end{equation}
We first discuss the case of $l\in(l_{0}-\delta_{1},l_{0} ]$. By the above-proved monotonicity, one can obtain that
\begin{equation*}
\lambda_{1}(l_{0})-\epsilon< \lambda_{1}(l_{0})\leq \lambda_{1}(l).
\end{equation*}
Therefore, it suffices to find $\delta_{1}>0$ such that $\lambda_{1}(l)<\lambda_{1}(l_{0})+\epsilon$ for all $l\in(l_{0}-\delta_{1},l_{0} ]$. It is easily seen that there exists some constant $\delta_{1}>0$ such that
\begin{equation*}
d_{1}\kappa_{1}(l)<d_{1}\kappa_{1}(l_{0})+\epsilon\text{~and~}d_{2}\kappa_{1}(l)<d_{2}\kappa_{1}(l_{0})+\epsilon
\end{equation*}
for all $l\in(l_{0}-\delta_{1},l_{0} ]$. We denote by $(\phi_{l_{0}}, \psi_{l_{0}})$ a strongly positive eigenfunction pair associated with $\lambda_{1}(l_{0})$. It follows from the method of separation of variables that
\begin{equation*}
\phi_{l_{0}}(x,t)=\cos\big(\frac{\pi}{2l_{0}}x\big)\Phi_{l_{0}}(t) \text{~and~}  \psi_{l_{0}}(x,t)=\cos\big(\frac{\pi}{2l_{0}}x\big)\Psi_{l_{0}}(t).
\end{equation*}
Let
\begin{equation}\label{3-23}
\phi(x,t)=\cos\big(\frac{\pi}{2l}x\big)\Phi_{l_{0}}(t) \text{~and~}  \psi(x,t)=\cos\big(\frac{\pi}{2l}x\big)\Psi_{l_{0}}(t).
\end{equation}
A regular calculation yields that the triple $(\phi, \psi, \lambda_{1}(l_{0})+\epsilon)$ satisfies \eqref{3-7}. Therefore, it follows from \autoref{lemma 3.1.1} that $\lambda_{1}(l)<\lambda_{1}(l_{0})+\epsilon$ for all $l\in(l_{0}-\delta_{1},l_{0} ]$.

Finally, we consider the case of $l\in[l_{0},l_{0}+\delta_{2} )$. It suffices to prove that $\lambda_{1}(l_{0})-\epsilon<\lambda_{1}(l)$ for the reason that $\lambda_{1}(l)\leq \lambda_{1}(l_{0})<\lambda_{1}(l_{0})+\epsilon$ can be obtained by using the monotonicity of $\lambda_{1}(l)$. The strong monotonicity and continuity of $\kappa_{1}(l)$ tell us that there exists a positive constant $\delta_{2}>0$ such that
\begin{equation*}
d_{1}\kappa_{1}(l_{0})<d_{1}\kappa_{1}(l)+\epsilon\text{~and~}d_{2}\kappa_{1}(l_{0})<d_{2}\kappa_{1}(l)+\epsilon
\end{equation*}
for all $l\in[l_{0},l_{0}+\delta_{2} )$. Through verification, it can be concluded that the triple $(\phi, \psi, \lambda_{1}(l_{0})-\epsilon)$ satisfies \eqref{3-14}, where $\phi$ and $\psi$ are defined in \eqref{3-23}. Then, \autoref{lemma 3.1.2} yields that
$\lambda_{1}(l_{0})-\epsilon<\lambda_{1}(l)$ for all $l\in[l_{0},l_{0}+\delta_{2} )$. By taking $\delta=\min\{\delta_{1},\delta_{2} \}$,  \eqref{3-22} holds true. We have now proved the continuity.

(3) We begin by establishing the monotonicity property of $\lambda_{1}(\tau)$. For any $0<\tau_{1}<\tau_{2}<T$, let $(\lambda_{1}(\tau_{1}),\phi, \psi)$ and $(\mu_{1}(\tau_{2}),\zeta, \eta)$ are the principal eigenpairs of the eigenvalue problem \eqref{3-3} with $\tau=\tau_{1}$ and the adjoint problem \eqref{3-6} with $\tau=\tau_{2}$, respectively. Multiplying the first equations in \eqref{3-3} and \eqref{3-6} by $\zeta$ and $\phi$, respectively, we have that
\begin{eqnarray}\label{3-24}
\left\{
\begin{array}{ll}
\phi_{t}\zeta=-\delta_{1} \phi\zeta+\lambda_{1}(\tau_{1}) \phi\zeta, \\[2mm]
\zeta_{t}\phi=\delta_{1} \zeta\phi-\mu_{1}(\tau_{2}) \zeta\phi.
\end{array}
\right.
\end{eqnarray}
Then, integrating both sides of the equations in \eqref{3-24} over $(-l,l)\times(0,\tau_{1}]$ and adding these together yield that
\begin{equation}\label{3-25}
\int_{-l}^{l}\phi(t)\zeta(t)|_{0}^{\tau_{1}}dx= [\lambda_{1}(\tau_{1})-\lambda_{1}(\tau_{2})]\int_{-l}^{l}\int_{0}^{\tau_{1}}\zeta\phi dtdx.
\end{equation}
By performing the same procedure to the third equations of \eqref{3-3} and \eqref{3-6}, except that the integration interval $(-l,l)\times(0,\tau_{1}]$ is replaced by $(-l,l)\times(\tau_{1}, T]$, we obtain that
\begin{equation}\label{3-26}
\begin{split}
\int_{-l}^{l}\phi(t)\zeta(t)|_{\tau_{1}}^{T}dx=& [\lambda_{1}(\tau_{1})-\lambda_{1}(\tau_{2})]\int_{-l}^{l}\int_{\tau_{1}}^{T}\zeta\phi dtdx+\int_{-l}^{l}\int_{\tau_{2}}^{T}(a_{12}\zeta\psi- f'(0)\eta\phi )dtdx\\
&+\int_{-l}^{l}\int_{\tau_{1}}^{\tau_{2}}(d_{1}\Delta\phi\zeta-a_{11}\phi\zeta+a_{12}\zeta\psi+\delta_{1}\phi\zeta)dtdx.
\end{split}
\end{equation}
Adding \eqref{3-25} and \eqref{3-26} gives that
\begin{equation}\label{3-27}
\begin{split}
(\lambda_{1}(\tau_{2})-\lambda_{1}(\tau_{1}))\int_{-l}^{l}\int_{0}^{T}\zeta\phi dtdx=&\int_{-l}^{l}\int_{\tau_{1}}^{\tau_{2}}(d_{1}\Delta\phi\zeta-a_{11}\phi\zeta+a_{12}\zeta\psi+\delta_{1}\phi\zeta)dtdx\\
&+\int_{-l}^{l}\int_{\tau_{2}}^{T}(a_{12}\zeta\psi- f'(0)\eta\phi )dtdx.
\end{split}
\end{equation}
Similarly, we have that
\begin{equation}\label{3-28}
\begin{split}
(\lambda_{1}(\tau_{2})-\lambda_{1}(\tau_{1}))\int_{-l}^{l}\int_{0}^{T}\eta\psi dtdx=&\int_{-l}^{l}\int_{\tau_{1}}^{\tau_{2}}(f'(0)\eta\phi-a_{22}\psi\eta+\delta_{2}\psi\eta)dtdx\\
&+\int_{-l}^{l}\int_{\tau_{2}}^{T}( f'(0)\eta\phi -a_{12}\zeta\psi)dtdx.
\end{split}
\end{equation}
By adding \eqref{3-27} and \eqref{3-28}, it follows that
\begin{equation}\label{3-29}
\begin{split}
(\lambda_{1}(\tau_{2})-\lambda_{1}(\tau_{1}))\int_{-l}^{l}\int_{0}^{T}(\zeta\phi+\eta\psi )dtdx>&\int_{-l}^{l}\int_{\tau_{1}}^{\tau_{2}}(d_{1}\Delta\phi\zeta-a_{11}\phi\zeta+\delta_{1}\phi\zeta)dtdx\\
&+\int_{-l}^{l}\int_{\tau_{1}}^{\tau_{2}}(-a_{22}\psi\eta+\delta_{2}\psi\eta)dtdx.
\end{split}
\end{equation}
According to the method of separation of variables, we let
\begin{equation*}
\phi(x,t)=\chi_{1}(x)\Phi(t),~ \psi(x,t)=\chi_{1}(x)\Psi(t), ~\zeta(x,t)=\chi_{1}(x)\widetilde{\Phi}(t), \text{~and~}\eta(x,t)=\chi_{1}(x)\widetilde{\Psi}(t),
\end{equation*}
where $\chi_{1}(x)$ is the eigenfunction corresponding to the principal eigenvalue $\kappa_{1}$ of problem \eqref{3-30} with $l_{1}=-l$ and $l_{2}=l$, $(\Phi(t), \Psi(t))$ is the strongly positive eigenfunction pair corresponding to the principal eigenvalue $\lambda_{1}(\tau_{1})$ of problem \eqref{3-15} with $\kappa_{1}=\frac{\pi^{2}}{4l^{2}}$, and $(\widetilde{\Phi},\widetilde{\Psi})$ is the strongly positive eigenfunction pair corresponding to the principal eigenvalue $\lambda_{1}(\tau_{2})$ of the following problem
\begin{eqnarray*}
\left\{
\begin{array}{ll}
-\widetilde{\Phi}'(t)=[\lambda-\delta_{1} ]\widetilde{\Phi}(t),\; &\,  t\in(0, \tau_{2}], \\[2mm]
-\widetilde{\Psi}'(t)=[\lambda-\delta_{2}-\frac{d_{2}\pi^{2}}{4l^{2}} ]\widetilde{\Psi}(t),\; &\,  t\in(0, \tau_{2}], \\[2mm]
-\widetilde{\Phi}'(t)=[\lambda-\frac{d_{1}\pi^{2}}{4l^{2}} -a_{11}]\widetilde{\Phi}(t)+f'(0)\widetilde{\Psi}(t),\; &\,  t\in(\tau_{2}, T], \\[2mm]
-\widetilde{\Psi}'(t)=[\lambda-\frac{d_{2}\pi^{2}}{4l^{2}} -a_{22}]\widetilde{\Psi}(t)+a_{12}\widetilde{\Phi}(t),\; &\,  t\in(\tau_{2}, T], \\[2mm]
\widetilde{\Phi}(0)=\frac{1}{H'(0)}\widetilde{\Phi}(T), ~~~~~~\widetilde{\Psi}(0)=\widetilde{\Psi}(T).
\end{array} \right.
\end{eqnarray*}
Then, using the assumption that $\delta_{1}\geq a_{11}+\frac{d_{1}\pi^{2}}{4l^{2}}$ and $\delta_{2}\geq a_{22}$, \eqref{3-29} can be further written as
\begin{equation*}
[\lambda_{1}(\tau_{2})-\lambda_{1}(\tau_{1})]\int_{0}^{T}\big[\widetilde{\Phi}(t)\Phi(t)+\widetilde{\Psi}(t)\Psi(t) \big]dt>0,
\end{equation*}
which implies that $\lambda_{1}(\tau_{2})>\lambda_{1}(\tau_{1})$. Additionally,  the continuity property of $\lambda_{1}(\tau)$ can be derived from \eqref{3-20} due to the continuous dependence of the solution with respect to the coefficients. This ends the proof.
\end{proof}
\end{lemma}
\autoref{lemma 3.1.3}\textcolor{blue}{(2)} shows that the principal eigenvalue $\lambda_{1}(l)$ of problem \eqref{3-3} is strongly decreasing with respect to $l$. To conclude this subsection, we shall investigate the limiting behavior of this principal eigenvalue as $l$ tends to infinity. For this purpose, we consider the following problem
\begin{eqnarray}\label{3-31}
\left\{
\begin{array}{ll}
\phi'(t)=[\nu-\delta_{1} ]\phi(t),\; &\,  t\in(0, \tau], \\[2mm]
\psi'(t)=[\nu-\delta_{2}]\psi(t),\; &\,  t\in(0, \tau], \\[2mm]
\phi'(t)=[\nu -a_{11}]\phi(t)+a_{12}\psi(t),\; &\,  t\in(\tau, T], \\[2mm]
\psi'(t)=[\nu -a_{22}]\psi(t)+f'(0)\phi(t),\; &\,  t\in(\tau, T], \\[2mm]
\phi(0)=H'(0)\phi(T), ~~~~~~\psi(0)=\psi(T).
\end{array} \right.
\end{eqnarray}
By using the same method as in \autoref{lemma 3.1.4}, it follows that problem \eqref{3-31} has a unique and algebraically simple eigenvalue $\nu_{1}$ with a strongly positive vector-valued eigenfunction $(\phi_{1}(t), \psi_{1}(t))$.
\begin{lemma}\label{lemma 3.1.7}
Suppose $l_{2}=-l_{1}=l$. Let $\lambda_{1}(l)$ and $\nu_{1}$ be the principal eigenvalues of \eqref{3-3} and \eqref{3-31}, respectively. Then, it follows that
\begin{equation*}
\lim\limits_{l\rightarrow +\infty} \lambda_{1}(l)=\nu_{1}.
\end{equation*}
\begin{proof}
By employing \textcolor{blue}{Lemmas} \ref{lemma 3.1.1}, \ref{lemma 3.1.2}, and \ref{lemma 3.1.3}\textcolor{blue}{(2)}, the desired conclusion can be established following a methodology similar to that in \cite[Lemma 3.5]{Zhou-Pedersen-Lin-ANA} or \cite[Lemma 3.9]{Wang-Du-JDE}. For brevity, we omit the proof here.
\end{proof}
\end{lemma}
\subsection{The long-time dynamics of the fixed boundary problem}\label{Section-3-2}
It follows from \autoref{lemma 3.0.1} that the fixed-boundary problem \eqref{3-1} admits a unique nonnegative solution for all $t>0$. In this subsection,
we always use $(u(x,t), v(x,t))$ to denote this unique solution.

For problem \eqref{3-1}, its steady state has the following form
\begin{eqnarray}\label{3-32}
\left\{
\begin{array}{ll}
w_{t}=-\delta_{1} w ,\; &\,x\in(l_{1},l_{2}), t\in(0^{+}, \tau],  \\[2mm]
z_{t}=d_{2}\Delta z-\delta_{2} z ,\; &\,x\in(l_{1},l_{2}), t\in(0^{+}, \tau],  \\[2mm]
w_{t}=d_{1}\Delta w-a_{11}w+a_{12}z ,\; &\,x\in(l_{1},l_{2}), t\in(\tau, T],  \\[2mm]
z_{t}=d_{2}\Delta z-a_{22}z+f(w) ,\; &\,x\in(l_{1},l_{2}), t\in(\tau, T],  \\[2mm]
w(x,t)=z(x,t)=0,\; &\, x\in\{l_{1},l_{2}\},  t\in(0, T],\\[2mm]
w(x,0^{+})=H(w(x,0)),\; &\, x\in(l_{1},l_{2}),  \\[2mm]
z(x,0^{+})=z(x,0),\; &\, x\in(l_{1},l_{2}),  \\[2mm]
w(x,0)=w(x,T),z(x,0)=z(x,T),\; &\,x\in[l_{1},l_{2}].
\end{array} \right.
\end{eqnarray}

Here, the method of upper and lower solutions will be employed to study the steady state problem \eqref{3-32}. To begin with, we first introduce the
definition of the upper and lower solutions to problem \eqref{3-32}.
\begin{definition}
The nonnegative vector-valued functions
\begin{equation*}
(\underline{w}, \underline{z}), (\overline{w}, \overline{z})\in \mathcal{PC}^{2,1}(l_{1},l_{2}, 0)
\end{equation*}
satisfying
\begin{equation*}
\underline{w}\leq \overline{w} \text{~and~}\underline{z}\leq \overline{z}
\end{equation*}
are the ordered lower and upper solutions of problem \eqref{3-32}, respectively, if $(\underline{w}, \underline{z})$ satisfies
\begin{eqnarray}\label{pc-1}
\left\{
\begin{array}{ll}
\underline{w}_{t}\leq -\delta_{1} \underline{w} ,\; &\,x\in(l_{1},l_{2}), t\in(0^{+}, \tau],  \\[2mm]
\underline{z}_{t}\leq d_{2}\Delta \underline{z}-\delta_{2} \underline{z} ,\; &\,x\in(l_{1},l_{2}), t\in(0^{+}, \tau],  \\[2mm]
\underline{w}_{t}\leq d_{1}\Delta \underline{w}-a_{11}\underline{w}+a_{12}\underline{z} ,\; &\,x\in(l_{1},l_{2}), t\in(\tau, T],  \\[2mm]
\underline{z}_{t}\leq d_{2}\Delta \underline{z}-a_{22}\underline{z}+f(\underline{w}) ,\; &\,x\in(l_{1},l_{2}), t\in(\tau, T],  \\[2mm]
\underline{w}(x,t)\leq0, ~\underline{z}(x,t)\leq0,\; &\, x\in\{l_{1},l_{2}\},  t\in(0, T],\\[2mm]
\underline{w}(x,0^{+})\leq H(\underline{w}(x,0)),\; &\, x\in(l_{1},l_{2}),  \\[2mm]
\underline{z}(x,0^{+})\leq \underline{z}(x,0),\; &\, x\in(l_{1},l_{2}),  \\[2mm]
\underline{w}(x,0)\leq \underline{w}(x,T),\underline{z}(x,0)\leq \underline{z}(x,T),\; &\,x\in[l_{1},l_{2}],
\end{array} \right.
\end{eqnarray}
and $(\overline{w}, \overline{z})$ satisfies \eqref{pc-1} with the inequality sign reversed ($\leq$ replaced by $\geq$), where $\mathcal{PC}^{2,1}(l_{1},l_{2}, 0)$ is defined in \eqref{pc}.
\end{definition}
Now, we present the main result of this subsection.
\begin{theorem}\label{theorem 3.2}
Let $\lambda_{1}$ be the principal eigenvalue of problem \eqref{3-3}. Denote by $(u,v)$ and $(w,z)$ the solutions to \eqref{3-1} and \eqref{3-32}, respectively. Then, the following assertions are valid:
\begin{enumerate}
\item[$(1)$]
If $\lambda_{1}< 0$, then the steady state problem \eqref{3-32} admits a unique positive periodic solution $(w,z)\in\mathcal{PC}^{2,1}(l_{1},l_{2}, 0)$, where $\mathcal{PC}^{2,1}(l_{1},l_{2}, 0)$ is defined in \eqref{pc}. Additionally, $(u,v)$ satisfies
\begin{equation}\label{3-38}
\lim\limits_{n\rightarrow+\infty}\big(u(x,t+nT),v(x,t+nT)\big)=\big(w(x,t), z(x,t)\big)
\end{equation}
uniformly for $(x,t)\in[l_{1},l_{2}]\times[0, +\infty)$;
\item[$(2)$]
If $\lambda_{1}\geq 0$, then $(0,0)$ is the only nonnegative solution of the steady state problem \eqref{3-32}. Additionally, $(u,v)$ satisfies
\begin{equation}\label{3-39}
\lim\limits_{t\rightarrow+\infty}\big(u(x,t),v(x,t)\big)=\big(0, 0\big) \text{~uniformly~for~}x\in[l_{1},l_{2}].
\end{equation}
\end{enumerate}
\begin{proof}
(1) We begin by proving the first assertion. Its proof is lengthy and has been divided into three parts for clarity.

\textbf{(a) The ordered upper and lower solutions}

We construct the following functional pair
\begin{equation*}
(\overline{w}(x,t), \overline{z}(x,t))=(C_{2}, C_{3}) \text{~for~all~}(x,t)\in[l_{1},l_{2}]\times[0,T],
\end{equation*}
where $C_{2}$ and $C_{3}$ are defined in \eqref{2-2}. For any given $\lambda_{1}< 0$, there exists a positive constant $\upsilon$ such that $\lambda_{1}+\upsilon< 0$. Then, we also construct the following functions
\begin{eqnarray*}
\underline{w}(x,t)=
\left\{
\begin{array}{ll}
\epsilon \phi_{1}(x,t),\; &\, (x,t)\in [l_{1},l_{2}]\times\{0\}, \\[2mm]
\epsilon e^{(\lambda_{1}+\upsilon)T} \phi_{1}(x,t),\; &\, (x,t)\in [l_{1},l_{2}]\times\{0^{+}\}, \\[2mm]
\epsilon e^{(\lambda_{1}+\upsilon)(T-t)}\phi_{1}(x,t),\; &\, (x,t)\in [l_{1},l_{2}]\times(0^{+}, T],
\end{array} \right.
\end{eqnarray*}
and
\begin{eqnarray*}
\underline{z}(x,t)=
\left\{
\begin{array}{ll}
\epsilon \psi_{1}(x,t),\; &\, (x,t)\in [l_{1},l_{2}]\times\{0\}, \\[2mm]
\epsilon e^{(\lambda_{1}+\upsilon)T} \psi_{1}(x,t),\; &\, (x,t)\in [l_{1},l_{2}]\times\{0^{+}\}, \\[2mm]
\epsilon e^{(\lambda_{1}+\upsilon)(T-t)}\psi_{1}(x,t),\; &\, (x,t)\in [l_{1},l_{2}]\times(0^{+}, T],
\end{array} \right.
\end{eqnarray*}
where $(\phi_{1}, \psi_{1})$ is the eigenfunction pair corresponding to the principal eigenvalue $\lambda_{1}$ of problem \eqref{3-2}, and $\epsilon$ is a positive constant such that
\begin{equation*}
\epsilon e^{(\lambda_{1}+\upsilon)T}H'(0)\phi_{1}(x,0)\leq H(\epsilon\phi_{1}(x,0))
\end{equation*}
for all $x\in(l_{1}, l_{2})$ and
\begin{equation*}
\epsilon e^{(\lambda_{1}+\upsilon)(T-t)}[f'(0)\phi_{1}(x,t)-\upsilon\psi_{1}(x,t)] \leq f(\epsilon e^{(\lambda_{1}+\upsilon)(T-t)}\phi_{1}(x,t))
\end{equation*}
for all $(x,t)\in(l_{1},l_{2})\times(\tau,T]$. By verification, the constructed pairs $(\overline{w}(x,t), \overline{z}(x,t))$ and $(\underline{w}(x,t), \underline{z}(x,t))$ are the ordered lower and upper solutions of the steady state problem \eqref{3-32}, respectively.

\textbf{(b)  The existence of the positive periodic solution}

For the convenience of writing, we let
\begin{equation*}
\begin{split}
&H_{1}(w,z)=(\delta_{2}+a_{11}+a_{22})w, ~~~~~~~~~~H_{2}(w,z)=(\delta_{1}+a_{11}+a_{22})z, ~\\
&H_{3}(w,z)=(\delta_{1}+\delta_{2}+a_{22})w+a_{12}z, ~H_{4}(w,z)=(\delta_{1}+\delta_{2}+a_{11})z+f(w).
\end{split}
\end{equation*}
It is easily seen that $H_{i}$ is monotonically nondecreasing with respect to $w$ and $z$ for $i=1,2,3,4$. Then, we construct the following iteration rule
\begin{eqnarray}\label{3-33}
\left\{
\begin{array}{ll}
w^{i}_{t}+\gamma w^{i}=H_{1}(w^{i-1},z^{i-1}) ,\; &\,x\in(l_{1},l_{2}), t\in((mT)^{+}, mT+\tau],  \\[2mm]
z^{i}_{t}-d_{2}\Delta z^{i}+\gamma z^{i}=H_{2}(w^{i-1},z^{i-1}) ,\; &\,x\in(l_{1},l_{2}), t\in((mT)^{+}, mT+\tau],  \\[2mm]
w^{i}_{t}-d_{1}\Delta w^{i}+\gamma w^{i}=H_{3}(w^{i-1},z^{i-1}) ,\; &\,x\in(l_{1},l_{2}), t\in(mT+\tau, (m+1)T],  \\[2mm]
z^{i}_{t}-d_{2}\Delta z^{i}+\gamma z^{i}=H_{4}(w^{i-1},z^{i-1}) ,\; &\,x\in(l_{1},l_{2}), t\in(mT+\tau, (m+1)T],  \\[2mm]
w^{i}(x,t)=z^{i}(x,t)=0,\; &\, x\in\{l_{1},l_{2}\},  t\in(0, (m+1)T],  \\[2mm]
w^{i}(x,(mT)^{+})=H(w^{i-1}(x,(m+1)T)),\; &\, x\in(l_{1},l_{2}),  \\[2mm]
z^{i}(x,(mT)^{+})=z^{i-1}(x,(m+1)T),\; &\, x\in(l_{1},l_{2}),  \\[2mm]
w^{i}(x,mT)=w^{i-1}(x,(m+1)T),\; &\,x\in[l_{1},l_{2}],\\[2mm]
z^{i}(x,mT)=z^{i-1}(x,(m+1)T),\; &\,x\in[l_{1},l_{2}],
\end{array} \right.
\end{eqnarray}
where $m=0$ and
\begin{equation}\label{XJD}
\gamma=\delta_{1}+\delta_{2}+a_{11}+a_{22}.
\end{equation}
Choosing the upper and lower solutions constructed in \textbf{(a)} as initial iteration, we can construct the iteration sequences $\big\{(\overline{w}^{i}, \overline{z}^{i})\big\}_{i=0}^{n}$ and $\big\{(\underline{w}^{i}, \underline{z}^{i})\big\}_{i=0}^{n}$ by the above iteration rule \eqref{3-33}. The method of induction and the positivity lemma for parabolic problem yield the monotone property
\begin{equation*}
(\underline{w}^{i}, \underline{z}^{i})\leq (\underline{w}^{i+1}, \underline{z}^{i+1})\leq (\overline{w}^{i+1}, \overline{z}^{i+1})\leq (\overline{w}^{i}, \overline{z}^{i})
\end{equation*}
for $x\in[l_{1},l_{2}]$ and $t\in[0,T]$. Next, it follows from the monotone bounded convergence theorem that
\begin{equation}\label{3-5}
\lim\limits_{n\rightarrow+\infty}(\overline{w}^{n}, \overline{z}^{n})= (w_{1},z_{1})\text{~~and~~}\lim\limits_{n\rightarrow+\infty}(\underline{w}^{n}, \underline{z}^{n})= (w_{2},z_{2}),
\end{equation}
which implies, after applying standard arguments, that $ (w_{1},z_{1})$ and $ (w_{2},z_{2})$ are the solutions of problem \eqref{3-32}. Additionally, employing a proof technique analogous to \cite[Theorem 4.2]{Zhou-Lin-Santos-CFS}, the pairs $ (w_{1},z_{1})$ and $ (w_{2},z_{2})$ are the maximal and minimal positive periodic solutions of problem \eqref{3-32}, respectively.

\textbf{(c)  The uniqueness of the positive periodic solution}

Finally, we show that problem \eqref{3-32} admits a unique positive periodic solution. Arguing indirectly, we assume that $ (w_{1},z_{1})$ and $ (w_{2},z_{2})$ are two
distinct solutions of problem \eqref{3-32}. Since
$(w_{1},z_{1})$ and $(w_{2},z_{2})$ are both strongly positive,
\begin{equation*}
\zeta_{0}=:\inf\big\{\zeta\geq 1 | \zeta(w_{1},z_{1})\succeq(w_{2},z_{2})\text{~for~}(x,t)\in[l_{1},l_{2}]\times[0,T] \big\}
\end{equation*}
is well-defined and finite, which implies that $\zeta_{0}(w_{1},z_{1})\succeq(w_{2},z_{2})$ for $(x,t)\in[l_{1},l_{2}]\times[0,T]$. We claim $\zeta_{0}=1$. If not, $\zeta_{0}>1$,
there would exist some $(x_{0},t_{0})\in[l_{1},l_{2}]\times(0,T]$ such that
\begin{equation*}
\mbox{(i)}~\zeta_{0}w_{1}(x_{0},t_{0}) =w_{2}(x_{0},t_{0}) \text{~or~} \mbox{(ii)}~\zeta_{0}z_{1}(x_{0},t_{0}) =z_{2}(x_{0},t_{0}).
\end{equation*}
For the convenience of writing, let
\begin{equation*}
U(x,t)=\zeta_{0}w_{1}- w_{2}\text{~and~}V(x,t)=\zeta_{0}z_{1}- z_{2} .
\end{equation*}
If $\mbox{(i)}$ holds, then the direct calculation yields that
\begin{eqnarray}\label{3-34}
\left\{
\begin{array}{ll}
U_{t}+\gamma U=(\delta_{2}+a_{11}+a_{22} )U\geq 0 ,\; &\,x\in(l_{1},l_{2}), t\in(0^{+}, \tau],  \\[2mm]
U_{t}-d_{1}\Delta U+\gamma U=(\gamma-a_{11} )U+a_{12}V\geq 0 ,\; &\,x\in(l_{1},l_{2}), t\in(\tau, T],  \\[2mm]
U(x,t)=0,\; &\, x\in\{l_{1},l_{2}\},  t\in(0^{+}, T],\\[2mm]
U(x,0^{+})\geq, \not \equiv 0,\; &\, x\in[l_{1},l_{2}],  \\[2mm]
\end{array} \right.
\end{eqnarray}
where $\gamma$ is defined in \eqref{XJD}. In virtue of the strong maximum principle, one of the following two cases must hold:

\textbf{Case 1}: $U(x,t)>0$ holds for $(x,t)\in(l_{1},l_{2})\times(\tau,T]$. Additionally, according to Hopf boundary lemma, it follows that
\begin{equation*}
\frac{\partial U}{\partial n}>0\text{~for~all~}(x,t)\in\{l_{1},l_{2}\}\times(\tau,T],
\end{equation*}
where $\partial / \partial n$ denotes the outward directional derivative. Therefore, there exists a sufficiently small $\epsilon_{0}>0$ such that $U(x,t)\geq\epsilon_{0} w_{1}(x,t)$ for $(x,t)\in[l_{1},l_{2}]\times[0,T]$. This shows a  contradiction with the assumption that $\mbox{(i)}$ holds.

\textbf{Case 2}: $U(x,t)=0$ holds for $(x,t)\in[l_{1},l_{2}]\times (\tau,T]$. This case is impossible for the reason that $U(x,\tau)>0$ for $x\in(l_{1},l_{2})$.

When $\mbox{(ii)}$ holds, we have that
\begin{eqnarray}\label{3-35}
\left\{
\begin{array}{ll}
V_{t}-d_{2}\Delta V+\gamma V=(\delta_{1}+a_{11}+a_{22} )V\geq 0 ,\; &\,x\in(l_{1},l_{2}), t\in(0^{+}, \tau],  \\[2mm]
V_{t}-d_{2}\Delta V+\gamma V=(\gamma-a_{22} )V> 0 ,\; &\,x\in(l_{1},l_{2}), t\in(\tau, T],  \\[2mm]
V(x,t)=0,\; &\, x\in\{l_{1},l_{2}\},  t\in(0^{+}, T],\\[2mm]
V(x,0^{+})\geq 0,\; &\, x\in[l_{1},l_{2}]. \\[2mm]
\end{array} \right.
\end{eqnarray}
Then, we can also conclude that $\mbox{(ii)}$ does not hold by an argument analogous to $\mbox{(i)}$. Therefore, we must have $\zeta_{0}=1$. Similarity, we can also obtain the conclusion that $(w_{2},z_{2})\succeq(w_{1},z_{1})$ for $(x,t)\in[l_{1},l_{2}]\times[0,T]$. Now, we complete the proof of the first assertion.

Next, we prove the second assertion. We periodically extend the constructed pairs $(\overline{w}(x,t), \overline{z}(x,t))$ and $(\underline{w}(x,t), \underline{z}(x,t))$ to $\mathcal{R}^{+}$ such that $(\overline{w}(x,t+T), \overline{z}(x,t+T))=(\overline{w}(x,t), \overline{z}(x,t))$ and $(\underline{w}(x,t+T), \underline{z}(x,t+T))=(\underline{w}(x,t), \underline{z}(x,t))$ for $x\in[l_{1},l_{2}]$ and $t\in \mathcal{R}^{+}$. We still use $(\overline{w}(x,t), \overline{z}(x,t))$ and $(\underline{w}(x,t), \underline{z}(x,t))$ to denote such extended pairs. The standard calculation yields that the pair $(\overline{w}(x,t), \overline{z}(x,t))$ is the upper solution of problem \eqref{3-1}. Since $(u_{0}(x), v_{0}(x))\succ (0,0)$ for all $x\in(l_{1},l_{2})$, we have that the pair $(\underline{w}(x,t), \underline{z}(x,t))$ is the lower solution of problem \eqref{3-1} by decreasing the $\epsilon$ if necessary.

By choosing the pairs $(\overline{w}(x,t), \overline{z}(x,t))$ and $(\underline{w}(x,t), \underline{z}(x,t))$ defined above as the initial iteration of \eqref{3-33} with $m=0$ replaced by $m=0, 1,\cdots$, the iteration sequences $\big\{(\overline{w}^{i}, \overline{z}^{i})\big\}_{i=0}^{n}$ and $\big\{(\underline{w}^{i}, \underline{z}^{i})\big\}_{i=0}^{n}$ can be obtained, respectively. We claim that
\begin{equation}\label{3-36}
\big(\underline{w}^{n}(x,t), \underline{z}^{n}(x,t)\big)\preceq\big(u(x, t+nT), v(x, t+nT)\big)\preceq\big(\overline{w}^{n}(x,t), \overline{z}^{n}(x,t)\big)
\end{equation}
for $(x,t)\in[l_{1},l_{2}]\times[0,\infty)$. First, we note that \eqref{3-36} holds for $n=0$ for the reason that the pairs $(\overline{w}(x,t), \overline{z}(x,t))$ and $(\underline{w}(x,t), \underline{z}(x,t))$ are the ordered upper and lower solutions of problem \eqref{3-1}, respectively. We assume that \eqref{3-36} holds for $n=k$, that is,
\begin{equation}\label{3-37}
\big(\underline{w}^{k}(x,t), \underline{z}^{k}(x,t)\big)\preceq\big(u(x, t+kT), v(x, t+kT)\big)\preceq\big(\overline{w}^{k}(x,t), \overline{z}^{k}(x,t)\big)
\end{equation}
for $(x,t)\in[l_{1},l_{2}]\times[0,\infty)$. Next, we will prove that \eqref{3-36} holds for $n=k+1$. It is easily seen from \eqref{3-37} that
\begin{equation*}
\big(\underline{w}^{k+1}(x,t), \underline{z}^{k+1}(x,t)\big)\preceq\big(u(x, t+(k+1)T), v(x, t+(k+1)T)\big)\preceq\big(\overline{w}^{k+1}(x,t), \overline{z}^{k+1}(x,t)\big)
\end{equation*}
for $(x,t)\in[l_{1},l_{2}]\times\{0, (mT)^{+} \}$.  Then, the pairs $\big(\overline{w}^{k+1}(x,t), \overline{z}^{k+1}(x,t)\big)$ and $\big(\underline{w}^{k+1}(x,t), \underline{z}^{k+1}(x,t)\big)$ are the ordered upper and lower solutions of $\big(u(x, t+(k+1)T), v(x, t+(k+1)T)\big)$ by a standard calculation.   Therefore, it follows from \autoref{lemma 3.0.2} and \autoref{remark 3.0.1} that \eqref{3-36} holds for $n=k+1$. In virtue of the proof of the first assertion, the desired \eqref{3-38} holds for $(x,t)\in[l_{1},l_{2}]\times[0,T]$. Noting that the unique solution $(w,z)$ of problem \eqref{3-32} is periodic, the desired \eqref{3-38} also holds for $(x,t)\in[l_{1},l_{2}]\times[0,\infty)$. This ends the proof of (1).

(2) We first prove the first assertion. It is easily seen that $(0,0)$ is the nonnegative solution of the steady state problem \eqref{3-32}. We will prove that it is the unique nonnegative solution. Arguing by contradiction, we may suppose that there would exist a pair $(w,z)$ being a strongly positive solution of problem \eqref{3-32}. It follows from problem \eqref{3-32} that
\begin{eqnarray*}
\left\{
\begin{array}{ll}
w_{t}=-\delta_{1} w ,\; &\,x\in(l_{1},l_{2}), t\in(0, \tau],  \\[2mm]
z_{t}=d_{2}\Delta z-\delta_{2} z ,\; &\,x\in(l_{1},l_{2}), t\in(0, \tau],  \\[2mm]
w_{t}=d_{1}\Delta w-a_{11}w+a_{12}z ,\; &\,x\in(l_{1},l_{2}), t\in(\tau, T],  \\[2mm]
z_{t}<d_{2}\Delta z-a_{22}z+f'(0)w ,\; &\,x\in(l_{1},l_{2}), t\in(\tau, T],  \\[2mm]
w(x,t)=z(x,t)=0,\; &\, x\in\{l_{1},l_{2}\},  t\in[0, T],\\[2mm]
w(x,0)\leq H'(0)w(x,T),\; &\, x\in[l_{1},l_{2}],  \\[2mm]
z(x,0)=z(x,T),\; &\, x\in[l_{1},l_{2}].
\end{array} \right.
\end{eqnarray*}
Then, \autoref{lemma 3.1.1} yields that $\lambda_{1}< 0$, which contradicts with $\lambda_{1}\geq 0$. This completes the proof of the first assertion.

Next, we prove the second assertion. Its proof is divided into two cases.

\textbf{Case 1}: $\lambda_{1}> 0$

Let $(\phi_{1}, \psi_{1})$ be the unique eigenfunction pair corresponding to the principal eigenvalue $\lambda_{1}$ of problem \eqref{3-3}. We periodically extend the pair $(\phi_{1}(x,t), \psi_{1}(x,t))$ to $\mathcal{R}^{+}$ such that $(\phi_{1}(x,t+T), \psi_{1}(x,t+T))=(\phi_{1}(x,t), \psi_{1}(x,t))$ for $x\in[l_{1},l_{2}]$ and $t\in \mathcal{R}^{+}$. We use $(\phi_{2}, \psi_{2})$ to denote such extended pair. Now, we construct the following functions
\begin{equation*}
\overline{u}(x,t)=Me^{-\lambda_{1}t}\phi_{2}(x,t)\text{~and~}\overline{v}(x,t)=Me^{-\lambda_{1}t}\psi_{2}(x,t)
\end{equation*}
for $(x,t)\in[l_{1},l_{2}]\times[0,\infty)$, where $M>0$ is a sufficiently large constant such that $(M\phi_{2}(x,0), M\psi_{2}(x,0))\succeq (u_{0}(x),v_{0}(x) )$ for any $x\in[l_{1},l_{2}]$. By a standard calculation, it follows from \autoref{lemma 3.0.2} that
\begin{equation*}
\lim\limits_{t\rightarrow+\infty}(u(x,t), v(x,t)) \leq\lim\limits_{t\rightarrow+\infty}(\overline{u}(x,t), \overline{v}(x,t))=(0,0)
\end{equation*}
for $x\in[l_{1},l_{2}]$ and $t\in \mathcal{R}^{+}$.

\textbf{Case 2}: $\lambda_{1}=0$

The direct calculation gives that the extended pair $(\overline{w}(x,t), \overline{z}(x,t))$ defined in (1) is an upper solution of problem \eqref{3-1}. By choosing $(\overline{w}(x,t), \overline{z}(x,t))$ as an initial iteration, a monotonically decreasing sequence $\big\{(\overline{w}^{i}, \overline{z}^{i})\big\}_{i=0}^{n}$ can be obtained with the help of the iteration rule \eqref{3-33} with $m=0$ replaced by $m=0, 1,\cdots$. Similar to the discussion of the second assertion in (1), we have that
\begin{equation}\label{3-40}
(0, 0)\preceq\big(u(x, t+kT), v(x, t+kT)\big)\preceq\big(\overline{w}^{k}(x,t), \overline{z}^{k}(x,t)\big)
\end{equation}
for $(x,t)\in[l_{1},l_{2}]\times[0,\infty)$. Recalling that $(0,0)$ is the only nonnegative solution of the steady state problem \eqref{3-32}, it follows from \eqref{3-40} that
\begin{equation*}
\lim\limits_{k\rightarrow+\infty}\big(u(x,t+kT),v(x,t+kT)\big)=(0,0) \text{~for~all~}(x,t)\in[l_{1},l_{2}]\times[0,+\infty).
\end{equation*}
This implies that
\begin{equation*}
\lim\limits_{t\rightarrow+\infty}(u(x,t),v(x,t))=(0,0) \text{~uniformly~for~}x\in[l_{1},l_{2}].
\end{equation*}
The proof is completed.
\end{proof}
\end{theorem}
We next present the ordinary differential equation problem corresponding to the fixed boundary problem \eqref{3-1}, which has the form
\begin{eqnarray}\label{3-41}
\left\{
\begin{array}{ll}
U_{t}=-\delta_{1} U ,\; &\, t\in((mT)^{+}, mT+\tau],  \\[2mm]
V_{t}=-\delta_{2} V ,\; &\,t\in((mT)^{+}, mT+\tau],  \\[2mm]
U_{t}=-a_{11}U+a_{12}V ,\; &\,t\in(mT+\tau, (m+1)T],  \\[2mm]
V_{t}=-a_{22}V+f(U) ,\; &\, t\in(mT+\tau, (m+1)T],
\end{array} \right.
\end{eqnarray}
with the impulsive condition $(U((mT)^{+}),V((mT)^{+}))=(H(U(mT)),V(mT))$ and the initial function $(U(0), V(0))=(U_{0}, V_{0})$. The steady state problem
corresponding to problem \eqref{3-41} takes the following form
\begin{eqnarray}\label{3-42}
\left\{
\begin{array}{ll}
W_{t}=-\delta_{1} W ,\; &\, t\in(0^{+}, \tau],  \\[2mm]
Z_{t}=-\delta_{2} Z ,\; &\,t\in(0^{+}, \tau],  \\[2mm]
W_{t}=-a_{11}W+a_{12}Z ,\; &\,t\in(\tau, T],  \\[2mm]
Z_{t}=-a_{22}Z+f(W) ,\; &\, t\in(\tau, T],  \\[2mm]
W(0^{+})=H(W(0)),\; &\,  \\[2mm]
Z(0^{+})=Z(0),\; &\,  \\[2mm]
W(0)=W(T),Z(0)=Z(T).\; &\
\end{array} \right.
\end{eqnarray}
For problem \eqref{3-41}, we have a result analogous to \autoref{theorem 3.2}. We omit the proof here as the method is similar.
\begin{lemma}\label{lemma 3.1.5}
Let $\nu_{1}$ be the principal eigenvalue of problem \eqref{3-31}. Denote by $(U,V)$ and $(W,Z)$ the solutions to \eqref{3-41} and \eqref{3-42}, respectively. Then, the following assertions are valid:
\begin{enumerate}
\item[$(1)$]
If $\nu_{1}< 0$, then the steady state problem \eqref{3-42} admits a unique positive periodic solution $(W,Z)$. Additionally, $(U,V)$ satisfies
\begin{equation*}
\lim\limits_{n\rightarrow+\infty}\big(U(t+nT),V(t+nT)\big)=\big(W(t), Z(t)\big)
\end{equation*}
uniformly for $t\in[0, +\infty)$;
\item[$(2)$]
If $\nu_{1}\geq 0$, then $(0,0)$ is the only nonnegative solution of the steady state problem \eqref{3-42}. Additionally, $(U,V)$ satisfies
\begin{equation*}
\lim\limits_{t\rightarrow+\infty}\big(U(t),V(t)\big)=\big(0, 0\big).
\end{equation*}
\end{enumerate}
\end{lemma}
\autoref{theorem 3.2} and \autoref{lemma 3.1.5} above demonstrate the long-time dynamical behavior of problem \eqref{3-1} and its corresponding ordinary differential equation problem \eqref{3-41}, respectively. As the end of this subsection, we present the relationship between the steady state solutions of these two problems.
\begin{lemma}\label{lemma 3.1.6}
Suppose that $[l_{1},l_{2}]=[-s, s]$. Let $\nu_{1}$ and $(W(t),Z(t))$ be the principal eigenvalue of problem \eqref{3-31} and the solution to problem \eqref{3-42}, respectively. If $\nu_{1}<0$, then
\begin{equation*}
\lim\limits_{s\rightarrow\infty}(w_{s}(x,t),z_{s}(x,t))=(W(t),Z(t))
\end{equation*}
locally uniformly for $x\in\mathds{R}$ and uniformly for $t\in[0,T]$, where $(w_{s}(x,t),z_{s}(x,t))$ is the steady state solution to problem \eqref{3-32} with $[l_{1},l_{2}]$
replaced by $[-s, s]$.
\begin{proof}
Due to $\nu_{1}<0$, we have from \autoref{lemma 3.1.7} and (2) of \autoref{lemma 3.1.3} that there exists some positive constant $s^{*}$ such that $\lambda_{1}(s)<0$ holds for all $s>s^{*}$. Then, it follows from (1) of \autoref{theorem 3.2} that  problem \eqref{3-32} admits a unique positive periodic solution $(w_{s},z_{s})$ for $x\in[-s, s]$ and $t\in[0,T]$ when $s>s^{*}$.

Firstly, we claim that $(w_{s},z_{s})$ is nondecreasing with respect to $s\in(s^{*}, \infty)$. We suppose that $s_{2}>s_{1}>s^{*}$. Let $(u_{01},v_{01} )$
and $(u_{02},v_{02} )$ be the initial functions to problem \eqref{3-1} with $[l_{1},l_{2}]$ replaced by $[-s_{1}, s_{1}]$ and $[-s_{2}, s_{2}]$, respectively, and they satisfies
\begin{equation*}
0<u_{01}(x)\leq u_{02}(x)\leq C_{2} \text{~and~}0<v_{01}(x)\leq v_{02}(x)\leq C_{3}
\end{equation*}
for $x\in(-s_{1}, s_{1})$. We denote by $ (u_{s_{i}}, v_{s_{i}})$ the corresponding  solution of problem \eqref{3-1} with initial function $(u_{0i},v_{0i} )$ for $i=1,2$. Then, it follows from \autoref{lemma 3.0.2} and \autoref{remark 3.0.1} that
\begin{equation*}
(u_{s_{1}}, v_{s_{1}})\preceq (u_{s_{2}}, v_{s_{2}})\preceq(U, V)
\end{equation*}
for $x\in[-s_{1}, s_{1}]$ and $t\geq 0$, where $(U, V)$ is the solution to problem \eqref{3-41} with initial function $(U_{0}, V_{0})=(C_{2}, C_{3})$. We have from  (1) of \autoref{theorem 3.2} and \autoref{lemma 3.1.5} that
\begin{equation}\label{3-43}
(w_{s_{1}}, z_{s_{1}})\preceq (w_{s_{2}}, z_{s_{2}})\preceq(W, Z)
\end{equation}
for $x\in[-s_{1}, s_{1}]$ and $t\in[0, T]$.

By using parabolic regularity theory and a diagonal procedure, it is easily seen that $(w_{s}(x,t), z_{s}(x,t))$ converges to some $(\widetilde{w}(x,t), \widetilde{z}(x,t))$, in $\mathcal{C}^{2}_{\text{loc}}(\mathcal{R}\times [0,T])$, as  $s\rightarrow +\infty$, where $(\widetilde{w}(x,t), \widetilde{z}(x,t))$ satisfies
\begin{eqnarray*}
\left\{
\begin{array}{ll}
\widetilde{w}_{t}=-\delta_{1} \widetilde{w} ,\; &\,x\in\mathcal{R}, t\in(0^{+}, \tau],  \\[2mm]
\widetilde{z}_{t}=d_{2}\Delta \widetilde{z}-\delta_{2} \widetilde{z} ,\; &\,x\in\mathcal{R}, t\in(0^{+}, \tau],  \\[2mm]
\widetilde{w}_{t}=d_{1}\Delta \widetilde{w}-a_{11}\widetilde{w}+a_{12}\widetilde{z} ,\; &\,x\in\mathcal{R}, t\in(\tau, T],  \\[2mm]
\widetilde{z}_{t}=d_{2}\Delta \widetilde{z}-a_{22}\widetilde{z}+f(\widetilde{w}) ,\; &\,x\in\mathcal{R}, t\in(\tau, T],  \\[2mm]
\widetilde{w}(x,0^{+})=H(\widetilde{w}(x,0)),\; &\, x\in\mathcal{R},   \\[2mm]
\widetilde{z}(x,0^{+})=\widetilde{z}(x,0),\; &\, x\in\mathcal{R},   \\[2mm]
\widetilde{w}(x,0)=\widetilde{w}(x,T),\widetilde{z}(x,0)=\widetilde{z}(x,T),\; &\,x\in\mathcal{R}.
\end{array} \right.
\end{eqnarray*}
Then, it follows from \eqref{3-43} that $(0,0)\prec (\widetilde{w}(x,t), \widetilde{z}(x,t))\preceq (W(t), Z(t)) $ in $\mathcal{R}\times [0,T]$.

We next prove that the pair $(\widetilde{w}(x,t), \widetilde{z}(x,t))$ is independent of $x$. It suffices to show that for any given $x_{0}\in \mathcal{R}$, we have that $(\widetilde{w}(x_{0},t), \widetilde{z}(x_{0},t))=(\widetilde{w}(0,t), \widetilde{z}(0,t))$. Define $\tilde{s}=:|x_{0}|$, and
\begin{equation*}
(\widehat{w}(x,t), \widehat{z}(x,t))=:(w_{s-\tilde{s}}(x+x_{0},t), z_{s-\tilde{s}}(x+x_{0},t))
\end{equation*}
for $s>2\tilde{s}+s^{*}$. Noticing that
\begin{equation*}
[-s+2\tilde{s}, s-2\tilde{s}]\subset[-s+\tilde{s}-x_{0}, s-\tilde{s}-x_{0}]\subset [-s, s],
\end{equation*}
it follows from the claim above that
\begin{equation*}
(w_{s-2\tilde{s}}, z_{s-2\tilde{s}})\preceq (\widehat{w}(x,t), \widehat{z}(x,t)) \preceq(w_{s}, z_{s}).
\end{equation*}
Letting $s\rightarrow \infty$, we have that
\begin{equation}\label{3-44}
(\widetilde{w}(x,t), \widetilde{z}(x,t))\preceq (\widetilde{w}(x+x_{0},t), \widetilde{z}(x+x_{0},t)) \preceq(\widetilde{w}(x,t), \widetilde{z}(x,t)).
\end{equation}
Taking $x=0$, it follows from \eqref{3-44} that $(\widetilde{w}(x_{0},t), \widetilde{z}(x_{0},t))=(\widetilde{w}(0,t), \widetilde{z}(0,t))$.

Since $\nu_{1}<0$, it follows from (1) of \autoref{lemma 3.1.5} that $(W(t), Z(t))$ is the unique positive periodic solution to problem \eqref{3-42}, which is independent of $x$. Therefore, we must have $(\widetilde{w}(x,t), \widetilde{z}(x,t))=(W(t), Z(t)) $. This ends the proof.
\end{proof}
\end{lemma}
\section{Spreading-vanishing dichotomy and criteria}\label{Section-4}
This section shows that model \eqref{Fixed-1}-\eqref{Fixed-3} admits a dichotomy between spreading and vanishing, and derives a sharp criterion determining which outcome occurs. Throughout this section, we always denote the solution of model \eqref{Fixed-1}-\eqref{Fixed-3} by $(u,v,r,s)$ or $(u(x,t),v(x,t),r(t),s(t))$. Biologically, we assume that when $\nu_{1}=0$, the region at infinity is a high-risk environment for both the infectious agents and the infective human population, that is,
\begin{equation*}
(\textbf{G}): \text{If}~ \nu_{1}=0, \text{~then~} s_{\infty}-r_{\infty}<\infty,
\end{equation*}
where $\nu_{1}$ is defined in \autoref{lemma 3.1.7}.
\subsection{Spreading-vanishing dichotomy}\label{Section-4-1}
\begin{lemma}\label{lemma 4-1}
Let $\lambda_{1}(r_{\infty}, s_{\infty})$ be the principal eigenvalue of problem \eqref{3-2} with $l_{1}$ and $l_{2}$ replaced by $r_{\infty}$ and $s_{\infty}$. If $s_{\infty}-r_{\infty}<\infty$, then $\lambda_{1}(r_{\infty}, s_{\infty})\geq 0$ and
\begin{equation*}
\lim\limits_{t\rightarrow \infty}\|u(\cdot, t)\|_{\mathcal{C}([r(t), s(t)])}=\lim\limits_{t\rightarrow \infty}\|v(\cdot, t)\|_{\mathcal{C}([r(t), s(t)])}=0.
\end{equation*}
\begin{proof}
We first prove the first assertion. Arguing indirectly, we assume that $\lambda_{1}(r_{\infty}, s_{\infty})<0$. It follows from (2) of \autoref{lemma 3.1.3}
that there exists a sufficiently small constant $\epsilon$ such that $\lambda_{1}(r_{\infty}+\epsilon, s_{\infty}-\epsilon)<0$. For this fixed constant $\epsilon$,
there exists a sufficiently large integer $m_{0}$ such that $\lambda_{1}(r(t), s(t))<0$ and $[r(t), s(t)]\supset[r_{\infty}+\epsilon, s_{\infty}-\epsilon]$ for $t\geq m_{0}T$.

We now consider the following initial-boundary value problem
\begin{eqnarray}\label{4-29}
\left\{
\begin{array}{ll}
\underline{u}^{\epsilon}_{t}=-\delta_{1} \underline{u}^{\epsilon},\; &\,x\in\Omega_{1}(\epsilon), t\in((mT)^{+}, mT+\tau],  \\[2mm]
\underline{v}^{\epsilon}_{t}=d_{2}\Delta \underline{v}^{\epsilon}-\delta_{2} \underline{v}^{\epsilon} ,\; &\,x\in\Omega_{1}(\epsilon), t\in((mT)^{+}, mT+\tau],  \\[2mm]
\underline{u}^{\epsilon}_{t}=d_{1}\Delta \underline{u}^{\epsilon}-a_{11}\underline{u}^{\epsilon}+a_{12}\underline{v}^{\epsilon} ,\; &\,x\in\Omega_{1}(\epsilon), t\in(mT+\tau, (m+1)T],  \\[2mm]
\underline{v}^{\epsilon}_{t}=d_{2}\Delta \underline{v}^{\epsilon}-a_{22}\underline{v}^{\epsilon}+f(\underline{u}^{\epsilon}) ,\; &\,x\in\Omega_{1}(\epsilon), t\in(mT+\tau, (m+1)T],  \\[2mm]
\underline{u}^{\epsilon}(x,t)=\underline{v}^{\epsilon}(x,t)=0,\; &\, x\in\partial \Omega_{1}(\epsilon),  t\geq m_{0}T,\\[2mm]
\underline{u}^{\epsilon}(x,(mT)^{+})=H(\underline{u}^{\epsilon}(x,mT)),\; &\, x\in\Omega_{1}(\epsilon),  \\[2mm]
\underline{v}^{\epsilon}(x,(mT)^{+})=\underline{v}^{\epsilon}(x,mT),\; &\, x\in\Omega_{1}(\epsilon),  \\[2mm]
\underline{u}^{\epsilon}(m_{0}T,x)=u(m_{0}T, x),\underline{v}^{\epsilon}(m_{0}T,x)=v(m_{0}T, x),\; &\,x\in\overline{\Omega}_{1}(\epsilon), m=m_{0}, m_{0}+1, \cdots,
\end{array} \right.
\end{eqnarray}
where $\Omega_{1}(\epsilon)=(r_{\infty}+\epsilon, s_{\infty}-\epsilon)$. Then, it follows from \autoref{remark 2.1} that
\begin{equation}\label{4-1}
\underline{u}^{\epsilon}(x,t)\leq u(x,t)\text{~and~} \underline{v}^{\epsilon}(x,t)\leq v(x,t)\text{~for~} x\in\overline{\Omega}_{1}(\epsilon) \text{~and~}t\in[m_{0}T,\infty).
\end{equation}
Since $\lambda_{1}(\Omega_{1}(\epsilon))<0$, we have from (1) of \autoref{theorem 3.2} that
\begin{equation}\label{4-2}
\lim\limits_{m\rightarrow+\infty}\big(\underline{u}^{\epsilon}(x,t+mT),\underline{v}^{\epsilon}(x,t+mT)\big)=\big(\underline{w}^{\epsilon}(x,t), \underline{z}^{\epsilon}(x,t)\big)
\end{equation}
uniformly for $(x,t)\in\overline{\Omega}_{1}(\epsilon) \times[0, +\infty)$, where $\big(\underline{w}^{\epsilon}(x,t), \underline{z}^{\epsilon}(x,t)\big)$ denotes the solution of  problem \eqref{3-32} with $l_{1}$ and $l_{2}$ replaced by $r_{\infty}+\epsilon$ and $s_{\infty}-\epsilon$. It follows from \eqref{4-1} and \eqref{4-2} that
\begin{equation}\label{4-5}
\liminf\limits_{m\rightarrow+\infty}\big(u(x,t+mT),v(x,t+mT)\big)\geq \big(\underline{w}^{\epsilon}(x,t), \underline{z}^{\epsilon}(x,t)\big)
\end{equation}
uniformly for $(x,t)\in\overline{\Omega}_{1}(\epsilon) \times[0, +\infty)$.

Next, we consider the following problem
\begin{eqnarray*}
\left\{
\begin{array}{ll}
\overline{u}_{t}=-\delta_{1} \overline{u},\; &\,x\in\Omega_{2}, t\in((mT)^{+}, mT+\tau],  \\[2mm]
\overline{v}_{t}=d_{2}\Delta \overline{v}-\delta_{2} \overline{v} ,\; &\,x\in\Omega_{2}, t\in((mT)^{+}, mT+\tau],  \\[2mm]
\overline{u}_{t}=d_{1}\Delta \overline{u}-a_{11}\overline{u}+a_{12}\overline{v} ,\; &\,x\in\Omega_{2}, t\in(mT+\tau, (m+1)T],  \\[2mm]
\overline{v}_{t}=d_{2}\Delta \overline{v}-a_{22}\overline{v}+f(\overline{u}) ,\; &\,x\in\Omega_{2}, t\in(mT+\tau, (m+1)T],  \\[2mm]
\overline{u}(x,t)=\overline{v}(x,t)=0,\; &\, x\in\partial \Omega_{2},  t\geq 0,\\[2mm]
\overline{u}(x,(mT)^{+})=H(\overline{u}(x,mT)),\; &\, x\in\Omega_{2},  \\[2mm]
\overline{v}(x,(mT)^{+})=\overline{v}(x,mT),\; &\, x\in\Omega_{2} , m=0, 1, \cdots \\[2mm]
\end{array} \right.
\end{eqnarray*}
with the initial functions
\begin{eqnarray*}
\overline{u}(x,0)=
\left\{
\begin{array}{ll}
u_{0}(x),\; &\,-s_{0}\leq x\leq s_{0}, \\[2mm]
0,\; &\, x\in\overline{\Omega}_{2}\setminus[-s_{0}, s_{0}]
\end{array} \right.
\text{~and~~}
\overline{v}(x,0)=
\left\{
\begin{array}{ll}
v_{0}(x),\; &\, -s_{0}\leq x\leq s_{0}, \\[2mm]
0,\; &\, x\in\overline{\Omega}_{2}\setminus[-s_{0}, s_{0}],
\end{array} \right.
\end{eqnarray*}
where $\Omega_{2}=(r_{\infty}, s_{\infty})$. Then, it follows from \autoref{lemma 2-1} that
\begin{equation}\label{4-3}
u(x,t)\leq \overline{u}(x,t)\text{~and~} v(x,t)\leq \overline{v}(x,t)\text{~for~} x\in\overline{\Omega}_{2} \text{~and~}t\geq 0.
\end{equation}
Since $\lambda_{1}(\Omega_{2})<0$, we have from (1) of \autoref{theorem 3.2} that
\begin{equation}\label{4-4}
\lim\limits_{m\rightarrow+\infty}\big(\overline{u}(x,t+mT),\overline{v}(x,t+mT)\big)=\big(\overline{w}(x,t), \overline{z}(x,t)\big)
\end{equation}
uniformly for $(x,t)\in\overline{\Omega}_{2} \times[0, +\infty)$, where $\big(\overline{w}(x,t), \overline{z}(x,t)\big)$ is the solution of  problem \eqref{3-32} with $l_{1}$ and $l_{2}$ replaced by $r_{\infty}$ and $s_{\infty}$. It follows from \eqref{4-3} and \eqref{4-4} that
\begin{equation}\label{4-6}
\limsup\limits_{m\rightarrow+\infty}\big(u(x,t+mT),v(x,t+mT)\big)\leq \big(\overline{w}(x,t), \overline{z}(x,t)\big)
\end{equation}
uniformly for $(x,t)\in\overline{\Omega}_{2} \times[0, +\infty)$. From a standard compactness and uniqueness argument, it easily follows that
\begin{equation}\label{4-7}
\lim\limits_{\epsilon\rightarrow 0}\big(\underline{w}^{\epsilon}(x,t), \underline{z}^{\epsilon}(x,t)\big)=  \big(\overline{w}(x,t), \overline{z}(x,t)\big)
\end{equation}
in $\mathcal{C}^{1+\alpha, (1+\alpha)/2}_{\text{loc}}(\Omega_{2}\times [0,T];\mathcal{R}^{2} )$ for any $\alpha\in(0,1)$. By \eqref{4-5},  \eqref{4-6}, and the arbitrariness
of $\epsilon$, we have that
\begin{equation*}
\lim\limits_{m\rightarrow+\infty}\big(u(x,t+mT),v(x,t+mT)\big)= \big(\overline{w}(x,t), \overline{z}(x,t)\big)
\end{equation*}
uniformly for $(x,t)\in\Omega_{2} \times[0, T]$. This and the regularity of $u(x,t)$ and $v(x,t)$ also indicate that
\begin{equation*}
\lim\limits_{m\rightarrow+\infty}\big(u(x,t+mT),v(x,t+mT)\big)= \big(\overline{w}(x,t), \overline{z}(x,t)\big)
\end{equation*}
in $\mathcal{C}^{1+\alpha, (1+\alpha)/2}_{\text{loc}}(\Omega_{2}\times [0,T];\mathcal{R}^{2} )$ for any $\alpha\in(0,1)$. Therefore, we have that
\begin{equation}\label{4-8}
\lim\limits_{m\rightarrow+\infty}\big(u_{x}(s(t+mT), t+mT),v_{x}(s(t+mT), t+mT)\big)=\big(\overline{w}_{x}(s_{\infty},t), \overline{z}_{x}(s_{\infty},t)\big)<(0,0)
\end{equation}
and
\begin{equation}\label{4-9}
\lim\limits_{m\rightarrow+\infty}\big(u_{x}(r(t+mT), t+mT),v_{x}(r(t+mT), t+mT)\big)=\big(\overline{w}_{x}(r_{\infty},t), \overline{z}_{x}(r_{\infty},t)\big)>(0,0)
\end{equation}
uniformly on $[0,T]$. Recall that
\begin{equation*}
r'(t+mT)=-\mu_{1}u_{x}(r(t+mT),t+mT)-\mu_{2}v_{x}(r(t+mT),t+mT)
\end{equation*}
and
\begin{equation*}
s'(t+mT)=-\mu_{1}u_{x}(s(t+mT),t+mT)-\mu_{2}v_{x}(s(t+mT),t+mT)
\end{equation*}
for $t\in(\tau, T)$. This, together with \eqref{4-8} and \eqref{4-9}, yields that $s_{\infty}-r_{\infty}=\infty$, which show a contradiction with $s_{\infty}-r_{\infty}<\infty$.
Consequently, the first assertion holds.

Note that $\eqref{4-3}$ remains true. Since $s_{\infty}-r_{\infty}<\infty$, it follows from (2) of \autoref{theorem 3.2} that
\begin{equation*}
\lim\limits_{t\rightarrow+\infty}\big(u(x,t),v(x,t)\big)=\big(0, 0\big) \text{~uniformly~for~}x\in[r_{\infty},s_{\infty}].
\end{equation*}
Hence, we have
\begin{equation*}
\lim\limits_{t\rightarrow \infty}\|u(\cdot, t)\|_{\mathcal{C}([r(t), s(t)])}=\lim\limits_{t\rightarrow \infty}\|v(\cdot, t)\|_{\mathcal{C}([r(t), s(t)])}=0.
\end{equation*}
This ends the proof.
\end{proof}
\end{lemma}

Next, the long-time dynamic behavior produced by model \eqref{Fixed-1}-\eqref{Fixed-3} is discussed for each of the cases $\lambda_{1}(\infty)<0$, $\lambda_{1}(\infty)=0$, and $\lambda_{1}(\infty)>0$.
\begin{theorem}\label{theorem 4-1}
If $\lambda_{1}(\infty)>0$, then vanishing happens.
\begin{proof}
To derive the desired assertion, an upper solution will be constructed. Define
\begin{equation*}
\overline{s}(t)=\hat{s}(1+\eta-\frac{\eta}{2}e^{-\gamma t}) \text{~for~}t\geq0
\end{equation*}
and
\begin{equation}\label{4-42}
\overline{u}(x,t)=C_{0}e^{-\gamma t}\phi_{1}\bigg(\frac{x}{\tau(t)}, t\bigg)\text{~and~}\overline{v}(x,t)=C_{0}e^{-\gamma t}\psi_{1}\bigg(\frac{ x}{\tau(t)}, t\bigg)
\end{equation}
for $t\in [0, \infty)$ and $x\in [-\overline{s}(t), \overline{s}(t)]$, where $\tau(t)=1+\eta-\frac{\eta}{2}e^{-\gamma t}$. In the above $\hat{s}(\geq s_{0})$, $\eta$, $\gamma $, and $C_{0}$ are positive constants to be determined later, and $(\phi_{1}(x,t), \psi_{1}(x,t))$ is the strongly positive eigenfunction pair to the principle eigenvalue $\lambda_{1}(-\hat{s}, \hat{s})$ in problem \eqref{3-2} with $l_{1}$ and $l_{2}$ replaced by $-\hat{s}$ and $\hat{s}$. By using the method of separation of variables, $\phi_{1}$ and $\psi_{1}$ can be future written as
\begin{equation}\label{4-10}
\phi_{1}(x,t)=\chi^{\hat{s}}_{1}(x)\Phi_{\hat{s}}(t) \text{~and~} \psi_{1}(x,t)=\chi^{\hat{s}}_{1}(x)\Psi_{\hat{s}}(t),
\end{equation}
where $\chi^{\hat{s}}_{1}(x)$ is the strongly positive eigenfunction to the principle eigenvalue $\kappa^{\hat{s}}_{1}$ in problem \eqref{3-30} with $l_{1}$ and $l_{2}$ replaced by $-\hat{s}$ and $\hat{s}$, and $(\Phi_{\hat{s}}(t), \Psi_{\hat{s}}(t))$ is the positive eigenfunction pair to the principle eigenvalue $\lambda_{1}(-\hat{s}, \hat{s})$ in problem \eqref{3-15} with $\kappa_{1}$ replaced by $\kappa^{\hat{s}}_{1}$.

Next, we show that the quadruple $(\overline{u}, \overline{v}, -\overline{s}(t), \overline{s}(t))$ is an upper solution of model \eqref{Fixed-1}-\eqref{Fixed-3}. Fix an arbitrary integer $M\geq 1$. We first restrict $t\in((iT)^{+}, iT+\tau]$, where $i=0,1,2,\cdots, M-1$. Then, it follows from basic calculations that
\begin{equation*}
\begin{split}
\overline{u}_{t}+\delta_{1} \overline{u}=& C_{0}e^{-\gamma t}[-\gamma \phi_{1}+\phi_{1t} -x \tau^{-2}(t)\tau'(t)\phi_{1x}+\delta_{1} \phi_{1} ]\\
\geq & C_{0}e^{-\gamma t}[-\gamma \phi_{1}+\lambda_{1}(-\hat{s}, \hat{s})\phi_{1}-\delta_{1}\phi_{1} -x \tau^{-2}(t)\tau'(t)\phi_{1x}+\delta_{1} \phi_{1} ]\\
\geq & C_{0}e^{-\gamma t}[-\gamma \phi_{1}+\lambda_{1}(\infty)\phi_{1} -x \tau^{-2}(t)\tau'(t)\phi_{1x}]\\
\end{split}
\end{equation*}
for $x\in(-\overline{s}(t), \overline{s}(t))$. We notice that
\begin{equation}\label{4-11}
1+\frac{\eta}{2}\leq \tau(t)\leq 1+\eta\text{~and~}x \phi_{1x}=- \frac{\pi x\Phi_{\hat{s}}(t) }{2\hat{s}}\sin(\frac{\pi}{2\hat{s}}x) \leq 0.
\end{equation}
As a consequence, by \eqref{4-11} we find
\begin{equation}\label{4-12}
\overline{u}_{t}+\delta_{1} \overline{u}\geq 0
\end{equation}
for $t\in ((iT)^{+}, iT+\tau]$ and $x\in(-\overline{s}(t), \overline{s}(t))$, provided that
\begin{equation*}
\gamma\leq \gamma_{1}=:\frac{\lambda_{1}(\infty)}{2}.
\end{equation*}
Similarly, we have that
\begin{equation}\label{4-13}
\begin{split}
\overline{v}_{t}-d_{2}\Delta\overline{v}+\delta_{2} \overline{v} \geq & C_{0}e^{-\gamma t}\Big[\lambda_{1}(\infty)\psi_{1}-\gamma \psi_{1} +d_{2}(1-(1+\frac{\eta}{2})^{-2}) \Delta\psi_{1}  \Big]\\
\geq & C_{0}e^{-\gamma t}\chi^{\hat{s}}_{1}(x)\Psi_{\hat{s}}(t) \Big[\lambda_{1}(\infty)-\gamma -\frac{d_{2}\pi^{2}}{4s^{2}_{0}}(1-(1+\frac{\eta}{2})^{-2})  \Big]\\
\geq &0
\end{split}
\end{equation}
for $t\in ((iT)^{+}, iT+\tau]$ and $x\in(-\overline{s}(t), \overline{s}(t))$, provided that
\begin{equation*}
\gamma\leq \gamma_{1}\text{~and~}\eta\leq \eta_{1},
\end{equation*}
where $\eta_{1}$ is a sufficiently small positive constant such that
\begin{equation*}
1-(1+\frac{\eta}{2})^{-2}\leq\frac{\lambda_{1}(\infty)}{2}\frac{4s^{2}_{0}}{(d_{1}+d_{2})\pi^{2}}.
\end{equation*}

Secondly, we restrict $t\in( iT+\tau, (i+1)T]$, where $i=0,1,2,\cdots, M-1$. It follows from calculations analogous to \eqref{4-13} that
\begin{equation}\label{4-14}
\overline{u}_{t}-d_{1}\Delta \overline{u}+a_{11}\overline{u}-a_{12}\overline{v},~\overline{v}_{t}-d_{2}\Delta \overline{v}+a_{22}\overline{v}-f(\overline{u})\geq 0, \\
\end{equation}
for $t\in( iT+\tau, (i+1)T]$ and $x\in(-\overline{s}(t), \overline{s}(t))$, provided that
\begin{equation*}
\gamma\leq \gamma_{1}\text{~and~}\eta\leq \eta_{1}.
\end{equation*}
Therefore, \eqref{4-12}, \eqref{4-13}, and \eqref{4-14} hold if we take $\gamma=\frac{\gamma_{1}}{2}$ and $\eta=\frac{\eta_{1}}{2}$.

Thirdly, we restrict $t=(iT)^{+}$, where $i=0,1,2,\cdots, M-1$. By using assumption \textbf{(H)}, one can obtain that
\begin{equation*}
\begin{split}
\overline{u}(x, (iT)^{+})&=C_{0}e^{-\gamma iT}\phi_{1}\bigg(\frac{x}{\tau(iT)}, (iT)^{+}\bigg)=C_{0}H'(0)e^{-\gamma iT}\phi_{1}\bigg(\frac{x}{\tau(iT)}, iT\bigg)\\
&\geq H\Bigg(C_{0}e^{-\gamma iT}\phi_{1}\bigg(\frac{x}{\tau(iT)}, iT\bigg)\Bigg)=H(\overline{u}(x, iT))
\end{split}
\end{equation*}
for $x\in(-\overline{s}(t), \overline{s}(t))$. Similarly, we have that $\overline{v}(x, (iT)^{+})\geq \overline{v}(x, iT)$ for $x\in(-\overline{s}(t), \overline{s}(t))$.

Finally, we restrict $t=0$. \autoref{lemma 3.1.4} and the standard calculation yield that
\begin{equation*}
\begin{split}
\overline{u}(x,0)=C_{0}\phi_{1}\bigg(\frac{x}{\tau(0)}, 0\bigg)\geq C_{0}\chi^{\hat{s}}_{1}(\frac{x}{\tau(0)})\Phi_{\hat{s}}(0) \geq n_{1}C_{0}\chi^{s_{0}}_{1}(x)\geq u_{0}(x),\\
\overline{v}(x,0)=C_{0}\psi_{1}\bigg(\frac{x}{\tau(0)}, 0\bigg)\geq C_{0}\chi^{\hat{s}}_{1}(\frac{x}{\tau(0)})\Psi_{\hat{s}}(0) \geq n_{3}C_{0}\chi^{s_{0}}_{1}(x)\geq v_{0}(x)
\end{split}
\end{equation*}
for $x\in [-s_{0}, s_{0}]$, provided the sufficiently large $C_{0}$. For the boundary conditions, it is easily seen that
\begin{equation*}
\overline{u}(\pm \overline{s}(t),t)=C_{0}e^{-\gamma t}\phi_{1}(\pm \hat{s}, t)=0,~~\overline{v}(\pm \overline{s}(t),t)=C_{0}e^{-\gamma t}\psi_{1}(\pm \hat{s}, t)=0
\end{equation*}
for $t\in(0, MT]$.

In the preceding proof, constants $\eta$, $\gamma $, and $C_{0}$ are determined such that the defined function $(\overline{u},\overline{v})$ satisfies the differential equations, impulse conditions, and initial-boundary conditions required for an upper solution. We now proceed to determine constant $\hat{s}$ such that the free boundaries $-\overline{s}(t)$ and $\overline{s}(t)$ also satisfy the equation characterizing the upper solution. By using \eqref{4-11} and \autoref{lemma 3.1.4}, a direct calculation gives that
\begin{equation*}
\overline{s}'(t)=\frac{1}{2}\hat{s}\eta\gamma e^{-\gamma t}
\end{equation*}
and
\begin{equation*}
\begin{split}
-\mu_{1}\overline{u}_{x}(\overline{s}(t),t)-\mu_{2}\overline{v}_{x}(\overline{s}(t),t)=&\frac{C_{0}}{\tau(t)}e^{-\gamma t}[-\mu_{1}\phi_{1x}( \hat{s}, t)-\mu_{2}\psi_{1x}( \hat{s}, t)]\\
\leq &\frac{C_{0}\pi}{(2+\eta)\hat{s}}e^{-\gamma t}[\mu_{1}n_{2}+\mu_{2}n_{4}].
\end{split}
\end{equation*}
By setting
\begin{equation*}
\hat{s}=\max\Bigg\{s_{0}, \sqrt{\frac{2C_{0}\pi[\mu_{1}n_{2}+\mu_{2}n_{4}]}{(2+\eta)\eta\gamma}}\Bigg\},
\end{equation*}
it then follows that
\begin{equation*}
-\overline{s}'(t)\leq -\mu_{1}\overline{u}_{x}(-\overline{s}(t),t)-\mu_{2}\overline{v}_{x}(-\overline{s}(t),t)\text{~and~}\overline{s}'(t)\geq -\mu_{1}\overline{u}_{x}(\overline{s}(t),t)-\mu_{2}\overline{v}_{x}(\overline{s}(t),t)
\end{equation*}
for $t\in(iT+\tau, (i+1)T)$, where $i=0,1,2,\cdots, M-1$. For $t\in[iT, iT+\tau]$, it is easily seen that
\begin{equation*}
-\overline{s}(t)\leq -\overline{s}(iT)\text{~and~}\overline{s}(t)\geq\overline{s}(iT).
\end{equation*}

As a result, $(\overline{u},\overline{v}, -\overline{s}(t),\overline{s}(t))$ is an upper solution for model \eqref{Fixed-1}-\eqref{Fixed-3}, and we can apply \autoref{lemma 2-1} to assert that
\begin{equation*}
r(t)\geq -\overline{s}(t),~s(t)\leq \overline{s}(t)\text{~in~}[0, MT],
\end{equation*}
and
\begin{equation*}
u(x,t)\leq \overline{u}(x,t), v(x,t)\leq \overline{v}(x,t)\text{~for~}x\in[r(t), s(t)]\text{~and~}t\in[0, MT].
\end{equation*}
Note that $M$ is arbitrary and the constants $\hat{s}$, $\eta$, $\gamma $, and $C_{0}$ are independent of $M$. It then follows that
\begin{equation}\label{4-44}
r(t)\geq -\overline{s}(t), s(t)\leq \overline{s}(t)\text{~in~}[0, \infty),
\end{equation}
and
\begin{equation}\label{4-43}
u(x,t)\leq \overline{u}(x,t), v(x,t)\leq \overline{v}(x,t)\text{~for~}x\in[r(t), s(t)]\text{~and~}t\in[0, \infty).
\end{equation}
Therefore, we have that
\begin{equation*}
s_{\infty}-r_{\infty}\leq 2\lim\limits_{t\rightarrow\infty}\overline{s}(t)=2\hat{s}(1+\eta)
\end{equation*}
by using \eqref{4-44} and the construction of $\overline{s}(t)$, and
\begin{equation*}
\lim\limits_{t\rightarrow\infty}\|u(x,t)\|_{\mathcal{C}([r_{\infty}, s_{\infty}])}+\lim\limits_{t\rightarrow\infty}\|v(x,t)\|_{\mathcal{C}([r_{\infty}, s_{\infty}])}=0
\end{equation*}
by using \eqref{4-42} and \eqref{4-43}. This ends the proof.
\end{proof}
\end{theorem}

\begin{theorem}\label{theorem 4-2}
Assume that (\textbf{G}) holds. If $\lambda_{1}(\infty)=0$, then vanishing happens.
\begin{proof}
Let $(U, V)$ be the unique solution of problem \eqref{3-41} with $(U_{0}, V_{0})$ replaced by $(C_{2}, C_{3})$, where $C_{2}$ and $C_{3}$ are defined in \eqref{2-2}.
Then, it follows from the usual comparison principle that
\begin{equation}\label{4-15}
u\leq U\text{~and~}v\leq V\text{~for~}(x,t)\in \mathds{R}\times \mathds{R}^{+}.
\end{equation}
By using condition $\lambda_{1}(\infty)=0$, we have from \autoref{lemma 3.1.7} that $\nu_{1}=0$, which combined with \autoref{lemma 3.1.5} deduces that
\begin{equation}\label{4-16}
\lim\limits_{t\rightarrow+\infty}\big(U(t),V(t)\big)=\big(0, 0\big).
\end{equation}
Finally, \eqref{4-15} and \eqref{4-16} yield that
\begin{equation*}
\lim\limits_{t\rightarrow\infty}\|u(x,t)\|_{\mathcal{C}([r_{\infty}, s_{\infty}])}+\lim\limits_{t\rightarrow\infty}\|v(x,t)\|_{\mathcal{C}([r_{\infty}, s_{\infty}])}=0.
\end{equation*}
This proof is completed.
\end{proof}
\end{theorem}

When $\lambda_{1}(\infty)\geq0$, \textcolor{blue}{Theorems} \ref{theorem 4-1} and \ref{theorem 4-2} above demonstrate that vanishing occurs. Below we will discuss the case when $\lambda_{1}(\infty)<0$.
\begin{lemma}\label{lemma 4-2}
If $\lambda_{1}(\infty)<0$, then $r_{\infty}>-\infty$ if and only if $s_{\infty}<+\infty$.
\begin{proof}
Without loss of generalization, we may assume that $s_{\infty}<+\infty$. Let us prove $s_{\infty}-r_{\infty}< \infty$. Otherwise, $s_{\infty}-r_{\infty}= \infty$,
then there would exist an sufficiently large integer $m_{0}$ such that $\lambda_{1}(r(m_{0}T), s(m_{0}T))<0$.

Consider the following eigenvalue problem
\begin{eqnarray}\label{4-17}
\left\{
\begin{array}{ll}
\phi_{t}=-\delta_{1} \phi+\lambda \phi ,\; &\,x\in(r(m_{0}T), s(m_{0}T)), t\in(0^{+}, \tau],  \\[2mm]
\psi_{t}=d_{2}\Delta \psi+\frac{d_{2}\epsilon}{d_{1}}\psi_{x}-\delta_{2} \psi+\lambda \psi ,\; &\,x\in(r(m_{0}T), s(m_{0}T)), t\in(0^{+}, \tau],  \\[2mm]
\phi_{t}=d_{1}\Delta \phi+\epsilon\phi_{x}-a_{11}\phi+a_{12}\psi+\lambda \phi ,\; &\,x\in(r(m_{0}T), s(m_{0}T)), t\in(\tau, T],  \\[2mm]
\psi_{t}=d_{2}\Delta \psi+\frac{d_{2}\epsilon}{d_{1}}\psi_{x}-a_{22}\psi+f'(0) \phi+\lambda \psi  ,\; &\,x\in(r(m_{0}T), s(m_{0}T)), t\in(\tau, T],  \\[2mm]
\phi(x,t)=\psi(x,t)=0,\; &\, x\in\{r(m_{0}T), s(m_{0}T)\},  t\in(0, T],\\[2mm]
\phi(x,0^{+})=H'(0)\phi(x,0),\; &\, x\in(r(m_{0}T), s(m_{0}T)),  \\[2mm]
\psi(x,0^{+})=\psi(x,0),\; &\, x\in(r(m_{0}T), s(m_{0}T)),  \\[2mm]
\phi(0,x)=\phi(T,x),\psi(0,x)=\psi(T,x),\; &\,x\in[r(m_{0}T), s(m_{0}T)].
\end{array} \right.
\end{eqnarray}
Since $\lambda_{1}(r(m_{0}T), s(m_{0}T))<0$, there exists some sufficiently small constant $\epsilon>0$ such that the eigenvalue problem \eqref{4-17} has one real eigenvalue $\lambda^{\epsilon}_{1}$ with $\lambda^{\epsilon}_{1}<0$, and the corresponding eigenfunction pair $(\phi^{\epsilon}_{1}, \psi^{\epsilon}_{1})$ can be represented as
\begin{equation*}
\begin{aligned}
&\phi^{\epsilon}_{1}=\Phi(t)e^{-\frac{\epsilon}{2d_{1}}x}\sin\Big[\frac{\pi}{s(m_{0}T)-r(m_{0}T)}\Big(x-r(m_{0}T)\Big)\Big],  \\
&\psi^{\epsilon}_{1}=\Psi(t)e^{-\frac{\epsilon}{2d_{1}}x}\sin\Big[\frac{\pi}{s(m_{0}T)-r(m_{0}T)}\Big(x-r(m_{0}T)\Big)\Big],
\end{aligned}
\end{equation*}
where $(\Phi(t),\Psi(t))$ satisfies problem \eqref{3-15} with $\kappa^{l}_{1}$ replaced by $(\frac{\pi}{s(m_{0}T)-r(m_{0}T)})^{2}+(\frac{\epsilon}{2d_{1}})^{2}$.
A direct calculation derives that
\begin{equation}\label{4-20}
\begin{aligned}
&\phi^{\epsilon}_{1x}=-\frac{\epsilon}{2d_{1}}\phi^{\epsilon}_{1}+\frac{\pi}{s(m_{0}T)-r(m_{0}T)}e^{-\frac{\epsilon}{2d_{1}}x}\cos\Big[\frac{\pi}{s(m_{0}T)
-r(m_{0}T)}\Big(x-r(m_{0}T)\Big)\Big]\Phi(t),  \\
&\psi^{\epsilon}_{1x}=-\frac{\epsilon}{2d_{1}}\psi^{\epsilon}_{1}+\frac{\pi}{s(m_{0}T)-r(m_{0}T)}e^{-\frac{\epsilon}{2d_{1}}x}\cos\Big[\frac{\pi}{s(m_{0}T)
-r(m_{0}T)}\Big(x-r(m_{0}T)\Big)\Big]\Psi(t),
\end{aligned}
\end{equation}
which yield that there exists constant $x_{0}$ such that $\phi^{\epsilon}_{1x}(x_{0},t)=\psi^{\epsilon}_{1x}(x_{0},t)=0$, $\phi^{\epsilon}_{1x}(x,t)$, $\psi^{\epsilon}_{1x}(x,t)>0$ for $x\in[r(m_{0}T), x_{0})$, and  $\phi^{\epsilon}_{1x}(x,t)$, $\psi^{\epsilon}_{1x}(x,t)<0$ for $x\in( x_{0}, s(m_{0}T)]$. Specifically, we have that $x_{1}
=:\frac{s(m_{0}T)+r(m_{0}T)}{2}\in( x_{0}, s(m_{0}T)]$ for the reason that $\phi^{\epsilon}_{1x}(x_{1},t)$, $\psi^{\epsilon}_{1x}(x_{1},t)<0$. And, by calculation, we have that
\begin{equation}\label{4-19}
\begin{aligned}
&\phi^{\epsilon}_{1xx}=-\frac{\epsilon}{d_{1}}\phi^{\epsilon}_{1x}-\Big[\frac{\epsilon^{2}}{4d^{2}_{1}}+\Big(\frac{\pi}{s(m_{0}T)-r(m_{0}T)}\Big)^{2}\Big]\phi^{\epsilon}_{1},  \\
&\psi^{\epsilon}_{1xx}=-\frac{\epsilon}{d_{1}}\psi^{\epsilon}_{1x}-\Big[\frac{\epsilon^{2}}{4d^{2}_{1}}+\Big(\frac{\pi}{s(m_{0}T)-r(m_{0}T)}\Big)^{2}\Big]\psi^{\epsilon}_{1},
\end{aligned}
\end{equation}
which will be used in the sequel.

Subsequently, we aim to construct a lower solution $(\underline{u}, \underline{v})$ of model \eqref{Fixed-1}-\eqref{Fixed-3} for $x\in[r(m_{0}T), s(m_{0}T)]$ and $t\geq m_{0}T$. Let $(\lambda_{1}, (\phi_{1}, \psi_{1}))$ be the principal eigenpair of problem \eqref{4-17} with $\epsilon=0$. For such $\lambda_{1}$, there exists a sufficiently
small constant $\alpha_{0}>0$ such that $\lambda_{1}+\alpha_{0}<0$. It follows from \eqref{3-46} that there exists a large positive constant $A$ such that
\begin{equation}\label{4-18}
\frac{1}{A}\leq \frac{\phi_{1}}{\psi_{1}}\leq A
\end{equation}
for $x\in[r(m_{0}T), s(m_{0}T)]$ and $t\geq m_{0}T$.
Next, we construct the following functions
\begin{eqnarray*}
\underline{u}(x,t)=
\left\{
\begin{array}{ll}
\epsilon_{0} \phi_{1}(x,t),\; &\, (x,t)\in [r(m_{0}T), s(m_{0}T)]\times\{mT\}, \\[2mm]
\epsilon_{0} e^{(\lambda_{1}+\alpha_{0})T} \phi_{1}(x,t),\; &\, (x,t)\in [r(m_{0}T), s(m_{0}T)]\times\{(mT)^{+}\}, \\[2mm]
\epsilon_{0} e^{(\lambda_{1}+\alpha_{0})((m+1)T-t)}\phi_{1}(x,t),\; &\, (x,t)\in [r(m_{0}T), s(m_{0}T)]\times((mT)^{+}, (m+1)T]
\end{array} \right.
\end{eqnarray*}
and
\begin{eqnarray*}
\underline{v}(x,t)=
\left\{
\begin{array}{ll}
\epsilon_{0} \psi_{1}(x,t),\; &\, (x,t)\in [r(m_{0}T), s(m_{0}T)]\times\{mT\}, \\[2mm]
\epsilon_{0} e^{(\lambda_{1}+\alpha_{0})T} \psi_{1}(x,t),\; &\, (x,t)\in [r(m_{0}T), s(m_{0}T)]\times\{(mT)^{+}\}, \\[2mm]
\epsilon_{0} e^{(\lambda_{1}+\alpha_{0})((m+1)T-t)}\psi_{1}(x,t),\; &\, (x,t)\in [r(m_{0}T), s(m_{0}T)]\times((mT)^{+}, (m+1)T],
\end{array} \right.
\end{eqnarray*}
where $m=m_{0}, m_{0}+1, \cdots$, and positive constant $\epsilon_{0} $ will be determined later.

Before applying the usual comparison principle, we first check some conditions. Clearly,
\begin{equation*}
\underline{u}(r(m_{0}T),t)=\underline{u}(s(m_{0}T),t)=\underline{v}(r(m_{0}T),t)=\underline{v}(s(m_{0}T),t)=0
\end{equation*}
for $t\geq m_{0}T$. Standard calculations yield that
\begin{equation*}
\begin{split}
\underline{u}_{t}+\delta_{1} \underline{u}=&-\alpha_{0} \epsilon_{0} e^{(\lambda_{1}+\alpha_{0})((m+1)T-t)}  \phi_{1}(x,t) \leq 0
\end{split}
\end{equation*}
and
\begin{equation*}
\begin{split}
\overline{v}_{t}-d_{2}\Delta \overline{v}+\delta_{2} \overline{v}=&-\alpha_{0} \epsilon_{0} e^{(\lambda_{1}+\alpha_{0})((m+1)T-t)}  \psi_{1}(x,t) \leq 0
\end{split}
\end{equation*}
for $x\in(r(m_{0}T), s(m_{0}T))$ and $t\in((mT)^{+}, mT+\tau]$, and when $\epsilon_{0}$ is sufficiently small
\begin{equation*}
\begin{split}
\underline{u}_{t}-d_{1}\Delta \underline{u}+a_{11}\underline{u}-a_{12}\underline{v}=-\alpha_{0} \epsilon_{0} e^{(\lambda_{1}+\alpha_{0})((m+1)T-t)}  \phi_{1}(x,t) \leq 0
\end{split}
\end{equation*}
and
\begin{equation*}
\begin{split}
&\underline{v}_{t}-d_{2}\Delta \underline{v}+a_{22}\underline{v}-f(\underline{u})\\
=&  \epsilon_{0} e^{(\lambda_{1}+\alpha_{0})((m+1)T-t)}  \phi_{1}(x,t)(f'(0)-\alpha_{0}A)-f\big(\epsilon_{0} e^{(\lambda_{1}+\alpha_{0})((m+1)T-t)}  \phi_{1}(x,t)\big)\\
\leq & 0
\end{split}
\end{equation*}
for $x\in(r(m_{0}T), s(m_{0}T))$ and $t\in( mT+\tau, (m+1)T]$ by using \eqref{4-18} and assumption \textbf{(F)}. Based on assumption \textbf{(H)}, we have that for $x\in(r(m_{0}T), s(m_{0}T))$,
\begin{eqnarray*}
\underline{u}(x,(mT)^{+})=\epsilon_{0} H'(0) e^{(\lambda_{1}+\alpha_{0})T} \phi_{1}(x,mT)\leq H(\underline{u}(x,mT))
\end{eqnarray*}
and
\begin{eqnarray*}
\underline{v}(x,(mT)^{+})=\epsilon_{0} e^{(\lambda_{1}+\alpha_{0})T} \psi_{1}(x,mT)\leq\epsilon_{0} \psi_{1}(x,mT)\leq \underline{v}(x,mT)
\end{eqnarray*}
hold for $\epsilon_{0}$ is sufficiently small. For the initial functions, we can choose a sufficiently small $\epsilon_{0}$ such that
\begin{equation*}
\begin{aligned}
\underline{u}(x, m_{0}T)\leq u(x, m_{0}T)\text{~and~}\underline{v}(x, m_{0}T)\leq v(x, m_{0}T)
\end{aligned}
\end{equation*}
for $x\in[r(m_{0}T), s(m_{0}T)]$. Now, it follows from the usual comparison principle that
\begin{equation}\label{4-22}
\begin{aligned}
(\underline{u}, \underline{v})\leq (u, v)\text{~for~} (x,t)\in[r(m_{0}T), s(m_{0}T)]\times [m_{0}T, \infty).
\end{aligned}
\end{equation}

By using assumption $s_{\infty}<+\infty$, there exists a sufficiently large integer $m_{1}\geq m_{0}$ such that
\begin{eqnarray}\label{4-21}
s'(t)<\min\bigg\{\frac{d_{2}\epsilon[s(m_{0}T)-x_{1}]}{d_{1}(s_{\infty}-x_{1})}, \frac{\epsilon[s(m_{0}T)-x_{1}]}{s_{\infty}-x_{1}}\bigg\}
\end{eqnarray}
for $t\geq m_{1}T$, where $\epsilon$ and $x_{1}$ are defined above. Fix an arbitrary integer $M\geq  m_{1}+1$. We construct the following functions
\begin{eqnarray*}
\check{u}(x,t)=
\left\{
\begin{array}{ll}
\epsilon_{1} \phi^{\epsilon}_{1}(\bar{x},t),\; &\, (x,t)\in [x_{1}, s(t)]\times\{\bar{m}T\}, \\[2mm]
\epsilon_{1} e^{(\lambda^{\epsilon}_{1}+\alpha_{1})T} \phi^{\epsilon}_{1}(\bar{x},t),\; &\, (x,t)\in [x_{1}, s(t)]\times\{(\bar{m}T)^{+}\}, \\[2mm]
\epsilon_{1} e^{(\lambda^{\epsilon}_{1}+\alpha_{1})((\bar{m}+1)T-t)}\phi^{\epsilon}_{1}(\bar{x},t),\; &\, (x,t)\in [x_{1}, s(t)]\times((\bar{m}T)^{+}, (\bar{m}+1)T]
\end{array} \right.
\end{eqnarray*}
and
\begin{eqnarray*}
\check{v}(x,t)=
\left\{
\begin{array}{ll}
\epsilon_{1} \psi^{\epsilon}_{1}(\bar{x},t),\; &\, (x,t)\in [x_{1}, s(t)]\times\{\bar{m}T\}, \\[2mm]
\epsilon_{1} e^{(\lambda^{\epsilon}_{1}+\alpha_{1})T} \psi^{\epsilon}_{1}(\bar{x},t),\; &\, (x,t)\in [x_{1}, s(t)]\times\{(\bar{m}T)^{+}\}, \\[2mm]
\epsilon_{1} e^{(\lambda^{\epsilon}_{1}+\alpha_{1})((\bar{m}+1)T-t)}\psi^{\epsilon}_{1}(\bar{x},t),\; &\, (x,t)\in [x_{1}, s(t)]\times((\bar{m}T)^{+}, (\bar{m}+1)T],
\end{array} \right.
\end{eqnarray*}
where  $\bar{m}=m_{1}, m_{1}+1, \cdots, M-1$, $\alpha_{1}$ is a sufficiently small positive constant such that $\lambda^{\epsilon}_{1}+\alpha_{1}<0$,
\begin{eqnarray*}
\bar{x}=
\left\{
\begin{array}{ll}
x_{1}+\frac{s(m_{0}T)-x_{1}}{s(\bar{m}T)-x_{1}}(x-x_{1}),\; &\, (x,t)\in [x_{1}, s(\bar{m}T)]\times[\bar{m}T, \bar{m}T+\tau], \\[2mm]
x_{1}+\frac{s(m_{0}T)-x_{1}}{s(t)-x_{1}}(x-x_{1}),\; &\, (x,t)\in [x_{1}, s(t)]\times( \bar{m}T+\tau, (\bar{m}+1)T),
\end{array} \right.
\end{eqnarray*}
and the sufficiently small positive constant $\epsilon_{1} (\leq \epsilon_{0})$ will be determined later. Next, we show that pair $(\check{u}(x,t),\check{v}(x,t))$ is a lower solution of model \eqref{Fixed-1}-\eqref{Fixed-3} for $x\in[x_{1}, s(t)]$ and $t\in [m_{1}T, MT]$. Before applying \autoref{remark 2.1}, we first check some conditions. It is easily seen that
\begin{equation*}
\check{u}(s(t),t)=\check{v}(s(t),t)=0\text{~for~}t\in (m_{1}T, MT].
\end{equation*}
In virtue of \eqref{4-22}, one can obtain that
\begin{equation*}
\begin{aligned}
(\underline{u}(x_{1}, t), \underline{v}(x_{1}, t))\leq (u(x_{1}, t), v(x_{1}, t))\text{~for~} t\in (m_{1}T, MT],
\end{aligned}
\end{equation*}
which deduces that for $t\in(m_{1}T, MT]$,
\begin{equation*}
\begin{aligned}
(\check{u}(x_{1}, t), \check{v}(x_{1}, t))\leq(\underline{u}(x_{1}, t), \underline{v}(x_{1}, t))\leq (u(x_{1}, t), v(x_{1}, t))
\end{aligned}
\end{equation*}
holds for $\epsilon_{1}$ sufficiently small. A direct computation yields yield that
\begin{equation*}
\begin{split}
\check{u}_{t}+\delta_{1} \check{u}=& -\epsilon_{1} \alpha_{1} e^{(\lambda^{\epsilon}_{1}+\alpha_{1})((\bar{m}+1)T-t)}\phi^{\epsilon}_{1}(\bar{x},t)  \leq 0
\end{split}
\end{equation*}
and
\begin{equation*}
\begin{split}
&\check{v}_{t}-d_{2}\Delta \check{v}+\delta_{2} \check{v}\\
=&\epsilon_{1}  e^{(\lambda^{\epsilon}_{1}+\alpha_{1})((\bar{m}+1)T-t)}\bigg[\frac{d_{2}\epsilon}{d_{1}}\psi^{\epsilon}_{1\bar{x}}-\alpha_{1}\psi^{\epsilon}_{1}+d_{2}\bigg(1-\Big(\frac{s(m_{0}T)
-x_{1}}{s(\bar{m}T)-x_{1}}\Big)^{2}\bigg) \psi^{\epsilon}_{1\bar{x}\bar{x}}\bigg]\\
=&\Bigg\{-\bigg[\alpha_{1}+d_{2}\bigg(1-\Big(\frac{s(m_{0}T)-x_{1}}{s(\bar{m}T)-x_{1}}\Big)^{2}\bigg)\bigg(\frac{\epsilon^{2}}{4d^{2}_{1}}+\Big(\frac{\pi}{s(m_{0}T)
-r(m_{0}T)}\Big)^{2}\bigg)\bigg]\psi^{\epsilon}_{1}\\
&+\frac{d_{2}\epsilon}{d_{1}}\Big(\frac{s(m_{0}T)-x_{1}}{s(\bar{m}T)-x_{1}}\Big)^{2}\psi^{\epsilon}_{1\bar{x}}\Bigg\}\epsilon_{1}  e^{(\lambda^{\epsilon}_{1}
+\alpha_{1})((\bar{m}+1)T-t)}\\
\leq &0
\end{split}
\end{equation*}
for $x\in(x_{1}, s(\bar{m}T))$ and $t\in((\bar{m}T)^{+}, \bar{m}T+\tau]$ by using \eqref{4-20} and \eqref{4-19}, and when $\epsilon_{1}$ is sufficiently small
\begin{equation*}
\begin{split}
&\check{u}_{t}-d_{1}\Delta \check{u}+a_{11}\check{u}-a_{12}\check{v}\\
=&\epsilon_{1}  e^{(\lambda^{\epsilon}_{1}+\alpha_{1})((\bar{m}+1)T-t)}\bigg[-\alpha_{1}\phi^{\epsilon}_{1}+\bigg(\epsilon-\frac{(s(m_{0}T)-x_{1})(x-x_{1})s'(t)}
{(s(t)-x_{1})^{2}}\bigg)\phi^{\epsilon}_{1\bar{x}}\\
&+d_{1}\bigg(1-\Big(\frac{s(m_{0}T)-x_{1}}{s(t)-x_{1}}\Big)^{2}\bigg) \phi^{\epsilon}_{1\bar{x}\bar{x}}\bigg]\\
=&\Bigg\{-\bigg[\alpha_{1}+d_{1}\bigg(1-\Big(\frac{s(m_{0}T)-x_{1}}{s(t)-x_{1}}\Big)^{2}\bigg)\bigg(\frac{\epsilon^{2}}{4d^{2}_{1}}
+\Big(\frac{\pi}{s(m_{0}T)-r(m_{0}T)}\Big)^{2}\bigg)\bigg]\phi^{\epsilon}_{1}\\
&+\bigg[\epsilon\Big(\frac{s(m_{0}T)-x_{1}}{s(t)-x_{1}}\Big)^{2}-\frac{(s(m_{0}T)-x_{1})(x-x_{1})s'(t)}{(s(t)-x_{1})^{2}}\bigg]\phi^{\epsilon}_{1\bar{x}}\Bigg\}\epsilon_{1}  e^{(\lambda^{\epsilon}_{1}+\alpha_{1})((\bar{m}+1)T-t)}\\
\leq& 0
\end{split}
\end{equation*}
and
\begin{equation*}
\begin{split}
&\check{v}_{t}-d_{2}\Delta \check{v}+a_{22}\check{v}-f(\check{u})\\
=& \epsilon_{1}  e^{(\lambda^{\epsilon}_{1}+\alpha_{1})((\bar{m}+1)T-t)}\bigg[-\alpha_{1}\psi^{\epsilon}_{1}+f'(0)\phi^{\epsilon}_{1}+\bigg(\frac{d_{2}\epsilon}{d_{1}}
-\frac{(s(m_{0}T)-x_{1})(x-x_{1})s'(t)}{(s(t)-x_{1})^{2}}\bigg)\psi^{\epsilon}_{1\bar{x}}\\
&+d_{2}\bigg(1-\Big(\frac{s(m_{0}T)-x_{1}}{s(t)-x_{1}}\Big)^{2}\bigg) \psi^{\epsilon}_{1\bar{x}\bar{x}} \bigg]-f\bigg(\epsilon_{1} e^{(\lambda^{\epsilon}_{1}+\alpha_{1})((\bar{m}+1)T-t)}\phi^{\epsilon}_{1}\bigg)\\
\leq&\Bigg\{f'(0)\phi^{\epsilon}_{1}-\bigg[\alpha_{1}+d_{2}\bigg(1-\Big(\frac{s(m_{0}T)-x_{1}}{s(t)-x_{1}}\Big)^{2}\bigg)\bigg(\frac{\epsilon^{2}}{4d^{2}_{1}}
+\Big(\frac{\pi}{s(m_{0}T)-r(m_{0}T)}\Big)^{2}\bigg)\bigg]\psi^{\epsilon}_{1}\\
&+\bigg[\frac{d_{2}\epsilon}{d_{1}}\Big(\frac{s(m_{0}T)-x_{1}}{s(t)-x_{1}}\Big)^{2}-\frac{(s(m_{0}T)-x_{1})(x-x_{1})s'(t)}{(s(t)-x_{1})^{2}}\bigg]
\phi^{\epsilon}_{1\bar{x}}\Bigg\}\epsilon_{1}  e^{(\lambda^{\epsilon}_{1}+\alpha_{1})((\bar{m}+1)T-t)}\\
&-f\bigg(\epsilon_{1} e^{(\lambda^{\epsilon}_{1}+\alpha_{1})((\bar{m}+1)T-t)}\phi^{\epsilon}_{1}\bigg)\\
\leq &\epsilon_{1}  e^{(\lambda^{\epsilon}_{1}+\alpha_{1})((\bar{m}+1)T-t)}\phi^{\epsilon}_{1} (f'(0)-\alpha_{1}A^{-1})-f\bigg(\epsilon_{1} e^{(\lambda^{\epsilon}_{1}+\alpha_{1})((\bar{m}+1)T-t)}\phi^{\epsilon}_{1}\bigg)\\
\leq & 0
\end{split}
\end{equation*}
for $x\in(x_{1}, s(t))$ and $t\in( \bar{m}T+\tau, (\bar{m}+1)T]$ by using \eqref{4-20}, \eqref{4-19}, \eqref{4-18}, \eqref{4-21}, and assumption \textbf{(F)}. It follows from assumption \textbf{(H)} that for $x\in(x_{1}, s(t))$,
\begin{eqnarray*}
\check{u}(x,(\bar{m}T)^{+})=\epsilon_{1} H'(0) e^{(\lambda^{\epsilon}_{1}+\alpha_{1})T} \phi^{\epsilon}_{1}(\bar{x},mT)\leq H(\check{u}(x,\bar{m}T))
\end{eqnarray*}
and
\begin{eqnarray*}
\check{v}(x,(\bar{m}T)^{+})=\epsilon_{1}e^{(\lambda^{\epsilon}_{1}+\alpha_{1})T}\psi^{\epsilon}_{1}(\bar{x},\bar{m}T)\leq\epsilon_{1} \psi^{\epsilon}_{1}(\bar{x},\bar{m}T)= \check{v}(x,\bar{m}T)
\end{eqnarray*}
hold for $\epsilon_{1}$ sufficiently small. For the initial functions, we can choose a sufficiently small $\epsilon_{1}$ such that
\begin{equation*}
\begin{aligned}
\check{u}(x, m_{1}T)\leq u(x, m_{1}T)\text{~and~}\check{v}(x, m_{1}T)\leq v(x, m_{1}T)
\end{aligned}
\end{equation*}
for $x\in[x_{1}, s(m_{1}T)]$. Now, we apply \autoref{remark 2.1} to obtain that
\begin{equation*}
\begin{aligned}
\check{u}(x, t)\leq u(x, t)\text{~and~}\check{v}(x, t)\leq v(x, t)
\end{aligned}
\end{equation*}
for $x\in[x_{1}, s(t)]$ and $t\in [m_{1}T, MT]$. Note that $M$ is arbitrary and the above analysis is independent of $M$. Therefore, we have that
\begin{equation}\label{4-23}
\begin{aligned}
\check{u}(x, t)\leq u(x, t)\text{~and~}\check{v}(x, t)\leq v(x, t)
\end{aligned}
\end{equation}
holds for $x\in[x_{1}, s(t)]$ and $t\in [m_{1}T, +\infty)$. On the one hand, it follows from \eqref{4-23} that
\begin{equation}\label{4-24}
\begin{aligned}
\lim\limits_{t\rightarrow +\infty}u_{x}(s(t), t)\leq  \lim\limits_{t\rightarrow +\infty} \check{u}_{x}(s(t), t)<0\text{~and~}\lim\limits_{t\rightarrow +\infty}v_{x}(s(t), t)\leq  \lim\limits_{t\rightarrow +\infty} \check{v}_{x}(s(t), t)<0.
\end{aligned}
\end{equation}
On the other hand, we have from model \eqref{Fixed-1}-\eqref{Fixed-3} that
\begin{equation*}
\begin{aligned}
\lim\limits_{t\rightarrow +\infty}u_{x}(s(t), t)=-\frac{\mu_{2}}{\mu_{1}}\lim\limits_{t\rightarrow +\infty}v_{x}(s(t), t),
\end{aligned}
\end{equation*}
which creates a contradiction with \eqref{4-24}. Therefore, the desired assertion holds. This proof is completed.
\end{proof}
\end{lemma}
\autoref{lemma 4-2} above demonstrates that when $\lambda_{1}(\infty)<0$, the free boundaries $r(t)$ and $s(t)$ are either both finite or both infinite. \autoref{lemma 4-1}
shows the case of both being finite. The following lemma will address the scenario where both are infinite.
\begin{lemma}\label{lemma 4-3}
If $\lambda_{1}(\infty)<0$ and $s_{\infty}-r_{\infty}=+\infty$, then
\begin{equation*}
\lim\limits_{m\rightarrow+\infty}u(x,t+mT)=W(t)~\text{and}~\lim\limits_{m\rightarrow+\infty}v(x,t+mT)=Z(t)
\end{equation*}
locally uniformly in $\mathcal{R}$ and uniformly in $[0,T]$, where $(W(t),Z(t))$ denotes the solution of \eqref{3-42}.
\begin{proof}
In virtue of \autoref{lemma 4-2}, one can obtain that $r_{\infty}=-\infty$ and $s_{\infty}=+\infty$. We prove the assertion by showing
\begin{equation}\label{4-25}
\limsup\limits_{m\rightarrow +\infty}(u(x,t+mT), v(x,t+mT))\preceq(W(t),Z(t))
\end{equation}
uniformly in $\mathds{R}\times [0,T]$ and
\begin{equation}\label{4-26}
\liminf\limits_{m\rightarrow +\infty}(u(x,t+mT), v(x,t+mT))\succeq(W(t),Z(t))
\end{equation}
locally uniformly in $\mathds{R}\times[0,T]$. We first prove that \eqref{4-25} holds. By using the usual comparison principle, we have that
\begin{equation}\label{4-27}
(u(x,t), v(x,t))\preceq(U(t),V(t))
\end{equation}
for $x\in [r(t), s(t)]$ and $t\in\mathds{R}^{+}$, where $(U,V)$ is defined in the proof of \autoref{theorem 4-2}. Then, it follows from \eqref{4-27}, and \textcolor{blue}{Lemmas} \ref{lemma 3.1.7} and \ref{lemma 3.1.5} that
\begin{equation}\label{4-28}
\limsup\limits_{m\rightarrow +\infty}(u(x,t+mT), v(x,t+mT)\preceq\lim\limits_{m\rightarrow +\infty}(U(t+mT),V(t+mT))=(W(t),Z(t))
\end{equation}
uniformly in $\mathds{R}\times [0,T]$. Hence, \eqref{4-25} holds.

It remains to prove that \eqref{4-26} also holds. In virtue of $\lambda_{1}(\infty)<0$, (2) of \autoref{lemma 3.1.3} yields that there exists a sufficiently large positive
constant $l$ such that $\lambda_{1}(-l, l)<0$. Additionally, it follows from $r_{\infty}=-\infty$ and $s_{\infty}=+\infty$ that a sufficiently large positive integer $m^{*}$
can be found such that $s(t)\geq l$ and $r(t)\leq -l$ for $t\geq m^{*}T$. Let $(u_{l}(x,t),v_{l}(x,t))$ be the unique solution of the initial-boundary value
problem \eqref{4-29} with $m_{0}$, $r_{\infty}+\epsilon$, and $s_{\infty}-\epsilon$ replaced by $m^{*}$, $-l$, and $l$, respectively. Then, it follows from the usual
comparison principle that
\begin{equation}\label{4-30}
u_{l}(x,t)\leq u(x,t)\text{~and~} v_{l}(x,t)\leq v(x,t)
\end{equation}
for $x\in[-l, l]$ and $t\geq m^{*}T$. Since $\lambda_{1}(-l, l)<0$, we have from (1) of \autoref{theorem 3.2} that
\begin{equation}\label{4-31}
\lim\limits_{m\rightarrow+\infty}\big(u_{l}(x,t+mT),v_{l}(x,t+mT)\big)=\big(w_{l}(x,t), z_{l}(x,t)\big)
\end{equation}
uniformly for $(x,t)\in[-l, l] \times[0, +\infty)$, where $\big(w_{l}(x,t), z_{l}(x,t)\big)$ denotes the unique strongly positive periodic solution of  problem \eqref{3-32} with $l_{1}$ and $l_{2}$ replaced by $-l$ and $l$, respectively. It follows from \eqref{4-30} and \eqref{4-31} that
\begin{equation}\label{4-32}
\liminf\limits_{m\rightarrow+\infty}\big(u(x,t+mT),v(x,t+mT)\big)\geq \big(w_{l}(x,t), z_{l}(x,t)\big)
\end{equation}
uniformly for $(x,t)\in[-l, l] \times[0, +\infty)$. Finally, by taking $l\rightarrow +\infty$ in \eqref{4-32}, \autoref{lemma 3.1.6} yields that
\begin{equation*}
\liminf\limits_{m\rightarrow+\infty}\big(u(x,t+mT),v(x,t+mT)\big)\geq \big(W(t), Z(t)\big)
\end{equation*}
locally uniformly in $\mathcal{R}$ and uniformly in $[0,T]$, which implies that \eqref{4-26} also holds. The ends this proof.
\end{proof}
\end{lemma}
Combing \textcolor{blue}{Theorems} \ref{theorem 4-1} and \ref{theorem 4-2}, and \textcolor{blue}{Lemmas} \ref{lemma 4-1}-\ref{lemma 4-3}, we immediately obtain the following spreading-vanishing dichotomy.
\begin{theorem}\label{theorem 4-3}
Assume that (\textbf{G}) holds. The unique solution of model \eqref{Fixed-1}-\eqref{Fixed-3} is either vanishing or spreading.
\end{theorem}
\subsection{Spreading-vanishing criteria}\label{Section-4-2}
It should be acknowledged that the judgment conditions in \textcolor{blue}{Lemmas} \ref{lemma 4-1} and \ref{lemma 4-3} are not easily verifiable, since information about the free boundaries $s(t)$ and $r(t)$ cannot be known in advance. To address this, this subsection will provide readily verifiable alternative judgment condition by dividing the
discussion into two cases: $\lambda_{1}(-s_{0},s_{0}):=\lambda_{1}(s_{0})\leq 0$ and $\lambda_{1}(s_{0})>0$.
\begin{theorem}\label{theorem 4-4}
If $\lambda_{1}(\infty)<0$ and $\lambda_{1}(s_{0})\leq 0$, then the solution of model \eqref{Fixed-1}-\eqref{Fixed-3} is spreading.
\begin{proof}
Suppose by contradiction that $s_{\infty}-r_{\infty}<\infty$. By using \autoref{lemma 4-1}, we have that $\lambda_{1}(r_{\infty}, s_{\infty})\geq 0$, which combined with (2) of \autoref{lemma 3.1.3} deduces that
\begin{equation*}
\lambda_{1}(s_{0})>\lambda_{1}(r_{\infty}, s_{\infty})\geq 0.
\end{equation*}
This shows a contradiction with condition $\lambda_{1}(s_{0})\leq 0$. Therefore, we have that $s_{\infty}-r_{\infty}=\infty$, which combined with \autoref{lemma 4-3}
deduces that the assertion holds. This ends the proof.
\end{proof}
\end{theorem}

For case $\lambda_{1}(s_{0})>0$, let us discuss the effect of coefficients $\mu_{1}$ and $\mu_{2}$ on spreading and vanishing. As a first step, when $\mu_{1}+\mu_{2}$
is sufficiently large, we have the following assertion.
\begin{lemma}\label{lemma 4-4}
If $\lambda_{1}(\infty)<0$ and $\lambda_{1}(s_{0})>0$, then there exists a positive constant $\overline{\mu}$ such that spreading happens as $\mu_{1}+\mu_{2}\geq \overline{\mu}$.
\begin{proof}
To begin with, consider the auxiliary free boundary problem
\begin{eqnarray}\label{4-33}
\left\{
\begin{array}{ll}
\zeta_{t}=-\delta_{1} \zeta ,\; &\,x\in(g(t),h(t)), t\in((mT)^{+}, mT+\tau],  \\[2mm]
\eta_{t}= d_{2}\Delta \eta-\delta_{2} \eta ,\; &\,x\in(g(t),h(t)), t\in((mT)^{+}, mT+\tau],  \\[2mm]
\zeta_{t}= d_{1}\Delta \zeta-a_{11}\zeta ,\; &\,x\in(g(t),h(t)), t\in(mT+\tau, (m+1)T],  \\[2mm]
\eta_{t}= d_{2}\Delta \eta-a_{22}\eta ,\; &\,x\in(g(t),h(t)), t\in(mT+\tau, (m+1)T],  \\[2mm]
g(t)= g(mT),~h(t)=h(mT),\; &\,  t\in[mT, mT+\tau],  \\[2mm]
g'(t)= -\mu_{1}\zeta_{x}(g(t),t)-\mu_{2}\eta_{x}(g(t),t),\; &\,  t\in(mT+\tau, (m+1)T),  \\[2mm]
h'(t)=-\mu_{1}\zeta_{x}(h(t),t)-\mu_{2}\eta_{x}(h(t),t),\; &\,  t\in(mT+\tau, (m+1)T),  \\[2mm]
\zeta(x,t)= \eta(x,t)=0,\; &\, x\in\{g(t),h(t)\},  t\in(0,\infty),\\[2mm]
\zeta(x,(mT)^{+})= H(\zeta(x,mT)),\; &\, x\in(g(t),h(t)),  \\[2mm]
\eta(x,(mT)^{+})= \eta(x,mT),\; &\, x\in(g(t),h(t)),\\[2mm]
g(t)=-s_{0}~~\zeta(x,0)=u_{0}(x),\; &\, x\in[-s_{0},s_{0}],\\[2mm]
h(t)=s_{0}~~~~\eta(x,0)=v_{0}(x),\; &\, x\in[-s_{0},s_{0}].
\end{array} \right.
\end{eqnarray}
It follows from parabolic regularity theory and the contraction mapping principle, one can obtain that problem \eqref{4-33} admits a unique nonnegative classical
solution $(\zeta, \eta, g, h)$ which is well defined for all $t\geq 0$, and $\zeta_{x}(g(t),t), \eta_{x}(g(t),t), -\zeta_{x}(h(t),t), -\eta_{x}(h(t),t)>0$ for all $t\geq 0$.
To highlight the dependence of solution on $\mu=:(\mu_{1}, \mu_{2})$, we write $(u^{\mu}, v^{\mu}, r^{\mu}, s^{\mu})$ and $(\zeta^{\mu}, \eta^{\mu}, g^{\mu}, h^{\mu})$
instead of $(u, v, r, s)$ and $(\zeta, \eta, g, h)$, respectively. By using \autoref{remark 2.1}, we have that
\begin{equation}\label{4-39}
r^{\mu}(t) \leq g^{\mu}(t)\text{~and~} s^{\mu}(t) \geq h^{\mu}(t)
\end{equation}
for $t\geq 0$, and
\begin{equation*}
\zeta^{\mu}(x,t)\leq u^{\mu}(x,t)\text{~and~} \eta^{\mu}(x,t)\leq v^{\mu}(x,t)
\end{equation*}
for $x\in[g^{\mu}(t), h^{\mu}(t)]$ and $t\geq 0$.

We claim that there exists a positive constant $l^{*}$ such that $\lambda_{1}(r_{\infty}, s_{\infty})<0$ as $s_{\infty}-r_{\infty}>l^{*}$. Since $\lambda_{1}(\infty)<0$ and $\lambda_{1}(s_{0})>0$, it follows from (2) of \autoref{lemma 3.1.3} that there exists a positive constant $l^{*}$ such that $\lambda_{1}(-0.5l^{*}, 0.5l^{*})=0$. Using (2) of \autoref{lemma 3.1.3} again yields this claim holds.

In what follows, we will prove that there exists a positive constant $\overline{\mu}$ such that
\begin{equation}\label{4-40}
r^{\mu}(T)\leq -0.5l^{*} \text{~and~}s^{\mu}(T)\geq 0.5l^{*}
\end{equation}
for all $\mu_{1}+\mu_{2}\geq \overline{\mu}$. We consider the following free boundary problem
\begin{eqnarray}\label{4-34}
\left\{
\begin{array}{ll}
\underline{\zeta}_{t}=-\delta_{1} \underline{\zeta} ,\; &\,x\in(\underline{g}(t),\underline{h}(t)), t\in(0^{+}, \tau],  \\[2mm]
\underline{\eta}_{t}= d_{2}\Delta \underline{\eta}_{t}-\delta_{2} \underline{\eta}_{t} ,\; &\,x\in(\underline{g}(t),\underline{h}(t)), t\in(0^{+}, \tau],  \\[2mm]
\underline{\zeta}_{t}= d_{1}\Delta \underline{\zeta}-a_{11}\underline{\zeta} ,\; &\,x\in(\underline{g}(t),\underline{h}(t)), t\in(\tau, T],  \\[2mm]
\underline{\eta}_{t}= d_{2}\Delta \underline{\eta}-a_{22}\underline{\eta}  ,\; &\,x\in(\underline{g}(t),\underline{h}(t)), t\in(\tau, T],  \\[2mm]
\underline{\zeta}(x,t)=\underline{\eta}(x,t)=0,\; &\, x\in\{\underline{g}(t),\underline{h}(t)\},  t\in(0,T],\\[2mm]
\underline{\zeta}(x,0^{+})= H(\underline{\zeta}(x,0)),\; &\, x\in(\underline{g}(t),\underline{h}(t)),  \\[2mm]
\underline{\eta}(x,0^{+})= \underline{\eta}(x,0),\; &\, x\in(\underline{g}(t),\underline{h}(t)),
\end{array} \right.
\end{eqnarray}
where the chosen free boundaries $\underline{g}(t)$ and $\underline{h}(t)$ are smooth functions, which satisfy
\begin{eqnarray}\label{4-35}
\left\{
\begin{array}{ll}
\underline{g}(t)= -0.5s_{0},~~~~~~~~\underline{h}(t)=0.5s_{0},\; &\,  t\in[0, \tau],  \\[2mm]
\underline{g}(t)\text{~is ~decreasing~and~}\underline{g}(T)=-0.5l^{*},\; &\,  t\in(\tau, T),  \\[2mm]
\underline{h}(t)\text{~is ~increasing~and~}~\underline{h}(T)=0.5l^{*},\; &\,  t\in(\tau, T),
\end{array} \right.
\end{eqnarray}
and the initial functions $\underline{\zeta}_{0}(x)$ and $\underline{\eta}_{0}(x)$ satisfy
\begin{eqnarray}\label{4-36}
\left\{
\begin{array}{l}
\underline{\zeta}_{0}(\pm 0.5s_{0})=0,~~\underline{\zeta}'_{0}(-0.5s_{0})=-\underline{\zeta}'_{0}(0.5s_{0})>0, ~\text{and~}0<\underline{\zeta}_{0}(x)\leq u_{0}(x)~\text{in}~(-0.5s_{0},0.5s_{0}), \\[2mm]
\underline{\eta}_{0}(\pm 0.5s_{0})=0, ~\underline{\eta}'_{0}(-0.5s_{0})=-\underline{\eta}'_{0}(0.5s_{0})>0,~~\text{and~}0<\underline{\eta}_{0}(x)\leq v_{0}(x)~\text{in}~(-0.5s_{0},0.5s_{0}).
\end{array} \right.
\end{eqnarray}
It follows from parabolic regularity theory and the contraction mapping principle that problem \eqref{4-34}-\eqref{4-36} admits a unique nonnegative
classical solution $(\underline{\zeta}, \underline{\eta}, \underline{g}, \underline{h})$ which is well defined for all $t\in[0,T]$, and $\underline{\zeta}_{x}(\underline{g}(t),t), \underline{\eta}_{x}(\underline{g}(t),t), -\underline{\zeta}_{x}(\underline{h}(t),t), -\underline{\eta}_{x}(\underline{h}(t),t)>0$ for all $t\in[0,T]$. Therefore, there exists
a positive constant $\overline{\mu}$ such that when $\mu_{1}+\mu_{2}\geq \overline{\mu}$, we have that
\begin{eqnarray*}
\left\{
\begin{array}{ll}
\underline{g}'(t)\geq-\mu_{1}\underline{\zeta}_{x}(\underline{g}(t),t)-\mu_{2}\underline{\eta}_{x}(\underline{g}(t),t),\; &\,  t\in(\tau, T),  \\[2mm]
\underline{h}'(t)\leq-\mu_{1}\underline{\zeta}_{x}(\underline{h}(t),t)-\mu_{2}\underline{\eta}_{x}(\underline{h}(t),t),\; &\,  t\in(\tau, T).
\end{array} \right.
\end{eqnarray*}
For problem \eqref{4-33}, we can establish the comparison principle analogous to \autoref{lemma 2-1} and \autoref{remark 2.1}, which deduces when $\mu_{1}+\mu_{2}\geq \overline{\mu}$,  that
\begin{equation}\label{4-37}
\underline{g}(t)\geq g^{\mu}(t)\text{~and~}\underline{h}(t)\leq h^{\mu}(t)\text{~for~}t\in[0,T].
\end{equation}
Combing \eqref{4-39} and \eqref{4-37}, we have that
\begin{equation*}
\underline{g}(t)\geq r^{\mu}(t)\text{~and~}\underline{h}(t)\leq s^{\mu}(t)\text{~for~}t\in[0,T],
\end{equation*}
which combined with \eqref{4-35} deduces that \eqref{4-40} holds. Then, it follows from the claim that
\begin{equation}\label{4-41}
\lambda_{1}(r_{\infty}, s_{\infty})<0.
\end{equation}
Since $\lambda_{1}(\infty)<0$, one can obtain from \autoref{lemma 4-2} that $r_{\infty}$ and $s_{\infty}$ are either both finite or both infinite. If $r_{\infty}$ and $s_{\infty}$ are
both finite, then \autoref{lemma 4-1} yields that $\lambda_{1}(r_{\infty}, s_{\infty})\geq0$, which creates a contradiction with \eqref{4-41}. Therefore, $s_{\infty}=-r_{\infty}=+\infty$,
which combined with \autoref{lemma 4-3} gives that spreading happens. This proof is completed.
\end{proof}
\end{lemma}
Secondly, when $\mu_{1}+\mu_{2}$ is small enough, we have the following assertion.
\begin{lemma}\label{lemma 4-5}
If $\lambda_{1}(\infty)<0$ and $\lambda_{1}(s_{0})>0$, then there exists a positive constant $\underline{\mu}$ such that vanishing happens as $\mu_{1}+\mu_{2}\leq \underline{\mu}$.
\begin{proof}
To proof this lemma, we only need to make a minor modification in the proof of \autoref{theorem 4-1}. We construct the same quadruple $(\overline{u}(x,t), \overline{v}(x,t), -\overline{s}(t), \overline{s}(t))$. Then, we take $\hat{s}=s_{0}$, $\gamma=\frac{\lambda_{1}(s_{0})}{2}$, sufficiently large $C_{0}$, and $\eta$ such that
\begin{equation*}
1-(1+\frac{\eta}{2})^{-2}\leq \frac{\lambda_{1}(s_{0})}{2}\frac{4s^{2}_{0}}{(d_{1}+d_{2})\pi^{2}}.
\end{equation*}
We easily check that if
\begin{equation*}
\mu_{1}+\mu_{2}\leq \underline{\mu}=:\frac{\eta\gamma(2+\eta)s^{2}_{0}}{2C_{0}\pi \max\{n_{2}, n_{4}\}},
\end{equation*}
then quadruple $(\overline{u}(x,t), \overline{v}(x,t), -\overline{s}(t), \overline{s}(t))$ is an upper solution of model \eqref{Fixed-1}-\eqref{Fixed-3}, and hence vanishing happens.
This ends the proof.
\end{proof}
\end{lemma}
With the help of \textcolor{blue}{Lemmas} \ref{lemma 4-4} and \ref{lemma 4-5}, the same argument as in \cite[Theorem 3.9]{Du-Lin-SIAM-JMA} yields the following conclusion.
\begin{theorem}\label{theorem 4-5}
Assume that $\mu_{2}=\rho \mu_{1}$. If $\lambda_{1}(\infty)<0$ and $\lambda_{1}(s_{0})> 0$, then there exists a positive constant $\mu^{*}$ such that vanishing happens for $0<\mu_{1}\leq \mu^{*}$, and spreading happens for $\mu_{1}> \mu^{*}$.
\end{theorem}
Combing \textcolor{blue}{Theorems} \ref{theorem 4-1}, \ref{theorem 4-2}, \ref{theorem 4-4}, and \ref{theorem 4-5}, we immediately obtain the following criteria.
\begin{theorem}\label{theorem 4-6}
Assume that (\textbf{G}) holds and $\mu_{2}=\rho \mu_{1}$. We have the following assertions:
\begin{enumerate}
\item[$(1)$]
If $\lambda_{1}(\infty)\geq0$, then vanishing happens, see \textcolor{blue}{Theorems} \ref{theorem 4-1} and \ref{theorem 4-2};
\item[$(2)$]
If $\lambda_{1}(\infty)<0$ and $\lambda_{1}(s_{0})\leq0$, then spreading happens, see \autoref{theorem 4-4};
\item[$(3)$]
If $\lambda_{1}(\infty)<0$ and $\lambda_{1}(s_{0})>0$, then there exists a positive constant $\mu^{*}$ such that vanishing happens for $0<\mu_{1}\leq \mu^{*}$, and
spreading happens for $\mu_{1}> \mu^{*}$, see \autoref{theorem 4-5}.
\end{enumerate}
\end{theorem}
\section{Numerical simulation}\label{Section-5}
The theoretical analysis has been established in the preceding sections. This section employs three numerical examples to validate the theoretical findings and to visually illustrate the effects of impulsive intervention coefficient and duration of the dry season on the transmission dynamics of faecal-oral diseases.

The validation of our theoretical results requires prior knowledge of the sign of the principal eigenvalue. According to \autoref{lemma 3.1.3}\textcolor{blue}{(2)}, the principal eigenvalue exhibits strongly monotonic decrease with respect to the length of the infected interval. It is important to note that the boundaries of the infected area are not fixed: it remains constant during the dry season but expands continuously throughout the rainy season. Consequently, determining the sign via this monotonicity depends on the prior knowledge of the free boundaries. However, since these free boundaries cannot be known in advance, we adopt the simulated ones as approximations of the true boundaries for our numerical computations. Subsequently, the Laplace operator in the eigenvalue problem \eqref{3-3} is discretized using a second-order central difference scheme, which allows us to compute the value of the principal eigenvalue with the aid of mathematical software.

To conserve space, only the numerical simulations of infectious agents are presented here. In all simulations, we consider model \eqref{Fixed-1}-\eqref{Fixed-3} with Beverton-Holt function, that is
\begin{equation*}
f(u)=\frac{mu}{a+u}.
\end{equation*}
The initial region is chosen as $[-2, 2]$, and the initial functions are chosen as
\begin{equation*}
u_{0}(x)=0.4\cos\Big(\frac{\pi x}{4}\Big)\text{~and~}v_{0}(x)=0.1\cos\Big(\frac{\pi x}{4}\Big) \text {~for~}x\in[-2, 2].
\end{equation*}
The remaining parameters of model \eqref{Fixed-1}-\eqref{Fixed-3} are specified in the respective examples below to illustrate the varying dynamical behaviors.
\subsection{Monotonicity of the principal eigenvalue}
In this subsection, we employ numerical simulations to verify the monotonicity of the principal eigenvalue with respect to the length of the infected region, the coefficient of impulsive interventions, and the duration of the dry season. Some parameters are set to $d_{1}=5$, $d_{2}=40$, $\delta_{1}=0.6$, $\delta_{2}=0.9$, $a_{11}=0.2$, $a_{12}=0.8$, $a_{22}=0.3$, $a=1$, $m=1.5$, and $T=10$. Other of the parameters are set in each of the following simulations.
\begin{figure}[ht]
\centering
\subfigure{ {
\includegraphics[width=0.47\textwidth]{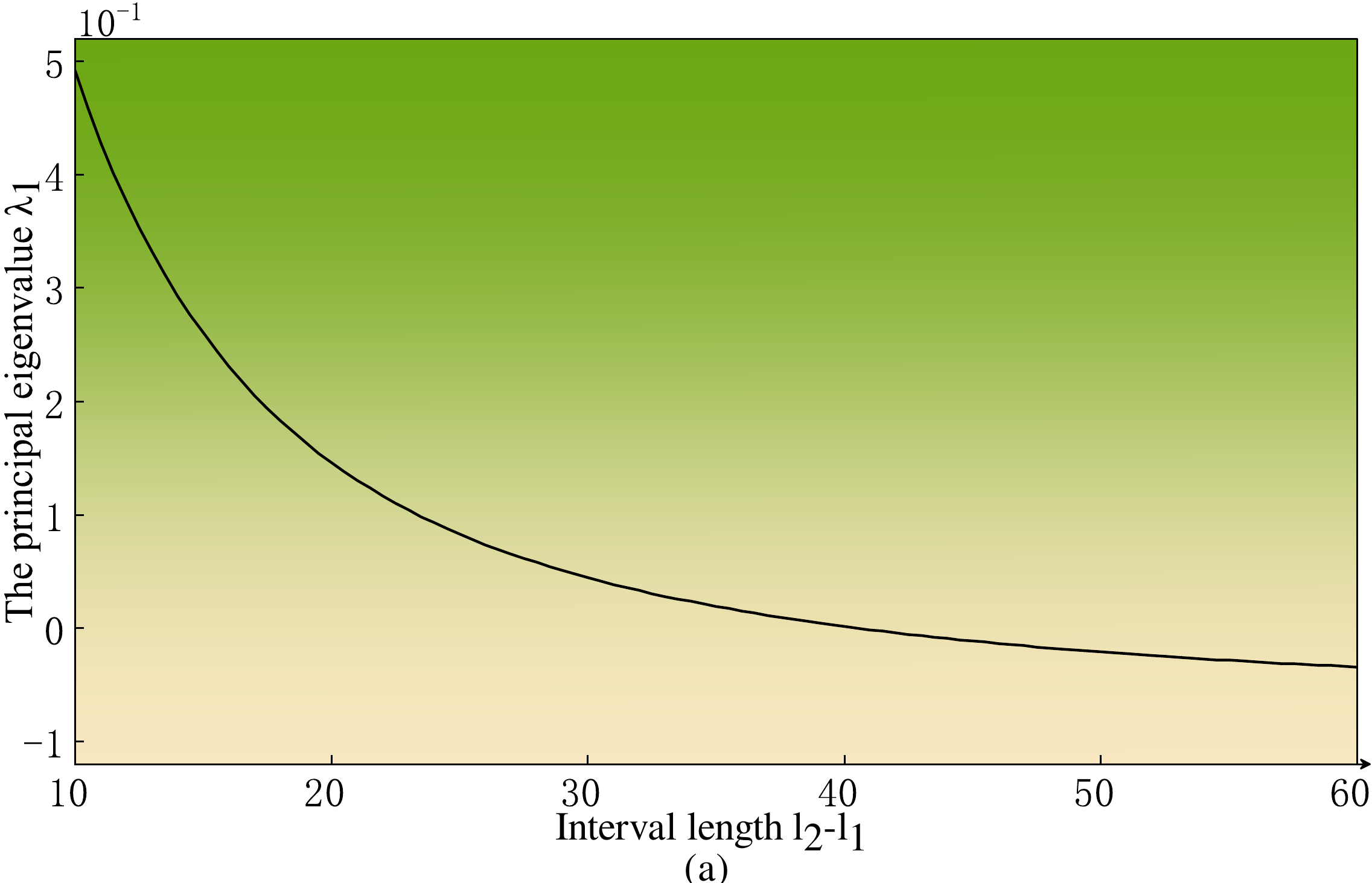}
} }\\
\subfigure{ {
\includegraphics[width=0.47\textwidth]{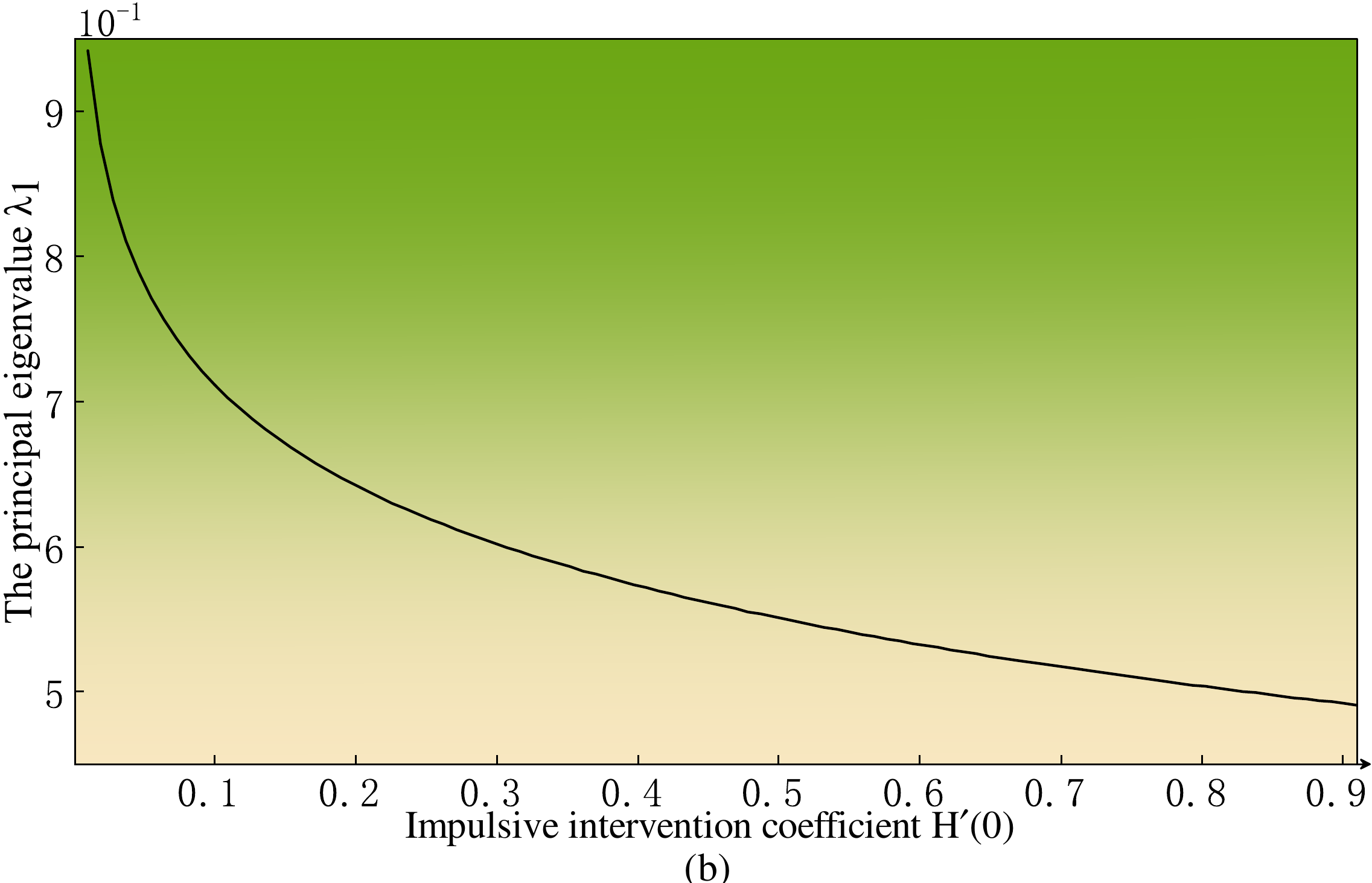}
} }
\subfigure{ {
\includegraphics[width=0.47\textwidth]{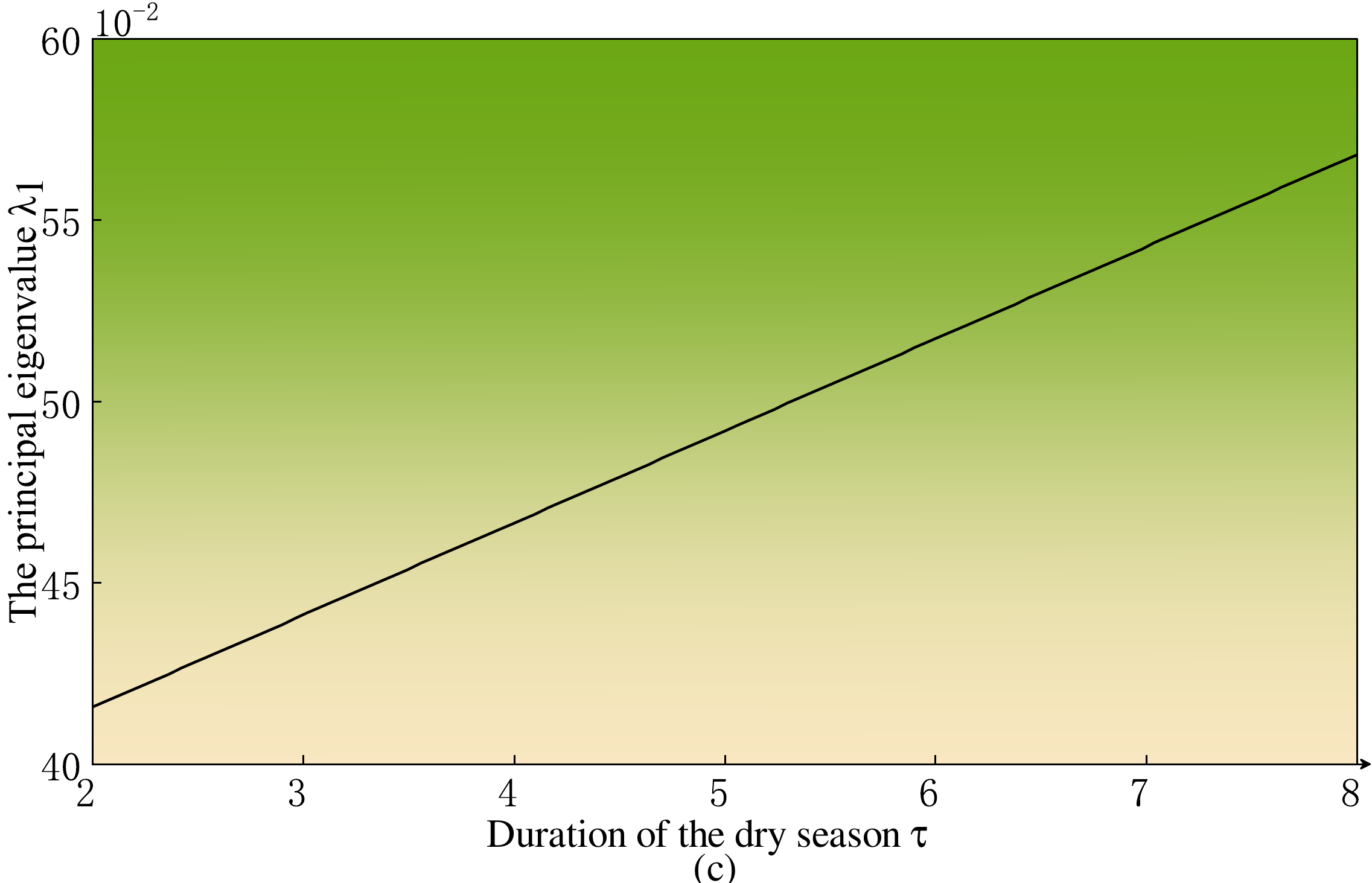}
} }
\captionsetup{justification=justified, singlelinecheck=false}
\caption{Monotonicity of the principal eigenvalue $\lambda_{1}$ with respect to $l_{2}-l_{1}$, $H'(0)$, and $\tau$, respectively. (a) decreasing of $\lambda_{1}$ with respect to $l_{2}-l_{1}$ under $H'(0)=0.9$ and $\tau=5$. (b) decreasing of $\lambda_{1}$ with respect to $H'(0)$ under $\tau=5$ and $l_{2}-l_{1}=20$. (c) increasing of $\lambda_{1}$ with respect to $\tau$ under $H'(0)=0.9$ and $l_{2}-l_{1}=20$.}
\label{B}
\end{figure}

We choose $H'(0)=0.9$ and $\tau=5$. The length of the infected region takes values from $10$ to $60$. It can be observed from \autoref{B}\textcolor[rgb]{0.00,0.00,1.00}{(a)} that the principal eigenvalue $\lambda_{1}$ is strongly monotonically decreasing with the length of the infected region $l_{2}-l_{1}$, in agreement with \autoref{lemma 3.1.3}\textcolor{blue}{(1)}.

We set $\tau=5$ and $l_{2}-l_{1}=20$, and vary the coefficient of impulsive interventions from $0.01$ to $0.91$. As shown in \autoref{B}\textcolor[rgb]{0.00,0.00,1.00}{(b)}, the principal eigenvalue $\lambda_{1}$ decreases strictly monotonically with respect to the coefficient of impulsive interventions $H'(0)$, which is consistent with the conclusion in \autoref{lemma 3.1.3}\textcolor{blue}{(2)}.

The duration of the dry season takes values from $2$ to $8$. By taking $H'(0)=0.9$ and $l_{2}-l_{1}=20$, the condition in \autoref{lemma 3.1.3}\textcolor{blue}{(3)} is satisfied. \autoref{lemma 3.1.3}\textcolor{blue}{(3)} indicates that the principal eigenvalue $\lambda_{1}$ is strictly monotonically increasing with respect to the duration of the dry season $\tau$, which is consistent with the observations from \autoref{B}\textcolor[rgb]{0.00,0.00,1.00}{(c)}.
\subsection{Effect of impulsive intervention coefficient}
This subsection graphically illustrates the impact of pulse intervention on diseases transmitted via the faecal-oral route by varying the intervention coefficient. Some of parameters are chosen as $d_{1}=0.5$, $d_{2}=0.5$, $a_{11}=0.8$, $a_{12}=1.67$, $a_{22}=0.8$, $\delta_{1}=1.5$, $\delta_{2}=1.5$, $a=1$, $m=1.7$, $\mu_{1}=6$, $\mu_{2}=8$, $\tau=6$, and $T=20$. Next, we utilize different impulsive intervention functions in order to present the distinct dynamical behaviors.
\begin{figure}[!ht]
\centering
\subfigure{ {
\includegraphics[width=0.45\textwidth]{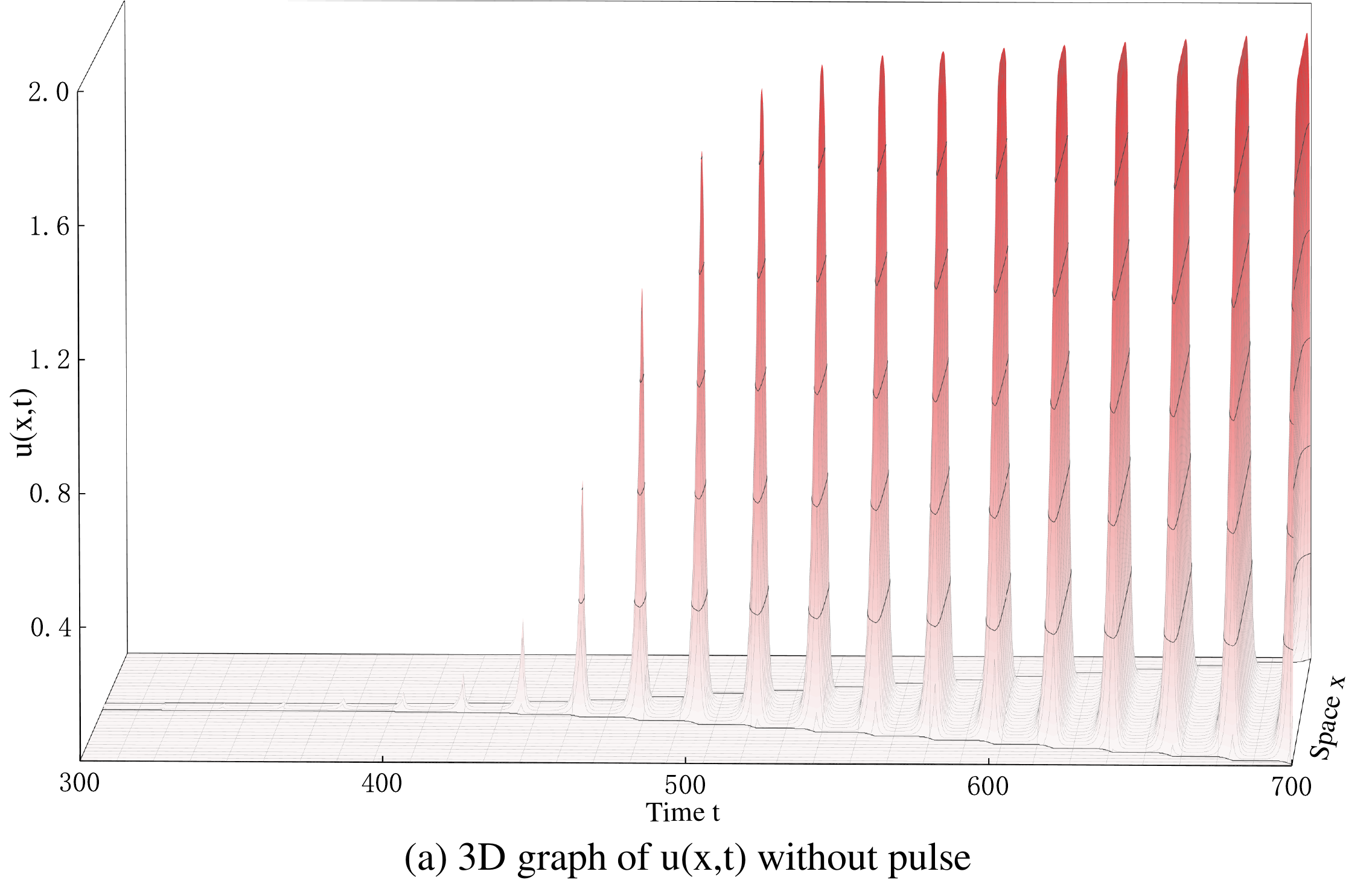}
} }
\subfigure{ {
\includegraphics[width=0.45\textwidth]{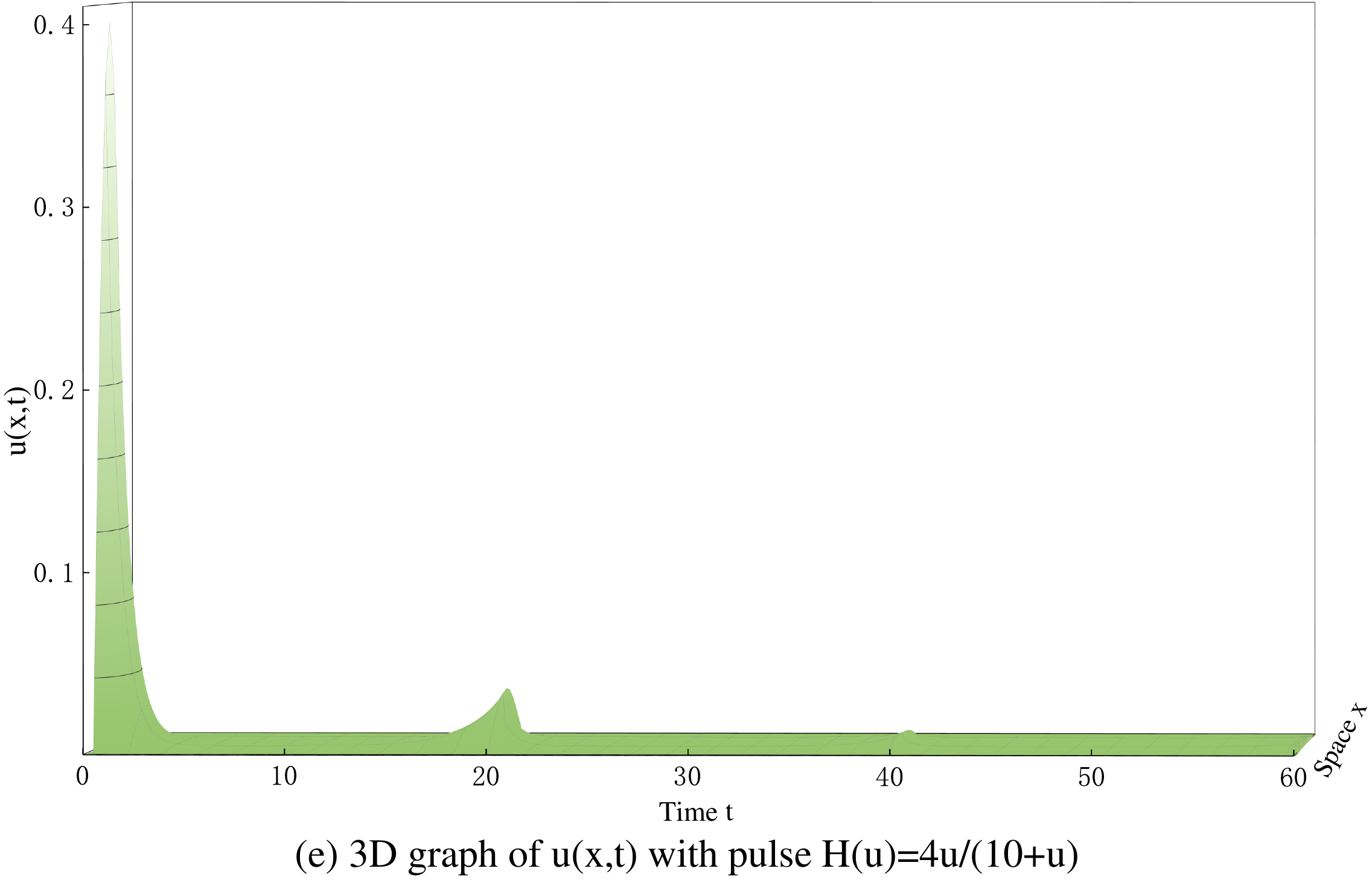}
} }
\subfigure{ {
\includegraphics[width=0.45\textwidth]{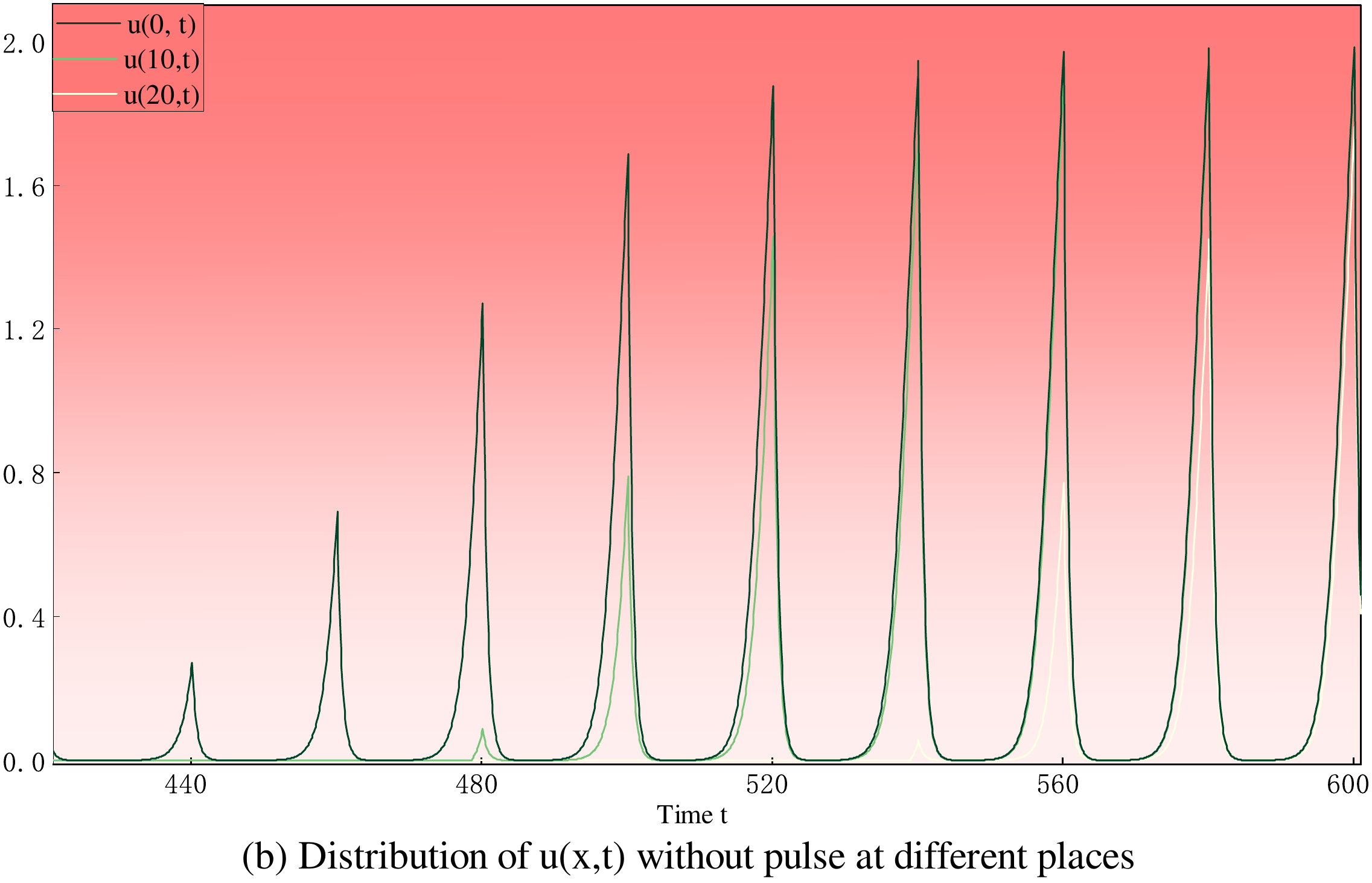}
} }
\subfigure{ {
\includegraphics[width=0.45\textwidth]{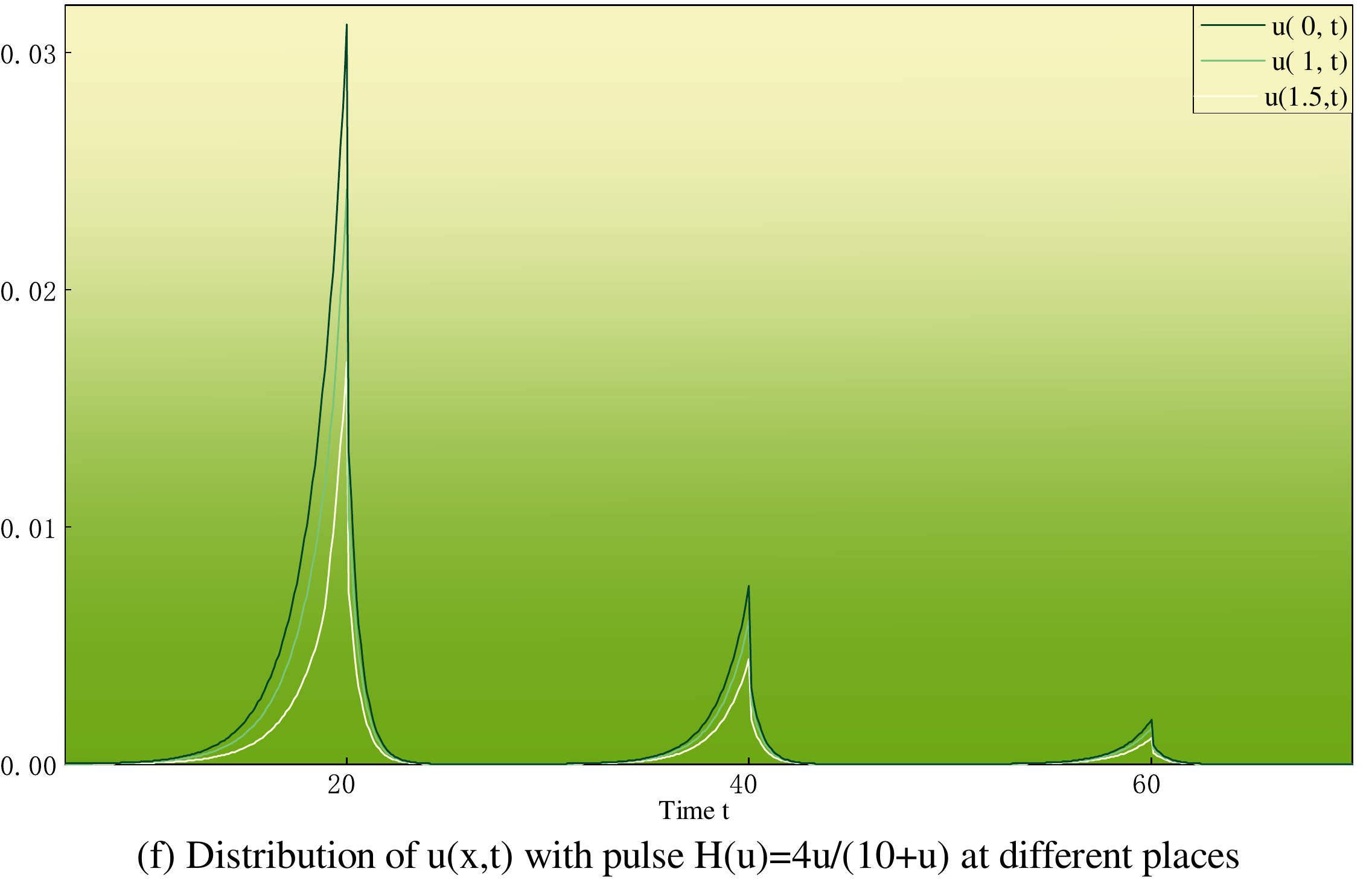}
} }
\subfigure{ {
\includegraphics[width=0.45\textwidth]{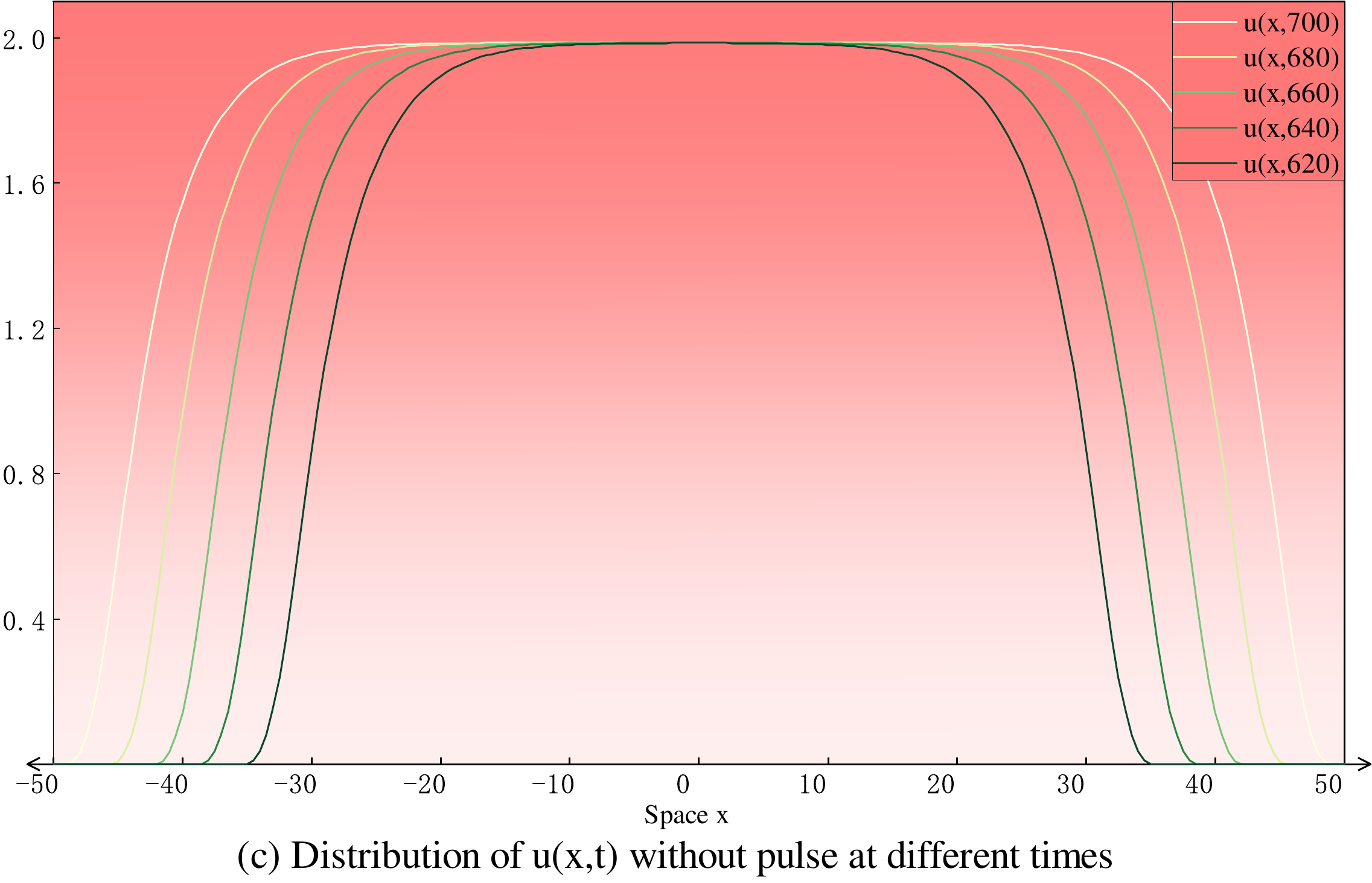}
} }
\subfigure{ {
\includegraphics[width=0.45\textwidth]{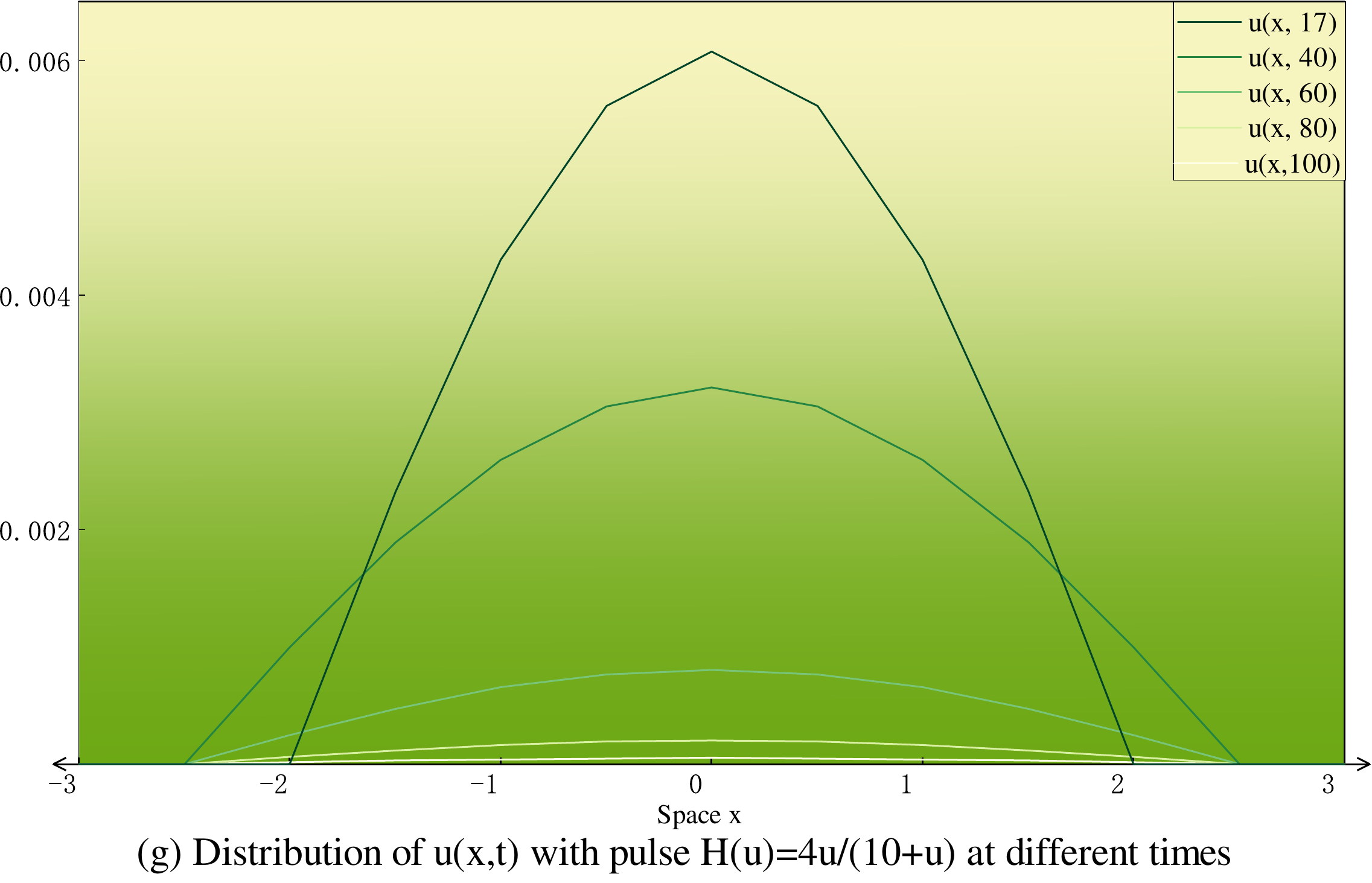}
} }
\subfigure{ {
\includegraphics[width=0.45\textwidth]{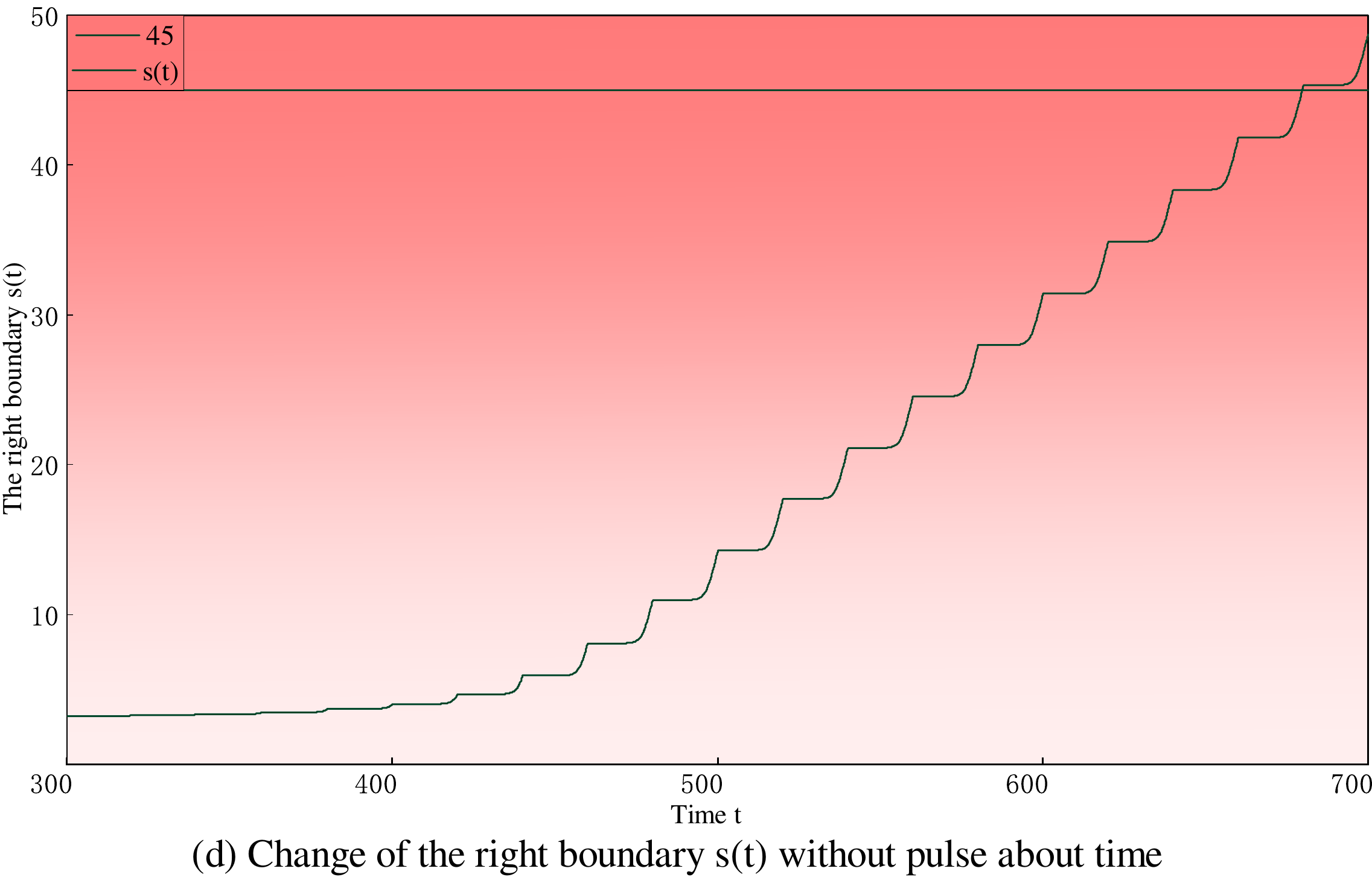}
} }
\subfigure{ {
\includegraphics[width=0.45\textwidth]{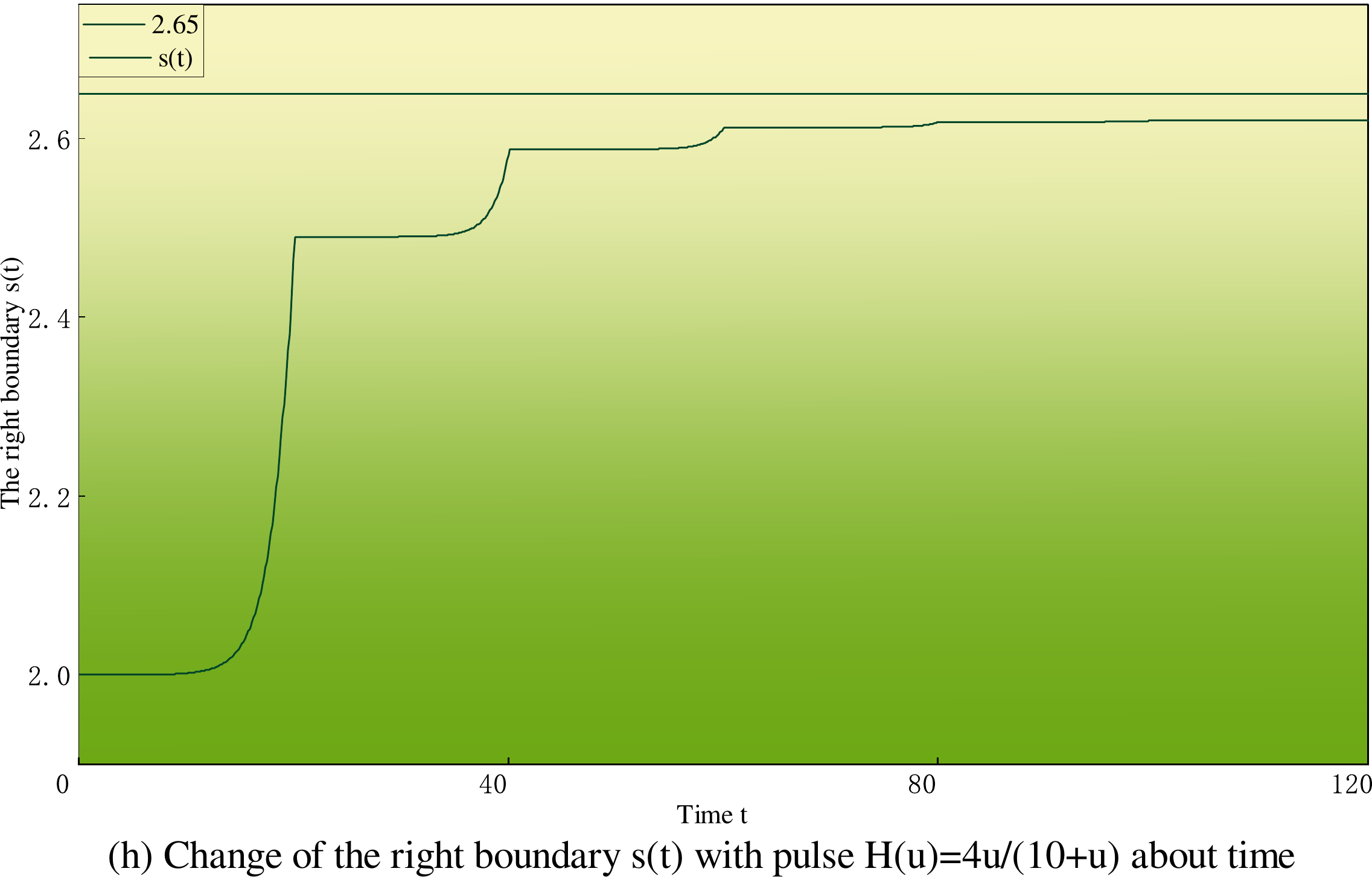}
} }
\captionsetup{justification=justified, singlelinecheck=false}
\caption{The left four images (a-d), in the absence of pulse intervention, show that the infectious agents $u$ is spreading. In contrast, the right four images (e-h), with the application of the impulsive function $H(u)=\frac{4u}{10+u}$, demonstrate that $u$ is vanishing.}
\label{C}
\end{figure}

When there is no the impulsive interventions, it can be seen from \autoref{C}\textcolor[rgb]{0.00,0.00,1.00}{(d)} that $s_{\infty}>s(700)>45$, which combined with \autoref{lemma 3.1.3}\textcolor{blue}{(1)} deduces that
\begin{equation*}
\lambda_{1}(r_{\infty}, s_{\infty})<\lambda_{1}(-45, 45)=-0.169<0.
\end{equation*}
It then follows from \autoref{lemma 4-1} that $s_{\infty}-r_{\infty}=\infty$, which combined with \autoref{lemma 4-3} deduces that
\begin{equation*}
\lim\limits_{m\rightarrow+\infty}u(x,t+mT)=W(t)
\end{equation*}
locally uniformly in $\mathcal{R}$ and uniformly in $[0,T]$, where $W(t)$ is defined in \autoref{lemma 3.1.5}. In fact, \autoref{C}\textcolor[rgb]{0.00,0.00,1.00}{(a)} shows that the spatial density of the infectious agents, denoted by $u(x,t)$, tends to a periodic steady state as time approaches infinity. \autoref{C}\textcolor[rgb]{0.00,0.00,1.00}{(b)} and \textcolor[rgb]{0.00,0.00,1.00}{(c)} demonstrate that this convergence is uniform in time and locally uniform in space, respectively. Furthermore, \autoref{C}\textcolor[rgb]{0.00,0.00,1.00}{(d)} reveals that the right boundary of the infected region, denoted by $s(t)$, continues to expand during the rainy season. These numerical results align consistently with the theoretical findings.

When the impulsive function $H(u)$ is taken as $\frac{4u}{10+u}$, one can obtain from \autoref{C}\textcolor[rgb]{0.00,0.00,1.00}{(h)} that $s_{\infty}<2.65$. Then, it follows from \autoref{lemma 3.1.3}\textcolor{blue}{(1)} that
\begin{equation*}
\lambda_{1}(r_{\infty}, s_{\infty})>\lambda_{1}(-2.65, 2.65)=0.003>0,
\end{equation*}
which combined with \autoref{theorem 4-1} deduces that vanishing happens. Actually, it is easily seen from \autoref{C}\textcolor[rgb]{0.00,0.00,1.00}{(e-g)} that the spatial density of the infectious agents converges to zero uniformly in apace $x$ and periodically in time $t$, which is consistent with the theoretical conclusion.

In conclusion, implementing impulsive interventions with appropriate intensity can drive a spreading diseases to vanishing.
\subsection{Effect of duration of the dry season}
This subsection presents the effect of the duration of the dry season on the development of the diseases. We choose some parameters as $d_{1}=0.5$, $d_{2}=0.5$, $a_{11}=0.8$, $a_{12}=1.7$, $a_{22}=0.8$, $\delta_{1}=0.9$, $\delta_{2}=0.9$, $a=1$, $m=1.7$, $\mu_{1}=6$, $\mu_{2}=8$, $H(u)=u$, and $T=10$. The duration of the dry season, denoted by $\tau$, is chosen in each of the following simulations.
\begin{figure}[!ht]
\centering
\subfigure{ {
\includegraphics[width=0.45\textwidth]{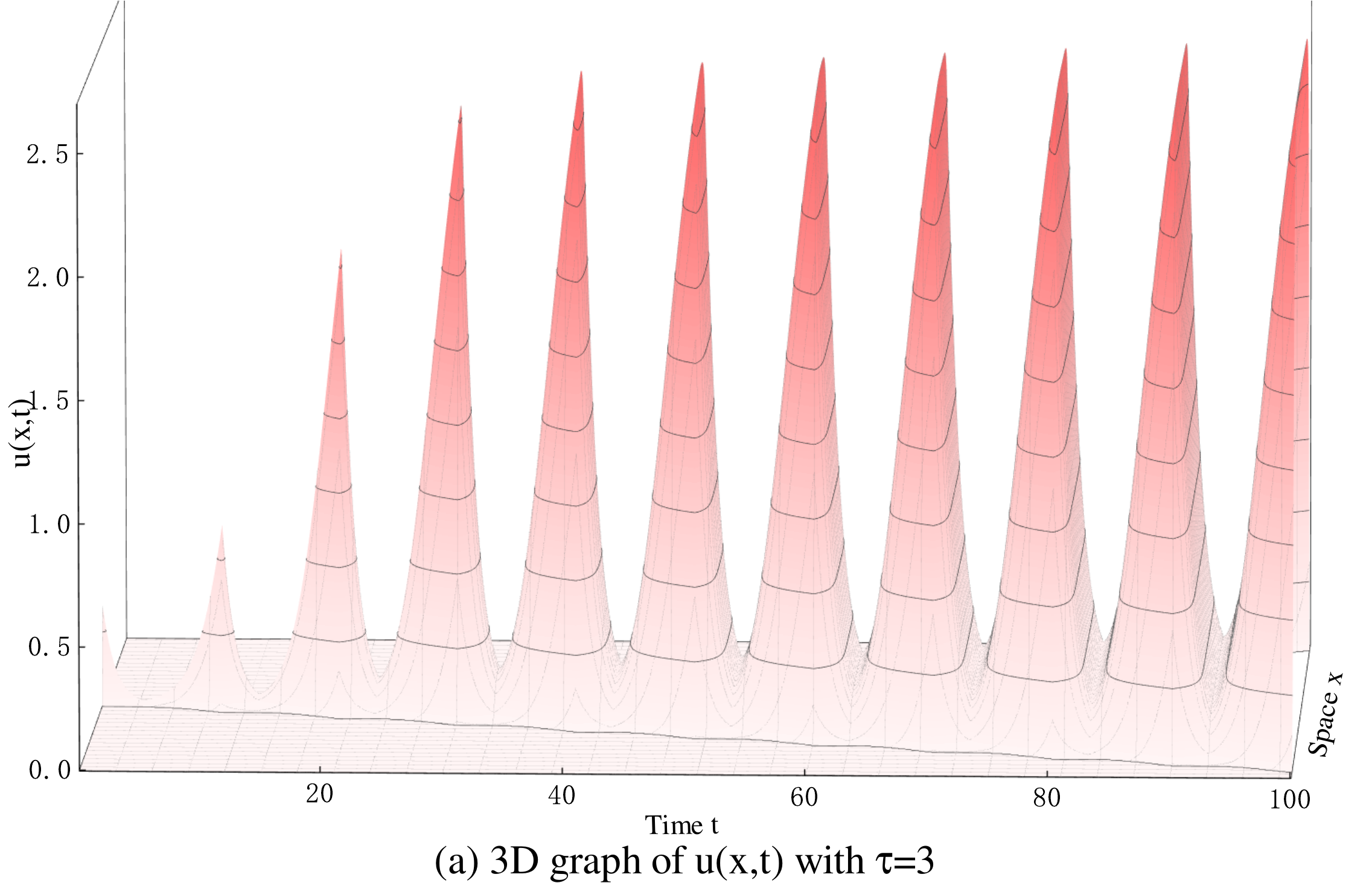}
} }
\subfigure{ {
\includegraphics[width=0.45\textwidth]{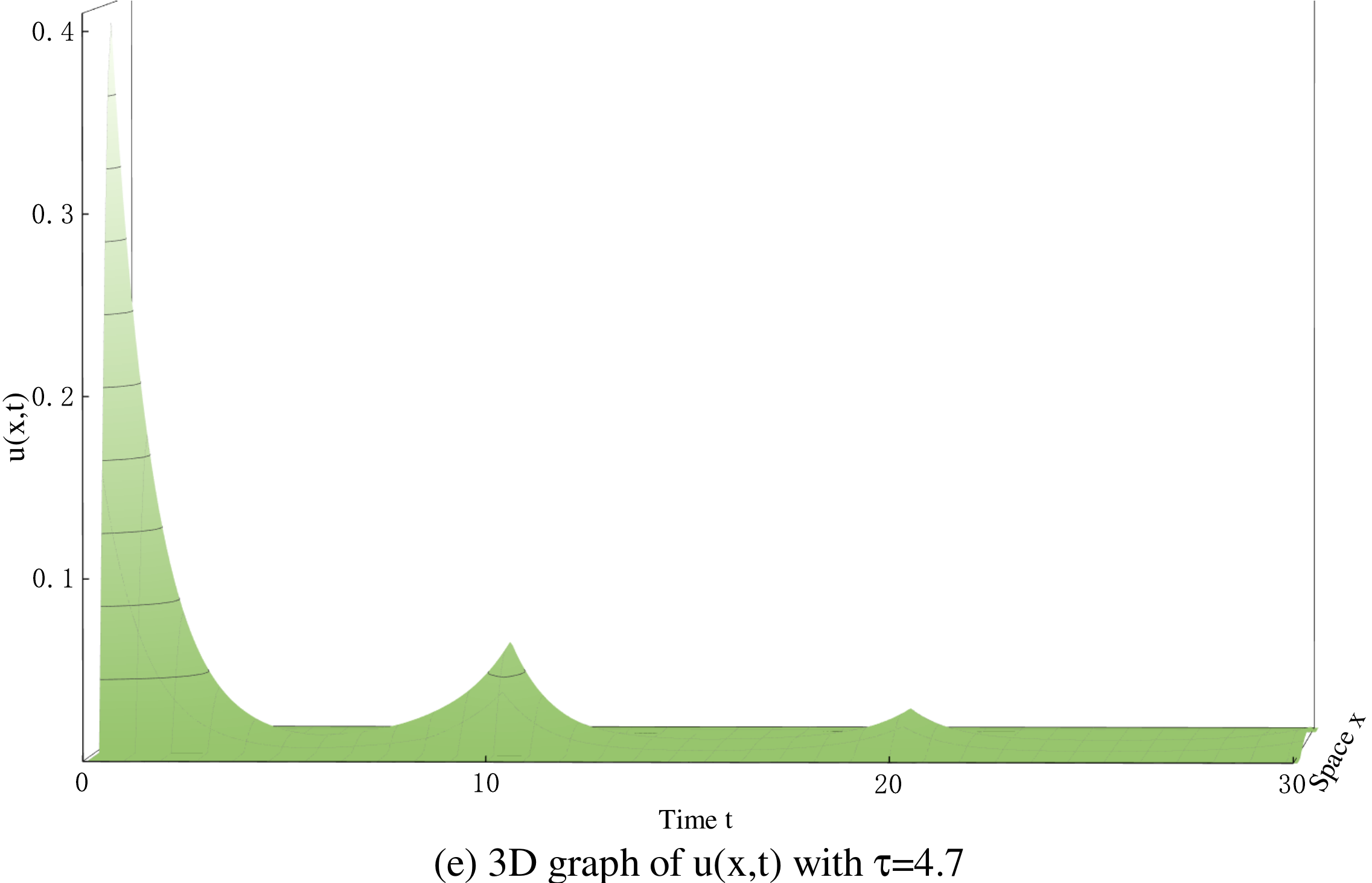}
} }
\subfigure{ {
\includegraphics[width=0.45\textwidth]{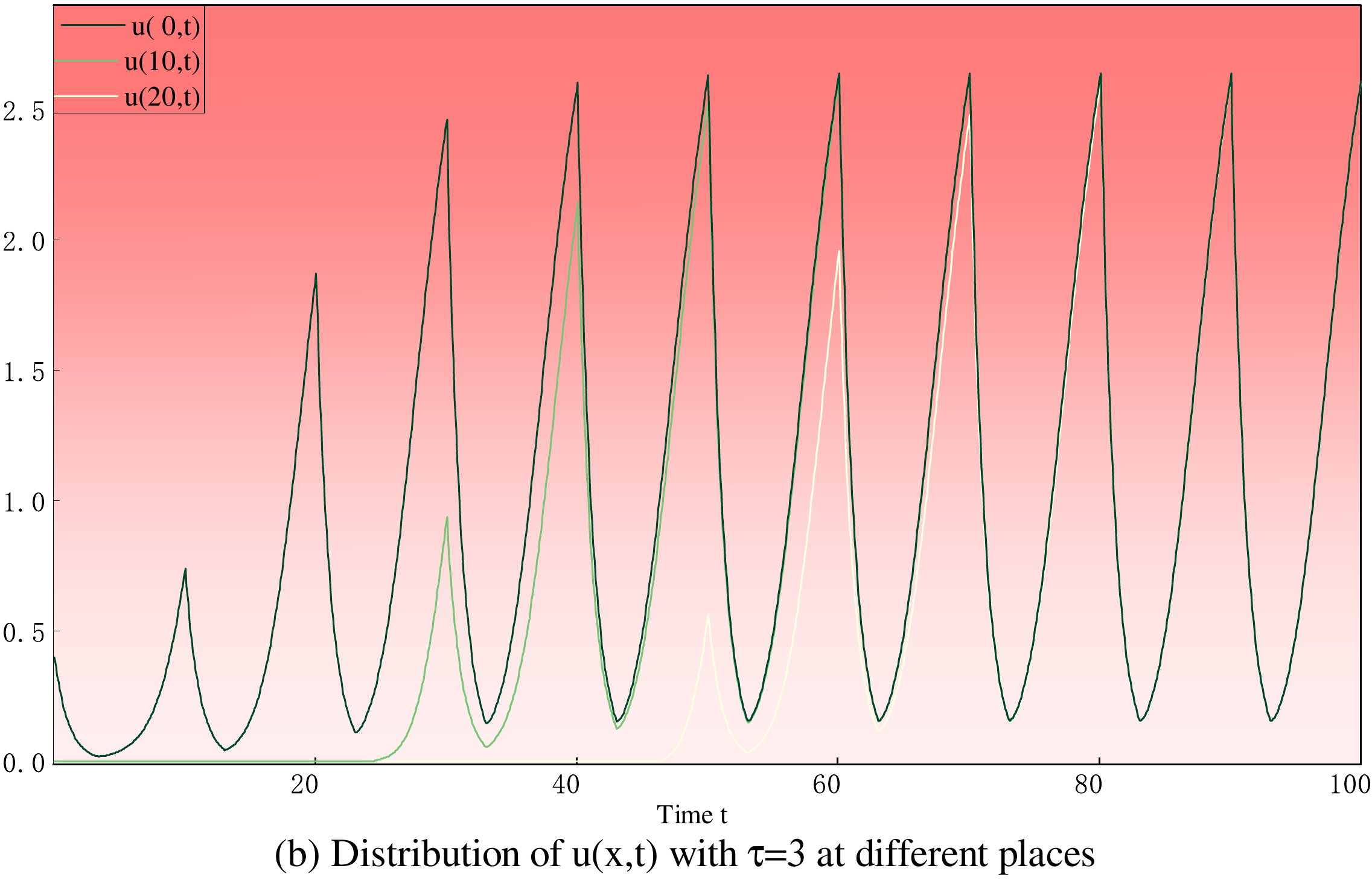}
} }
\subfigure{ {
\includegraphics[width=0.45\textwidth]{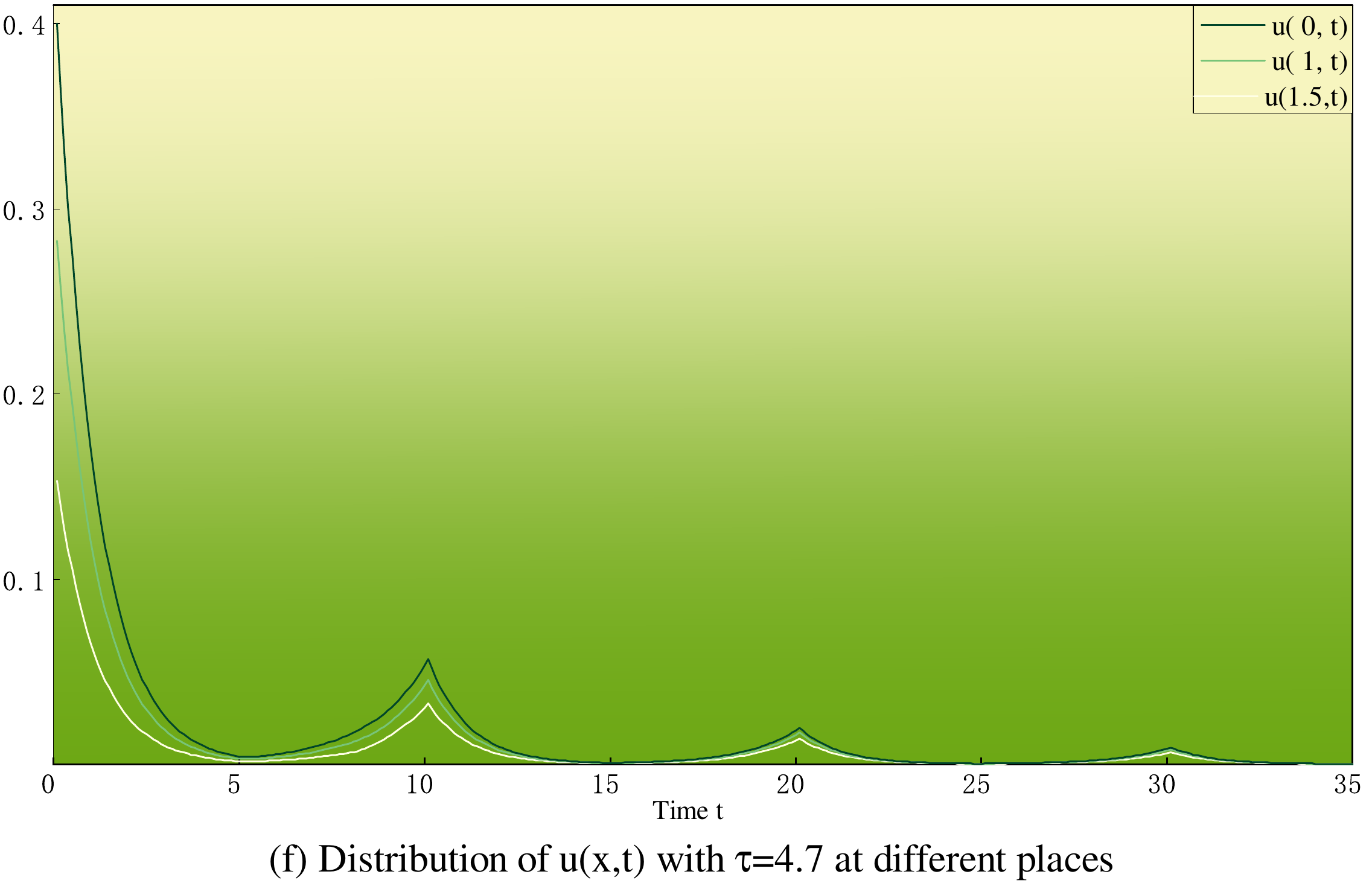}
} }
\subfigure{ {
\includegraphics[width=0.45\textwidth]{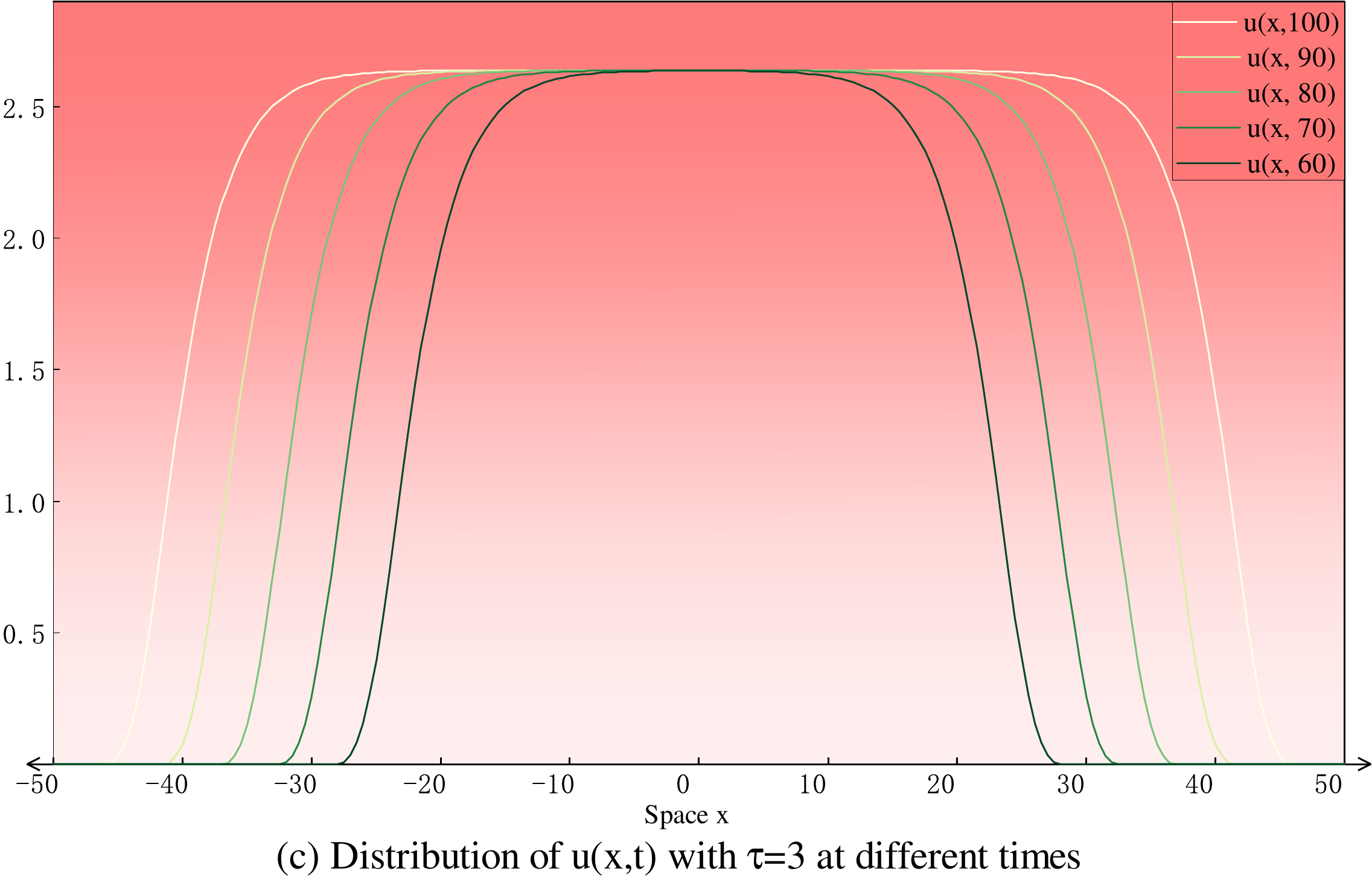}
} }
\subfigure{ {
\includegraphics[width=0.45\textwidth]{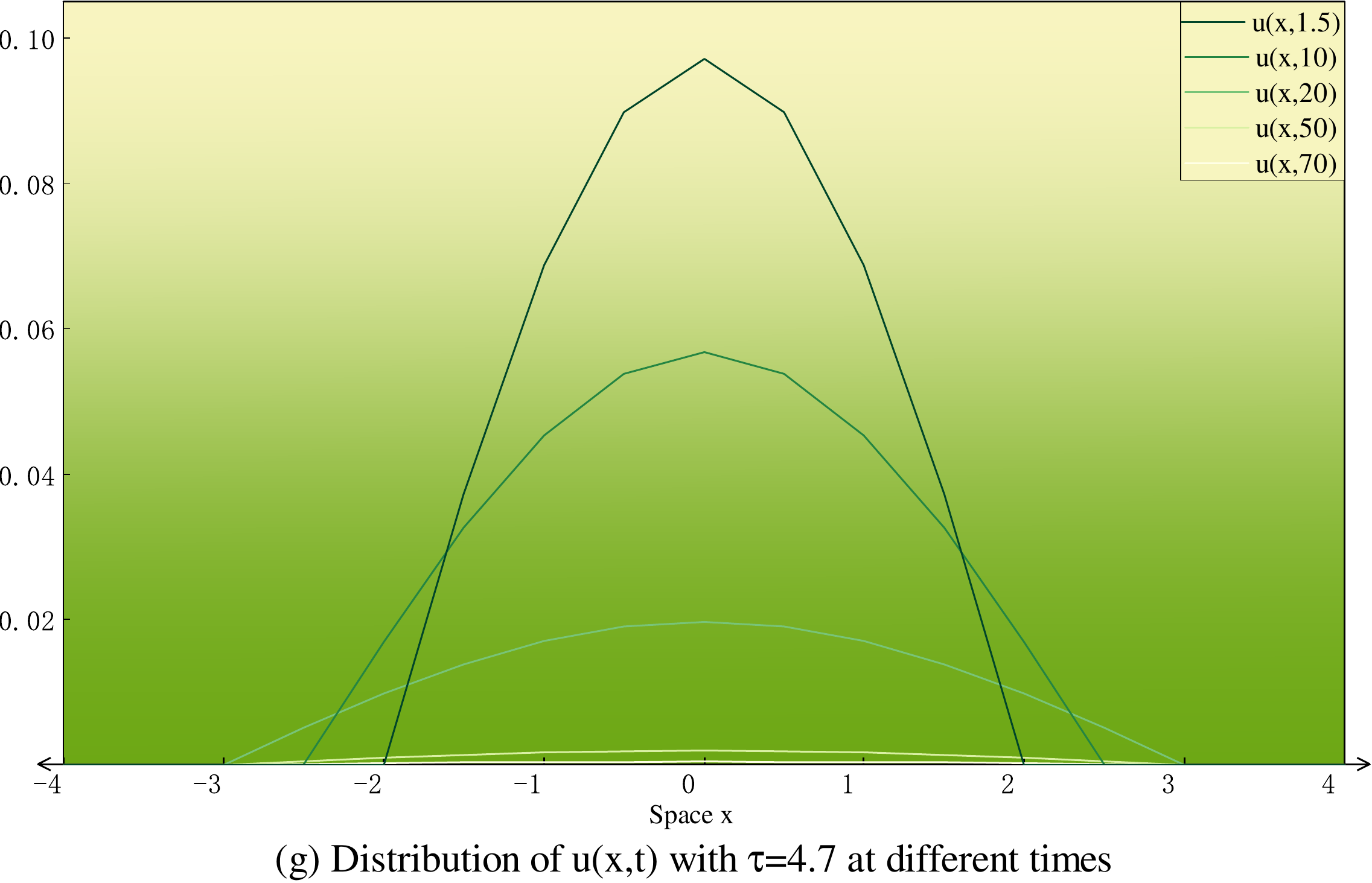}
} }
\subfigure{ {
\includegraphics[width=0.45\textwidth]{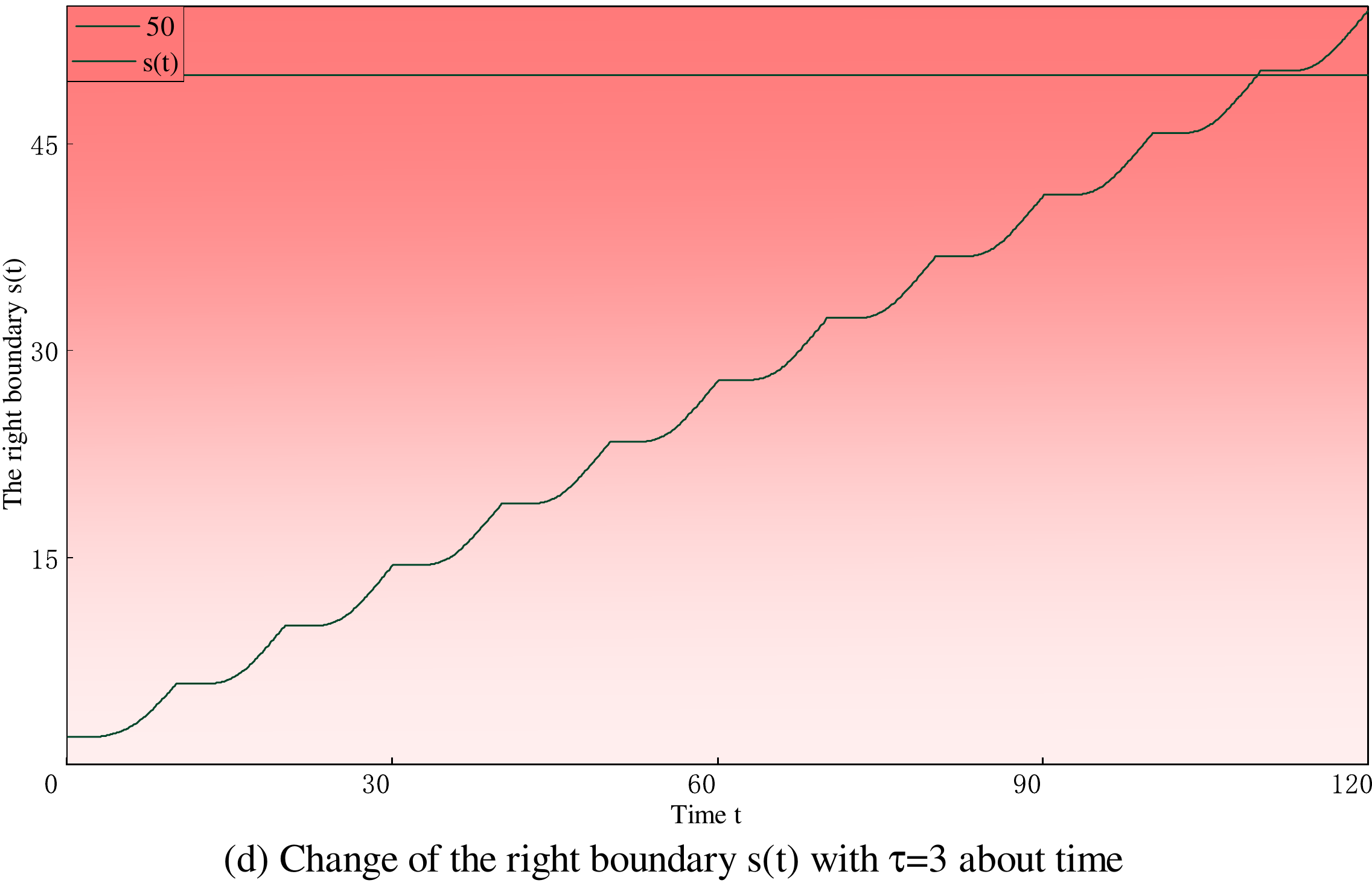}
} }
\subfigure{ {
\includegraphics[width=0.45\textwidth]{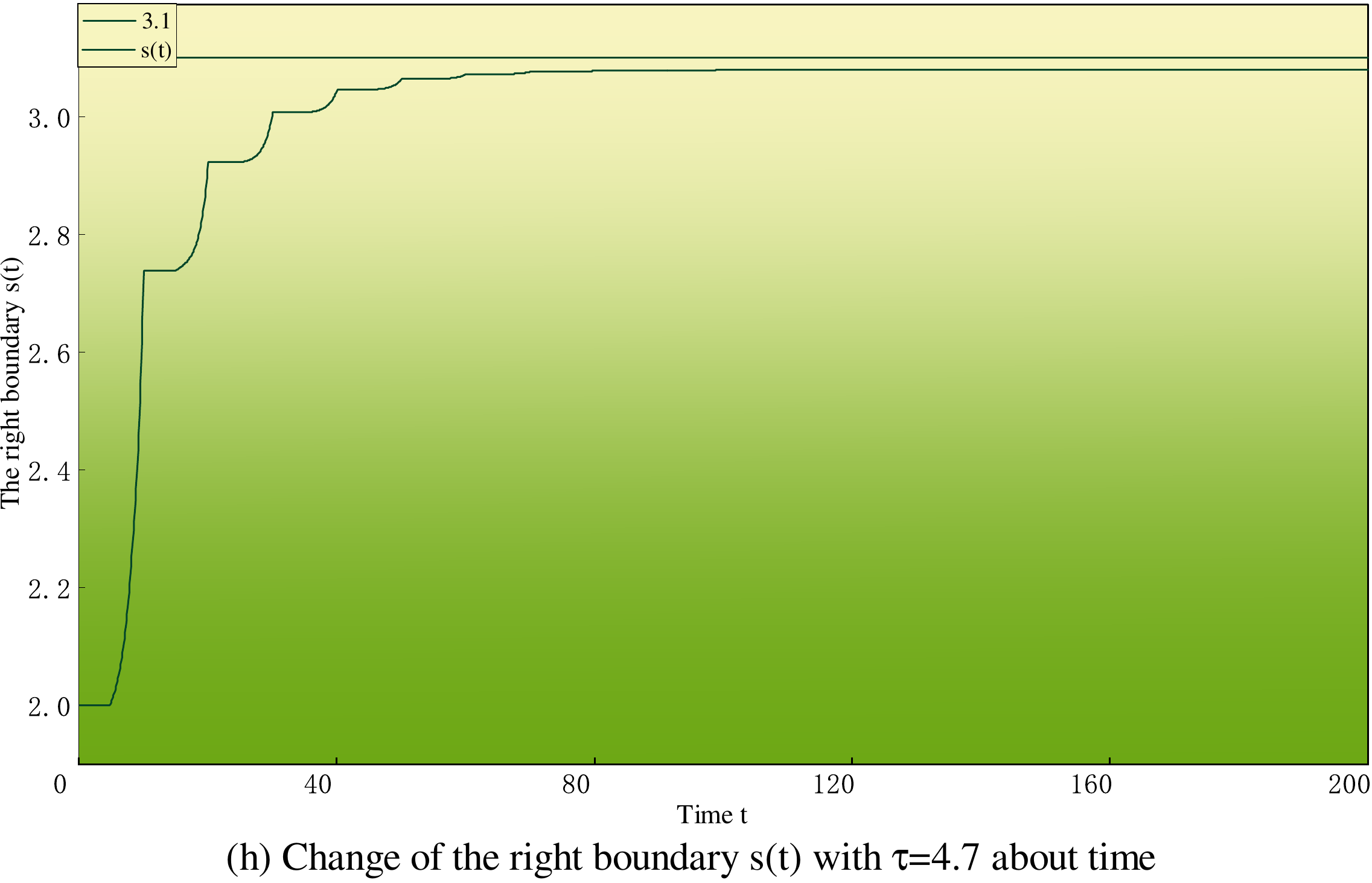}
} }
\captionsetup{justification=justified, singlelinecheck=false}
\caption{The left four images (a-d), with short dry season duration ($\tau=3$), show that the infectious agents $u$ is spreading. In contrast, the right four images (e-h), with long dry season duration ($\tau=4.7$), demonstrate that $u$ is vanishing.}
\label{D}
\end{figure}

We first consider a short dry season duration, set to $3$. By observing \autoref{D}\textcolor[rgb]{0.00,0.00,1.00}{(d)}, we have that $s_{\infty}>s(700)>50$, which combined with \autoref{lemma 3.1.3}\textcolor{blue}{(1)} yields that
\begin{equation*}
\lambda_{1}(r_{\infty}, s_{\infty})<\lambda_{1}(-50, 50)=-0.360<0.
\end{equation*}
Then, it follows from \autoref{lemma 4-1} that $s_{\infty}-r_{\infty}=\infty$, which combined with \autoref{lemma 4-3} deduces that
\begin{equation*}
\lim\limits_{m\rightarrow+\infty}u(x,t+mT)=W(t)
\end{equation*}
locally uniformly in $\mathcal{R}$ and uniformly in $[0,T]$, where $W(t)$ is defined in \autoref{lemma 3.1.5}. This is consistent with the observations from \autoref{D}\textcolor[rgb]{0.00,0.00,1.00}{(a-c)}.

Subsequently, we consider a long dry season duration, set to $4.7$. It follows from \autoref{D}\textcolor[rgb]{0.00,0.00,1.00}{(h)} that $s_{\infty}<3.1$. In virtue of the monotonicity of the principal eigenvalue with respect to the length of the infected interval, we have that
\begin{equation*}
\lambda_{1}(r_{\infty}, s_{\infty})>\lambda_{1}(-3.1, 3.1)=0.040>0.
\end{equation*}
This combined with \autoref{theorem 4-1} yields that vanishing happens, in agreement with the numerical simulation results of \autoref{D}\textcolor[rgb]{0.00,0.00,1.00}{(e-g)}.

To summarize, an increase in the dry season duration is positively correlated with the effectiveness of controlling faecal-oral diseases transmission.
\section{Conclusion and future work}\label{section-6}
In this paper, we develop a two-season faecal-oral model with impulsive intervention in a moving infected environment by simultaneously considering human  interventive behaviour and temporal variation in rainfall.
In the model, the fronts of the infection region during the dry season and rainy season are represented by fixed and moving boundaries, respectively. Our developed model extends the models in \cite{Ahn-Baek-Lin-AMM,Zhao-Li-DCDSB,Wang-Du-DCDSB,Zhou-Lin-Pedersen-ARXIV,7,Capasso-Maddalena-3}.

The simultaneous introduction of impulsive intervention and seasonal switching enables the model to capture more realistic phenomena, generating certain difficulties in mathematical analysis. To address these challenges, some new mathematical analyses include:
\begin{itemize}
\item{selecting appropriate function spaces and applying a strong version of the Krein-Rutman theorem to prove the existence, uniqueness, and algebraic simplicity of the associated principal eigenvalue;}
\item{establishing sharper estimates for the corresponding principal eigenfunction;}
\item{adapting the construction of comparison functions as upper and lower solutions.}
\end{itemize}
The techniques employed here are also applicable to other cooperative models with impulsive intervention and seasonal switching.

We obtain from parabolic regularity theory and the contraction mapping principle that the strong positive classical solution of the model exists and is unique, except at the points of impulse and seasonal switching. Our main findings are as follows:
\begin{itemize}
\item{a vanishing-spreading dichotomy is provided in \autoref{theorem 4-5}; }
\item{a sharp criteria governing this dichotomy is presented in \autoref{theorem 4-6}.}
\end{itemize}
Both theoretical analysis and numerical simulations indicate that the intensity of impulsive intervention and the duration of the dry season are positively correlated with the efficacy of disease control, which is consistent with observations within the real world.

When spreading occurs, estimating the spreading speed of the model constitutes an important and meaningful research topic, which we leave for future work.

\end{document}